\documentclass[12pt]{amsart}
\usepackage[margin=0.6in]{geometry}

\usepackage{setspace}
\usepackage[english]{babel}
\usepackage[hyperindex, linktocpage]{hyperref}
\usepackage[backrefs, msc-links]{amsrefs}
\usepackage{amsmath,amssymb,amsfonts,amsthm,enumerate}
\usepackage{url}
\usepackage{graphicx}
\usepackage{xcolor}
\usepackage{wrapfig}

\usepackage{wesa}
\usepackage[T1]{fontenc}

\usepackage[cal=boondox]{mathalfa}

\definecolor{darkblue}{RGB}{0,0,170}
\definecolor{brickred}{RGB}{200,0,0}
\hypersetup{colorlinks=true, linkcolor=darkblue, citecolor=blue} 

\numberwithin{equation}{section}

\newtheorem{difference}{Difference}[section]

\newcommand{\R}{\mathbb{R}}
\newcommand{\Z}{\mathbb{Z}}
\newcommand{\N}{\mathbb{N}}
\newcommand{\C}{\mathbb{C}}
\newcommand{\SURD}{\displaystyle{\surd}}

\usepackage{yfonts}
\newcommand{\LLL}{\textswab{L}}
\newcommand{\ddd}{\textswab{d}}

\def\Xint#1{\mathchoice
{\XXint\displaystyle\textstyle{#1}}%
{\XXint\textstyle\scriptstyle{#1}}%
{\XXint\scriptstyle\scriptscriptstyle{#1}}%
{\XXint\scriptscriptstyle\scriptscriptstyle{#1}}%
\!\int}
\def\XXint#1#2#3{{\setbox0=\hbox{$#1{#2#3}{\int}$ }
\vcenter{\hbox{$#2#3$ }}\kern-.6\wd0}}

\def\pashint{\Xint{\overline{\quad}}}

\renewcommand{\varpi}{\omega}

\renewcommand{\leq}{\leqslant}
\renewcommand{\le}{\leqslant}
\renewcommand{\geq}{\geqslant}
\renewcommand{\ge}{\geqslant}

\renewcommand{\limsup}{\varlimsup}

\renewcommand{\epsilon}{\varepsilon}
\newcommand{\eps}{\varepsilon}
\newcommand{\e}{\varepsilon}
\newcommand{\const}{\,{\rm const}\,}
\newcommand{\PV}{\,{\rm P.V.}\,}

\title[Getting acquainted with the fractional Laplacian]{Getting acquainted with the fractional Laplacian}\thanks{Supported by INdAM}

\author{Nicola Abatangelo}

\address{{\em Nicola Abatangelo:} 
D\'epartement de math\'ematique, Universit\'e Libre de Bruxelles CP 214,
Boulevard du Triomphe, 1050 Ixelles,
Belgium.}

\email{nicola.abatangelo@ulb.ac.be}

\author{Enrico Valdinoci}

\address{{\em Enrico Valdinoci:} 
Department of Mathematics and Statistics,
University of Western Australia, 35 Stirling Highway, Crawley WA 6009, Australia.}

\email{enrico.valdinoci@uwa.edu.au}

\keywords{Fractional calculus, functional analysis, applications.}
\subjclass[2010]{35R11, 34A08, 60G22.}

\begin{document}

\begin{abstract}
These are the handouts of an undergraduate minicourse at the Universit\`a di Bari
(see Figure~\ref{FOTO}),
in the context of the 2017
INdAM Intensive Period
``Contemporary Research in elliptic PDEs and related topics''.
Without any intention to serve as a throughout epitome to the subject,
we hope that these notes
can be of some help 
for a very initial introduction to a fascinating field of classical and
modern research.
\end{abstract}

\maketitle

\tableofcontents

\section*{Acknowledgements}
It is a great pleasure to thank the Universit\`a degli Studi di Bari
for its very warm hospitality and
the Istituto Nazionale di Alta Matematica for the strong financial and administrative
support which made this minicourse possible. And of course special thanks go to
all the participants, for their patience in attending the course, their competence,
empathy and
 contagious enthusiasm. 

\section{The Laplace operator}

The operator mostly studied in partial differential equations is likely the so-called
Laplacian, given by
\begin{equation} \label{DEL:CLA}
-\Delta u(x):=-\sum_{j=1}^n
\frac{\partial^2 u}{\partial x_j^2} (x)
=\lim_{r\searrow0} \frac{\const}{r^{n+2}} \,\int_{B_r(x)} \big( u(x)-u(y)\big)\,dy
{ 
=-\const\int_{\partial B_1} D^2u(x)\,\theta\cdot\theta\, d\theta}
\end{equation}
Of course, one may wonder
why mathematicians have a strong preference for
such kind of operators -- say, why not studying
$$\frac{\partial^7 u}{\partial x_1^7} (x)-
\frac{\partial^8 u}{\partial x_2^8} (x)+\frac{\partial^9 u}{\partial x_3^9} (x) 
-\frac{\partial^{10} u}{\partial x_1 \partial x_2^4 \partial x_3^5} (x)\quad\;
?$$
Since historical traditions, scientific legacies or
impositions from above by education systems
would not be enough
to justify such a strong interest in only one operator (plus all its modifications),
it may be worth to point out a simple geometric property enjoyed by 
the Laplacian
(and not by many other operators).
Namely,
equation~\eqref{DEL:CLA} somehow reveals that the fact that a function is harmonic
(i.e., that its Laplace operator
vanishes in some region) is deeply related
to the action of ``comparing with the surrounding values and reverting
to the averaged values in the neighborhood''.

\begin{figure}
    \centering
    \includegraphics[width=11cm]{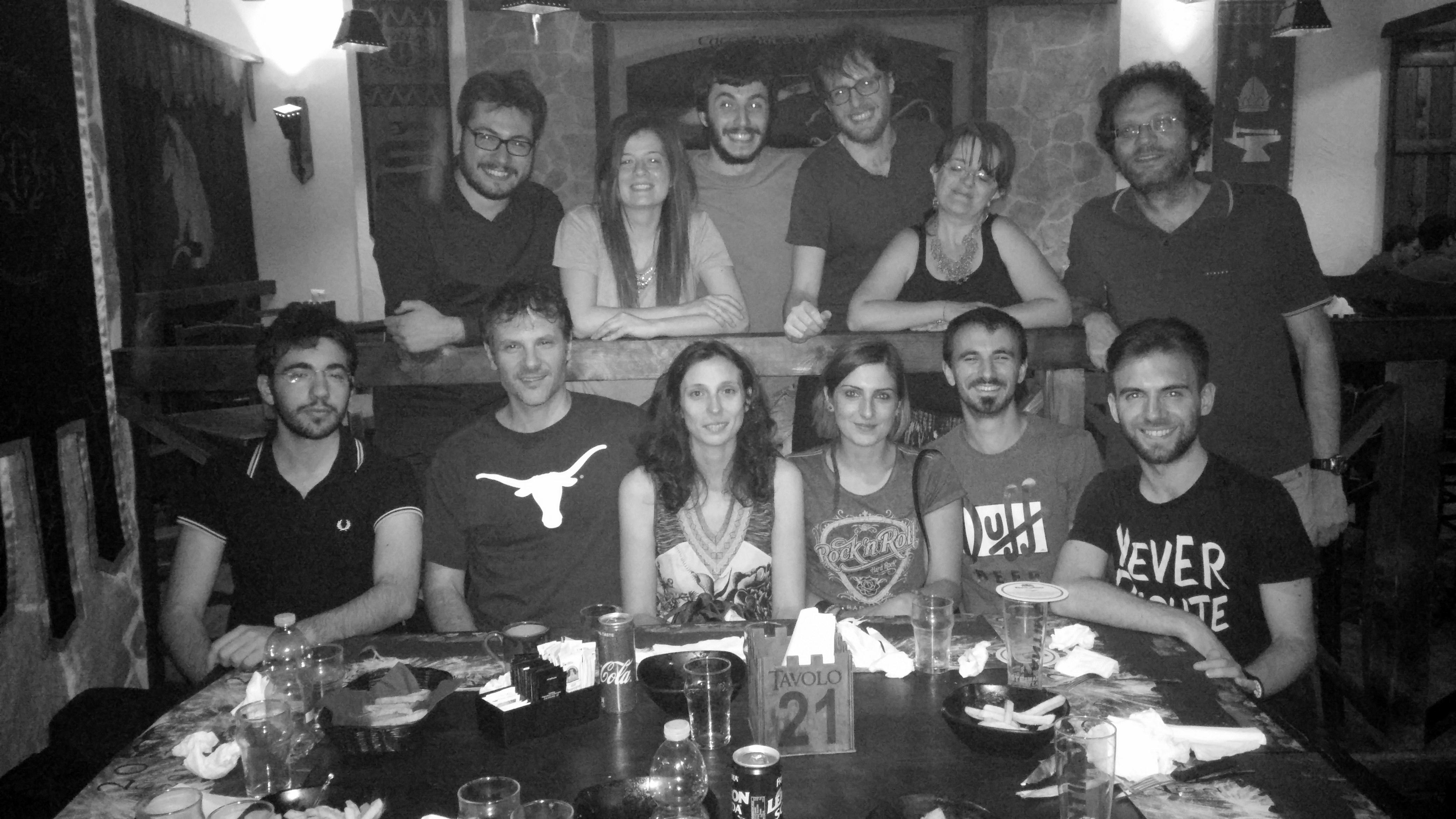}
    \caption{\it {{Working hard (and profitably) in Bari.}}}
    \label{FOTO}
\end{figure}

To wit, the idea behind the integral representation of the Laplacian in
formula~\eqref{DEL:CLA}
is that the Laplacian tries to model an ``elastic'' reaction:
the vanishing of such operator should try to ``revert the value
of a function at some point to the values nearby'', or,
in other words,
from a ``political'' perspective, the Laplacian is a very ``democratic''
operator, which aims at levelling out differences
in order to make things as uniform as possible.
In mathematical terms, one looks at the difference between
the values of a given function~$u$ and its average
in a small ball of radius~$r$, namely
$$ \ddd_r(x):= u(x)\,-\,\pashint_{B_r(x)} u(y)\,dy=
\pashint_{B_r(x)} \big( u(x)-u(y)\big)\,dy.$$
In the smooth setting,
a second order Taylor expansion of~$u$
and a cancellation in the integral due to odd symmetry show
that~$\ddd_r$ is quadratic in~$r$, hence, in order
to detect the ``elastic'', or
``democratic'', effect of the model at small scale,
one has to divide by~$r^2$ and take the limit as~$r\searrow0$.
This is exactly the procedure that we followed in
formula~\eqref{DEL:CLA}.

Other classical approaches to integral 
representations of elliptic operators come in view
of potential theory and inversion operators, see
e.g.~\cite{MR0030102}.\medskip

This tendency to revert to the surrounding mean suggests that harmonic equations,
or in general equations driven by operators ``similar to the Laplacian'', possess
some kind of rigidity or regularity properties
that prevents the solutions to oscillate too much
(of course, detecting and establishing these
properties is a marvelous, and technically extremely demanding,
success of modern mathematics, and
we do not indulge in this set of notes
on this topic of great beauty and outmost importance,
and we refer, e.g. to the classical books~\cites{MR1625845,
MR1814364,
MR1962933,
MR2777537}).\medskip

Interestingly, the Laplacian operator, in the perspective of~\eqref{DEL:CLA},
is the infinitesimal limit of integral operators. In the forthcoming sections,
we will discuss some other integral operators, which recover the Laplacian
in an appropriate limit, and which share the same property of averaging the values
of the function. Differently from what happens in~\eqref{DEL:CLA},
such averaging procedure will not be necessarily confined to a small neighborhood
of a given point, but will rather tend to comprise all the possible values
of a certain function, by possibly ``weighting more'' the close-by points
and ``less'' the contributions coming from far.

\section{Some fractional operators}\label{SOME}

We describe here the basics of some different fractional\footnote{The notion (or, better to say, several possible notions)
of fractional derivatives attracted the attention of many
distinguished mathematicians, such as
Leibniz, Bernoulli, Euler, Fourier, Abel, Liouville,
Riemann, Hadamard and Riesz, among the others.
A very interesting historical outline is given in pages~xxvii--xxxvi
of~\cite{MR1347689}.}
operators. The fractional exponent will be denoted
by~$s\in(0,1)$. For more exhaustive discussions and comparisons
see e.g.~\cites{MR0290095, MR0350027, MR1347689, MR2676137, MR2944369, MR3233760, MR3246044, MR3469920, MR3613319}.
For simplicity, we do not treat here the case of fractional
operators of order higher than~$1$ (see e.g.~\cites{GD, AJS1, AJS2, AJS3}).

\subsection{The fractional Laplacian}\label{DE:001}

A very popular nonlocal operator is given by the fractional Laplacian
\begin{equation} \label{FRAC LAP}
(-\Delta)^s u(x) := \PV \int_{\R^n} \frac{u(x)-u(y)}{|x-y|^{n+2s}}\,dy.\end{equation}
Here above, the notation ``$\PV$'' stands for ``in the Principal Value sense'', that is
$$ (-\Delta)^s u(x) := \lim_{\e\searrow0} \int_{\R^n\setminus B_\e(x)}
\frac{u(x)-u(y)}{|x-y|^{n+2s}}\,dy.$$
The definition in~\eqref{FRAC LAP} differs from others available in the literature since
a normalizing factor has been omitted for the sake of simplicity: this multiplicative
constant is only important in the limits as~$s\nearrow1$ and~$s\searrow0$, but plays no essential role
for a fixed fractional parameter~$s\in(0,1)$.

The operator in~\eqref{FRAC LAP} can be also conveniently written in the form
\begin{equation} \label{FRAC LAP:2}
-(-\Delta)^s u(x)= \frac12 \int_{\R^n} \frac{u(x+y)+u(x-y)-2u(x)}{|y|^{n+2s}}\,dy.\end{equation}
The expression in~\eqref{FRAC LAP:2} reveals that the fractional Laplacian
is a sort of second order difference operator, weighted by a measure supported in the
whole of~$\R^n$ and with a polynomial decay, namely
\begin{equation}\label{2PO}
\begin{split}&
-(-\Delta)^s u(x)= \frac12 \int_{\R^n} \delta_u(x,y)\,d\mu(y),\\
{\mbox{ where }}\quad&
\delta_u(x,y):=u(x+y)+u(x-y)-2u(x)\quad
{\mbox{ and }}\quad
d\mu(y):=\frac{dy}{|y|^{n+2s}}.
\end{split}
\end{equation}
Of course, one can give a pointwise meaning of~\eqref{FRAC LAP}
and~\eqref{FRAC LAP:2} if~$u$ is sufficiently smooth and with a controlled growth at infinity
(and, in fact, it is possible to set up a suitable notion of fractional Laplacian also
for functions that grow polynomially at infinity, see~\cite{POLYN}).
Besides, it is possible to provide a functional framework to define such operator in the weak sense
(see e.g.~\cite{MR2879266})
and a viscosity solution approach is often extremely appropriate
to construct general regularity theories (see e.g.~\cite{MR2494809}).

We refer to~\cite{MR2944369} for a gentle introduction to the fractional Laplacian.\medskip

{F}rom the point of view of the Fourier Transform,
denoted, as usual, by~$\widehat{\cdot}$ or by~${\mathcal{F}}$ (depending on the typographical
convenience), an instructive computation (see e.g.
Proposition~3.3 in~\cite{MR2944369}) shows that
$$ \widehat{ (-\Delta)^s u }(\xi)=c\,|\xi|^{2s} \,\widehat u(\xi),$$
for some~$c>0$. 
An appropriate choice of the normalization constant in~\eqref{FRAC LAP}
(also in dependence of~$n$ and~$s$)
allows us to take~$c=1$, and we will take this
normalization for the sake of simplicity
(and with the slight abuse of notation of dropping constants here and there).
With this choice, the fractional Laplacian in Fourier space is simply
the multiplication by the symbol~$|\xi|^{2s}$, consistently with the fact
that the classical Laplacian corresponds to the multiplication by~$|\xi|^2$.
In particular, the fractional Laplacian recovers\footnote{We think that it
is quite remarkable that the operator
obtained by the inverse Fourier Transform
of~$ \,|\xi|^{2} \,\widehat u$, the classical Laplacian,
reduces to a local operator. This is not true for
the inverse Fourier Transform
of~$ \,|\xi|^{2s} \,\widehat u$. In this spirit, it is interesting
to remark that the fact that
the classical Laplacian is a local operator is not immediate from
its definition in Fourier space, since computing
Fourier Transforms is always a nonlocal operation.}
the classical Laplacian as~$s\nearrow1$.
In addition, it satisfies the semigroup property, for any~$s$, $s'\in(0,1)$ with~$s+s'\le1$,
$$ {\mathcal{F}}{ (-\Delta)^s (-\Delta)^{s'}u } = |\xi|^{2s} \,{\mathcal{F}}((-\Delta)^{s'}u)=
|\xi|^{2s}\,|\xi|^{2s'}\, \widehat u=|\xi|^{2(s+s')}\,\widehat u=
{\mathcal{F}}{ (-\Delta)^{s+s'}u },
$$
that is
\begin{equation}\label{SEMI} (-\Delta)^s (-\Delta)^{s'}u = (-\Delta)^{s'} (-\Delta)^{s}u=(-\Delta)^{s+s'}u.\end{equation}
As a special case of~\eqref{SEMI}, when~$s=s'=1/2$, we have that the square root of the Laplacian
applied twice produces the classical Laplacian, namely
\begin{equation}\label{QUASDA}
\Big( (-\Delta)^{1/2}\Big)^2 =-\Delta.
\end{equation}
This observation gives that if~$U:\R^n\times[0,+\infty)\to\R$ is the harmonic
extension\footnote{Some \label{0qeiuru9467890000}
care has to be used with extension methods,
since the solution of~\eqref{EX:APP:E}
is not unique (if~$U$ solves~\eqref{EX:APP:E}, then so does~$U(x,y)+cy$
for any~$c\in\R$). The ``right'' solution of~\eqref{EX:APP:E}
that one has to take into account is the one with ``decay at infinity'',
or belonging to an ``energy space'', or obtained by convolution with
a Poisson-type kernel. See e.g.~\cite{MR3469920} for details.

Also, the extension method in~\eqref{EX:APP:E} and~\eqref{EX:APP:F}
can be related to an engineering application of the fractional
Laplacian motivated by the displacement of elastic membranes
on thin (i.e. codimension one) obstacles,
see~\cite{MR542512}. The intuition
for such application can be grasped from Figures~\ref{MK12},
\ref{SPIKES} and~\ref{MKXX}.
These pictures can be also useful to develop some
intuition about extension methods for fractional operators
and boundary reaction-diffusion equations.} of~$u:\R^n\to\R$, i.e. if
\begin{equation}\label{EX:APP:E}
\left\{
\begin{matrix}
\Delta U = 0 & {\mbox{ in }}\R^n\times[0,+\infty),\\
U(x,0)=u(x)& {\mbox{ for any }}x\in\R^n,
\end{matrix}
\right.\end{equation}
then
\begin{equation}\label{EX:APP:F}
-\partial_y U(x,0)=(-\Delta)^{1/2} u(x).
\end{equation}
See Appendix~\ref{EX:APP} for a confirmation of this.
In a sense, formula~\eqref{EX:APP:F} is a particular case of a general 
approach which reduces the fractional Laplacian to a local operator which is set
in a halfspace with an additional dimension and may be of singular or degenerate type,
see~\cite{MR2354493}.\medskip

As a rather approximative ``general nonsense'',
we may say that the fractional Laplacian shares some common feature with the classical
Laplacian.
In particular, both the classical and the fractional Laplacian are
invariant under translations and rotations. Moreover, a control
on the size of the fractional Laplacian of a function
translates, in view of~\eqref{2PO}, into a control of the oscillation of the function
(though in a rather ``global'' fashion): 
this ``democratic'' tendency
of the operator of ``averaging out'' any unevenness in the values of a function
is indeed typical of ``elliptic'' operators -- and the classical Laplacian is the prototype example
in this class of operators, while the fractional Laplacian is perhaps
the most natural fractional counterpart.

To make this counterpart more clear, we will say that a function~$u$
is $s$-harmonic in a set~$\Omega$ if~$(-\Delta)^s u=0$ at any point of~$\Omega$
(for simplicity, we take this notion in the ``strong'' sense, but equivalently
one could look at distributional definitions, see e.g. Theorem~3.12
in~\cite{MR1671973}). 

For example,
constant functions in $\R^n$ are $s$-harmonic in the whole space
for any $s\in(0,1)$, as both \eqref{FRAC LAP} and \eqref{FRAC LAP:2} imply.\medskip

Another similarity between classical and fractional Laplace equations is given
by the fact that notions like
those of fundamental solutions,
Green functions and Poisson kernels
are also well-posed in the fractional
case and somehow similar formulas
hold true, see e.g.
Definitions~1.7 and 1.8,
and Theorems~2.3, 2.10,
3.1
and~3.2 in~\cite{MR3461641}
(and related formulas hold true also for higher-order fractional operators,
see \cites{GD,AJS1,AJS2, AJS3}).\medskip

In addition, space inversions such as the Kelvin Transform also possess invariant properties
in the fractional framework, see e.g.~\cite{MR2256481}
(see also
Lemma~2.2 and Corollary~2.3 in~\cite{MR2964681}, and in addition
Proposition~A.1
on page~300 in~\cite{MR3168912} for a short proof).
Moreover, fractional
Liouville-type results hold under various assumptions,
see e.g.~\cite{MR3511811} and~\cite{POLYN}.\medskip

Another interesting link between classical and fractional operators
is given by subordination formulas which permit to reconstruct
fractional operators from the heat flow of classical operators,
such as
$$ (-\Delta)^s u=-\frac{s}{\Gamma(1-s)}\int_0^{+\infty}
t^{-1-s} \Big( e^{t\Delta}-1\Big)u\,dt,$$ 
see \cite{MR0115096}.\medskip

In spite of all these similarities, many important structural differences between
the classical and the fractional Laplacian arise. Let us list some of them.

\begin{difference}[Locality versus nonlocality]{\rm
The classical Laplacian of~$u$ at a point~$x$ only depends on the values of~$u$
in~$B_r(x)$, for any~$r>0$.

This is not true for the fractional Laplacian.
For instance, if~$u\in C^\infty_0(B_2, \,[0,1])$ with~$u=1$ in~$B_1$, we
have that, for any~$x\in\R^n\setminus B_4$,
\begin{equation}\label{DECAT}
-(-\Delta)^s u(x) = \PV \int_{\R^n} \frac{u(y)-u(x)}{|x-y|^{n+2s}}\,dy
=\int_{B_2} \frac{u(y)}{|x-y|^{n+2s}}\,dy\ge
\int_{B_1} \frac{dy}{\big(|x|+1\big)^{n+2s}}\ge\frac{\const}{|x|^{n+2s}}
\end{equation}
while of course~$\Delta u(x)=0$ in this setting.

It is worth remarking that the estimate in~\eqref{DECAT} is somewhat optimal.
Indeed, if~$u$ belongs to the Schwartz space (or
space of rapidly decreasing functions)
\begin{equation}\label{DECATS} {\mathcal{S}}:=
\left\{u\in C^{\infty }(\R^{n}){\mbox{ s.t. }}
\sup _{x\in \R^{n}}|x|^{\alpha }\left|D^{\beta }u(x)\right|
<+\infty {\mbox{ for all }} \alpha ,\beta \in \N^{n}\right\},
\end{equation}
we have that, for large~$|x|$,
\begin{equation}\label{DECAT:2}
\big|(-\Delta)^s u(x)\big|\le \frac{\const}{|x|^{n+2s}}
.\end{equation}
See Appendix~\ref{DECAT:3} for the proof of this fact. 
}\end{difference}

\begin{difference}[Summability assumptions]\rm
The pointwise computation of the classical Laplacian on a function $u$ 
does not require integrability properties on $u$.
Conversely, formula \eqref{FRAC LAP} for $u$ can make sense only when
\begin{equation*}
\int_{\R^n}\frac{|u(y)|}{1+{|y|}^{n+2s}}\;dy\ <\ +\infty
\end{equation*}
which can be read as a local integrability complemented by a growth condition at infinity.
This feature, which could look harmless at a first glance, can result problematic when 
looking for singular solutions to nonlinear problems 
(as, for example, in \cites{MR2985500,MR3393247} where there is 
an unavoidable integrability obstruction on a bounded domain)
or in ``blow-up'' type arguments (as mentioned in \cite{POLYN}, where the authors propose a way
to outflank this restriction).
\end{difference}

\begin{difference}[Computation along coordinate directions]{\rm
The classical Laplacian of~$u$ at the origin only depends on the values that~$u$
attains along the coordinate directions (or, up to a rotation, along a set of $n$ orthogonal directions).

\begin{figure}
    \centering
    \includegraphics[width=13cm]{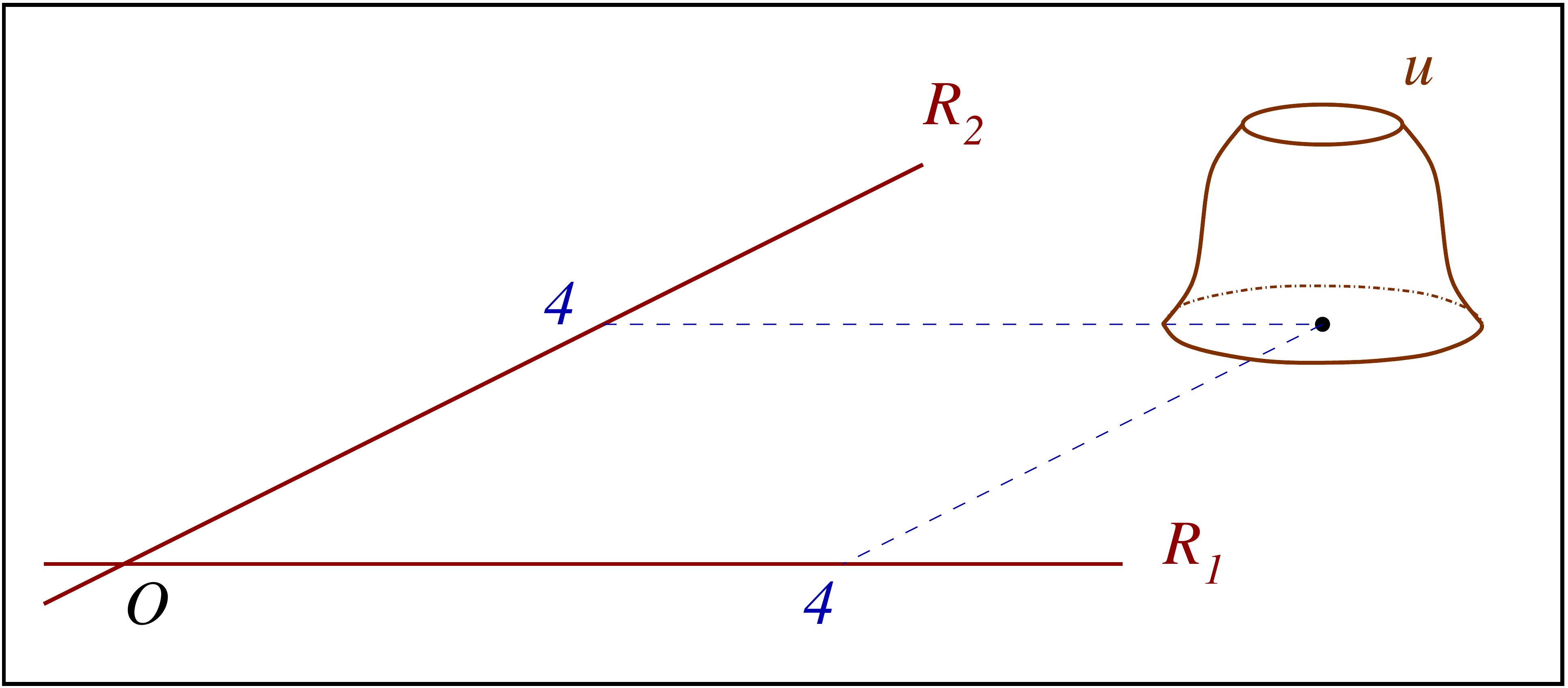}
    \caption{\it {{Coordinate directions not meeting a bump function.}}}
    \label{PUNTO}
\end{figure}

This is not true for the fractional Laplacian. As an example, let~$u\in C^\infty_0( B_2 (4e_1+4e_2),
\,[0,1])$, with~$u=1$ in~$B_1 (4e_1+4e_2)$. Let also~$R_j$ be the straight line in the $j$th coordinate
direction, that is
$$ R_j:=
\{ te_j,\; t\in\R\},$$
see Figure~\ref{PUNTO}.
Then
$$ R_j\cap B_2 (4e_1+4e_2)=\varnothing$$
for each~$j\in\{1,\dots,n\}$, and so~$u(te_j)=0$ for all~$t\in\R$ and~$j\in\{1,\dots,n\}$.
This gives that~$\Delta u(0)=0$.

On the other hand,
$$ \int_{\R^n} \frac{u(y)-u(0)}{|0-y|^{n+2s}}\,dy=
\int_{\R^n} \frac{u(y)}{|y|^{n+2s}}\,dy\ge
\int_{B_1 (4e_1+4e_2)} \frac{dy}{|y|^{n+2s}}>0,$$
which says that~$(-\Delta)^s u(0)\ne0$.
}\end{difference}

\begin{difference}[Harmonic versus $s$-harmonic functions]\label{DIALL}
{\rm If~$\Delta u(0)=1$, $\|u-v\|_{C^2(B_1)}\le\e$
and~$\e>0$ is sufficiently small
(see Figure~\ref{PUNTO2})
then~$\Delta v(0)\ge1-\const\e>0$, and in particular~$\Delta v(0)
\ne0$.

\begin{figure}
    \centering
    \includegraphics[width=9cm]{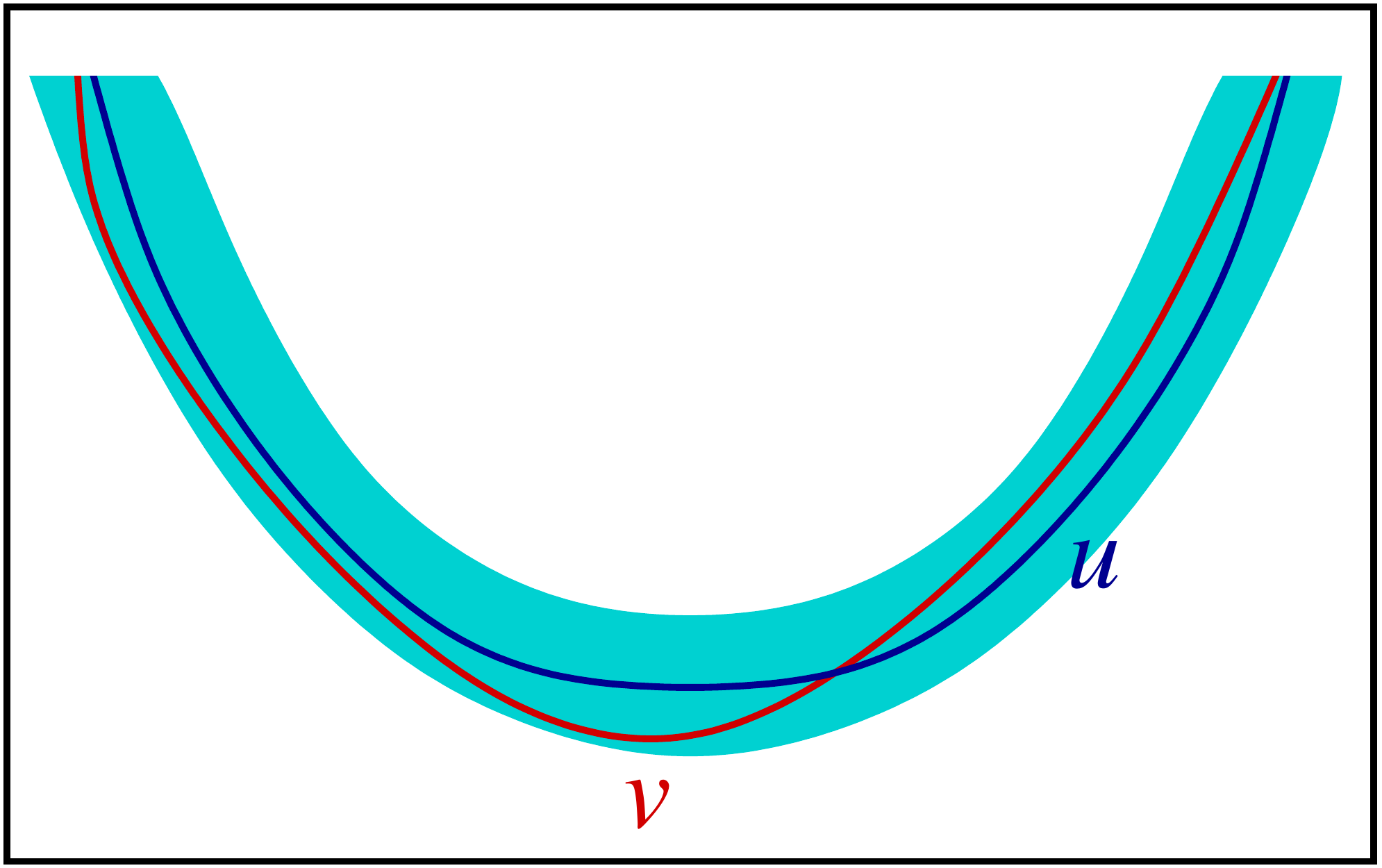}
    \caption{\it {{A function~$v$ which is ``close to~$u$''.}}}
    \label{PUNTO2}
\end{figure}

Quite surprisingly, this is not true for the fractional Laplacian.
More generally, in this case, as proved in~\cite{MR3626547},
for any~$\e>0$ and any (bounded, smooth)
function~$\bar u$, we can find~$v_\e$ such that
\begin{equation}\label{BARU}
\left\{\begin{matrix}
\|\bar u-v_\e\|_{C^2(B_1)}\le\e\\
{\mbox{and }} (-\Delta)^s v_\e=0{ \mbox{ in }}B_1.\end{matrix}\right.\end{equation}
A proof of this fact in dimension~$1$ for the
sake of simplicity is given in~\cite{HB} (the original paper~\cite{MR3626547}
presents a complete proof in any dimension).
See also~\cites{SALO1, SALO2, SALO3} for different approaches
to approximation methods in fractional settings which
lead to new proofs, and very refined and quantitative statements.\medskip

We also mention that the phenomenon
described in~\eqref{BARU} (which can
be summarized in the evocative statement
that
{\em all functions are locally $s$-harmonic (up to a small error)})
is very general, and it applies to
other nonlocal operators, also independently
from their possibly
``elliptic'' structure
(for instance all functions are locally $s$-caloric,
or $s$-hyperbolic, etc.). 
In this spirit, for completeness, in Section~\ref{CALORIC} we will establish the
density of fractional caloric functions in one space variable,
namely of the fact that for any~$\e>0$ and any (bounded, smooth)
function~$\bar u=\bar u(x,t)$, we can find~$v_\e=v_\e(x,t)$ such that
\begin{equation}\label{BARU:CALORIC}
\left\{\begin{matrix}
\|\bar u-v_\e\|_{C^2((-1,1)\times(-1,1))}\le\e\\
{\mbox{and }} \partial_t v_\e +(-\Delta)^s v_\e=0{ \mbox{ for any }}x\in(-1,1) {\mbox{ and any }}
t\in(-1,1).\end{matrix}\right.\end{equation}
We also refer to~\cite{SCALOR}
for a general approach and a series of general results on this type of approximation
problems with solutions of operators which are the
superposition of classical differential operators with fractional Laplacians.
Furthermore, 
similar results hold true for other nonlocal
operators with memory, see~\cite{ESAIM}. See in addition~\cites{2018arXiv180904005C,2018arXiv181007648K,2018arXiv181008448C}
for related results on higher order fractional operators.
}\end{difference}

\begin{difference}[Harnack Inequality]{\rm 
The classical Harnack Inequality says that if~$u$ is harmonic in~$B_1$
and~$u\ge0$ in~$B_1$ then
$$ \inf_{B_{1/2}} u\ge\const \sup_{B_{1/2}} u,$$
for a suitable universal constant, only depending on the dimension.

The same result is not true for $s$-harmonic functions.
To construct an easy counterexample, let~$\bar u(x)=|x|^2$
and, for a small~$\e>0$, let~$v_\e$ be as in~\eqref{BARU}.
Notice that, if~$x\in B_{1}\setminus B_{1/4}$
\begin{equation}\label{HANO} v_\e(x)\ge \bar u(x)-\| \bar u-v_\e\|_{L^\infty(B_1)}\ge \frac1{16}-\e>\frac1{32} \end{equation}
if~$\e$ is small enough, while
$$ v_\e(0)\le \bar u(0)+\| \bar u-v_\e\|_{L^\infty(B_1)}\le0+\e<\frac1{32}.$$
These observations imply that~$v_\e(0)<v_\e(x)$ for all~$x\in B_{1}\setminus B_{1/4}$
and therefore the infimum of~$v_\e$ in~$B_{1}$ is taken at some point~$\bar x$ in the closure of~$B_{1/4}$.
Then, we define
$$ u_\e(x):= v_\e(x)-\inf_{{B_1}} v_\e = v_\e(x)-v_\e(\bar x).$$
Notice that~$u_\e$ is $s$-harmonic in~$B_1$, since so is~$v_\e$,
and~$u_\e\ge0$ in~$B_1$. Also, $u_\e$ is strictly positive in~$B_{1}\setminus B_{1/4}$.
On the other hand, since~$\bar x\in B_{1/2}$
$$ \inf_{B_{1/2}} u_\e=u_\e(\bar x)=0,$$
which implies that~$u_\e$ cannot satisfy a Harnack Inequality
as the one in~\eqref{HANO}.
\medskip

In any case, it must be said that suitable Harnack Inequalities are valid also in the fractional case,
under suitable ``global'' assumptions on the solution: for instance, the
Harnack Inequality holds true for solutions that are positive in the whole of~$\R^n$
rather than in a given ball. We refer to~\cites{MR1941020, MR2817382}
for a comprehensive discussion on this topic and for recent developments.
}\end{difference}

\begin{difference}[Growth from the boundary]\label{009}{\rm
Roughly speaking,
solutions of Laplace equations have ``linear (i.e. Lipschitz) growth
from the boundary'', while solutions of fractional Laplace equations
have only H\"older growth from the boundary.
To understand this phenomenon, we point out that
if~$u$ is continuous in the closure of~$B_1$, with~$\Delta u=f$ in~$B_1$
and~$u=0$ on~$\partial B_1$, then
\begin{equation}\label{BOU:REG}
|u(x)|\le\const (1-|x|)\,\sup_{B_1}|f|.
\end{equation}
Notice that the term~$(1-|x|)$ represents the distance of
the point~$x\in B_1$ from~$\partial B_1$.
See e.g. Appendix~\ref{BOU:REG:S} for a proof of~\eqref{BOU:REG}.

The case of fractional equations is very different.
A first example which may be useful to keep in mind is that
the function
\begin{equation}\label{HHA}\begin{split}&
\R^n\ni x\mapsto (x_n)_+^s
\\&{\mbox{is $s$-harmonic
in the halfspace~$\{x_n>0\}$.}}
\end{split}\end{equation} 
For an elementary proof of this fact, see e.g. Section~2.4 in~\cite{MR3469920}.
Remarkably, the function in~\eqref{HHA} is only H\"older continuous
with H\"older exponent~$s$ near the origin.

Another interesting example is given by the function
\begin{equation}\label{CERCu12}
\R\ni x\mapsto u_{1/2}(x):=(1-|x|^2)_+^{1/2},\end{equation}
which satisfies
\begin{equation}\label{CERC}
(-\Delta)^{1/2} \,u_{1/2} =\const\quad{\mbox{ in }}\; (-1,1).
\end{equation}
A proof of~\eqref{CERC} based on extension methods and complex analysis
is given in Appendix~\ref{CERC:A}.

The identity in~\eqref{CERC}
is in fact a special case of a more general formula,
according to which the function
\begin{equation}\label{us}
\R^n\ni x\mapsto u_s(x):=(1-|x|^2)_+^s\end{equation}
satisfies
\begin{equation}\label{AHHA} (-\Delta)^s u_{s} =\const\quad{\mbox{ in }}\; B_1.\end{equation}
For this formula, and in fact
even more general ones, see~\cite{MR2974318}.
See also~\cite{MR0137148} for a probabilistic approach.

Interestingly, \eqref{HHA} can be obtained from~\eqref{AHHA}
by a blow-up at a point on the zero level set.
 
Notice also that
$$ \lim_{|x|\nearrow1}\frac{
|u_s(x)|}{1-|x|} =
\lim_{|x|\nearrow1}\frac{
(1-|x|^2)_+^s}{1-|x|} =
\lim_{|x|\nearrow1}\frac{
1}{(1-|x|)^{1-s}}=+\infty,$$
therefore, differently from the classical case, $u_s$ does not satisfy
an estimate like that in~\eqref{BOU:REG}.

It is also interesting to observe that the function~$u_{s}$ is related to the function~$x_+^{s}$
via space inversion (namely, a Kelvin transform) and integration, and indeed one can also deduce~\eqref{AHHA}
from~\eqref{HHA}: this fact was nicely
remarked to us by Xavier Ros-Oton and Joaquim Serra,
and the simple but instructive proof is sketched in Appendix~\ref{JSE}.

}\end{difference}

\begin{difference}[Global (up to the boundary) regularity]{\rm
Roughly speaking,
solutions of Laplace equations are ``smooth up to the boundary'', while
solutions of fractional Laplace equations
are not better than H\"older continuous at the boundary.
To understand this phenomenon, we point out that
if~$u$ is continuous in the closure of~$B_1$, 
\begin{equation}\label{LA:E10}
\left\{\begin{matrix}
\Delta u=f \;{\mbox{ in }} B_1,\\
u=0 \;{\mbox{ on }}
\partial B_1,\end{matrix}\right.\end{equation}
then
\begin{equation}\label{BOU:REG:LIP}
\sup_{x\in B_1} |\nabla u(x)|\le\const\,\sup_{B_1}|f|.
\end{equation}
See e.g. Appendix~\ref{BOU:REG:LIP:A} for a proof of this fact.

The case of fractional equations is very different since the function~$u_s$ in~\eqref{us}
is only H\"older continuous (with H\"older exponent~$s$) in~$B_1$, hence
the global Lipschitz estimate in~\eqref{BOU:REG:LIP} does not hold in this case.
This phenomenon can be seen as a counterpart of the one discussed in Difference~\ref{009}.
The boundary regularity for fractional Laplace problems is discussed in details in~\cite{MR3168912}.}\end{difference}

\begin{difference}[Explosive solutions]{\rm
Solutions of classical Laplace equations cannot attain
infinite values in the whole of the boundary. For instance,
if~$u$ is harmonic in~$B_1$, then
\begin{equation}\label{EXPNO}
\limsup_{\rho\nearrow1} \inf_{\partial B_\rho} u\le\const u(0).
\end{equation}
Indeed, by the Mean Value Property
for harmonic functions, for any~$\rho\in(0,1)$,
$$ u(0)=\frac{\const}{\rho^{n-1}} \int_{\partial B_\rho} u(x)\,d{\mathcal{H}}^{n-1}_x\ge
\inf_{\partial B_\rho} u,$$
from which~\eqref{EXPNO} plainly follows
(another proof follows by using the Maximum Principle instead of the
Mean Value Property).
On the contrary, and quite remarkably,
solutions of fractional Laplace equations
may ``explode'' at the boundary and~\eqref{EXPNO} can be violated
by $s$-harmonic functions in~$B_1$ which vanish outside~$B_1$.

For example, for 
\begin{equation}\label{1092373957349tufh39597306434676949767666794tfoe}
\R\ni x\mapsto u_{-1/2}(x):=\left\{\begin{matrix}
(1-|x|^2)^{-1/2} & {\mbox{ if }} x\in(-1,1),\\
0&{\mbox{otherwise}},\end{matrix}\right.\end{equation}
one has 
\begin{equation}\label{CERC-HARM}
(-\Delta)^{1/2} \,u_{-1/2} =0\quad{\mbox{ in }}\; (-1,1),
\end{equation}
and, of course, \eqref{EXPNO} is violated by~$u_{-1/2}$.
The claim in~\eqref{CERC-HARM} can be proven starting from \eqref{CERC} and
by suitably differentiating both sides of the equation:
the details of this computation can be found in Appendix \ref{CERC-HARM:A}.
For completeness,
we also give in Appendix~\ref{CERC-HARM:A:2} another proof of~\eqref{CERC-HARM}
based on complex variable and extension methods.

\begin{figure}
    \centering
    \includegraphics[width=12cm]{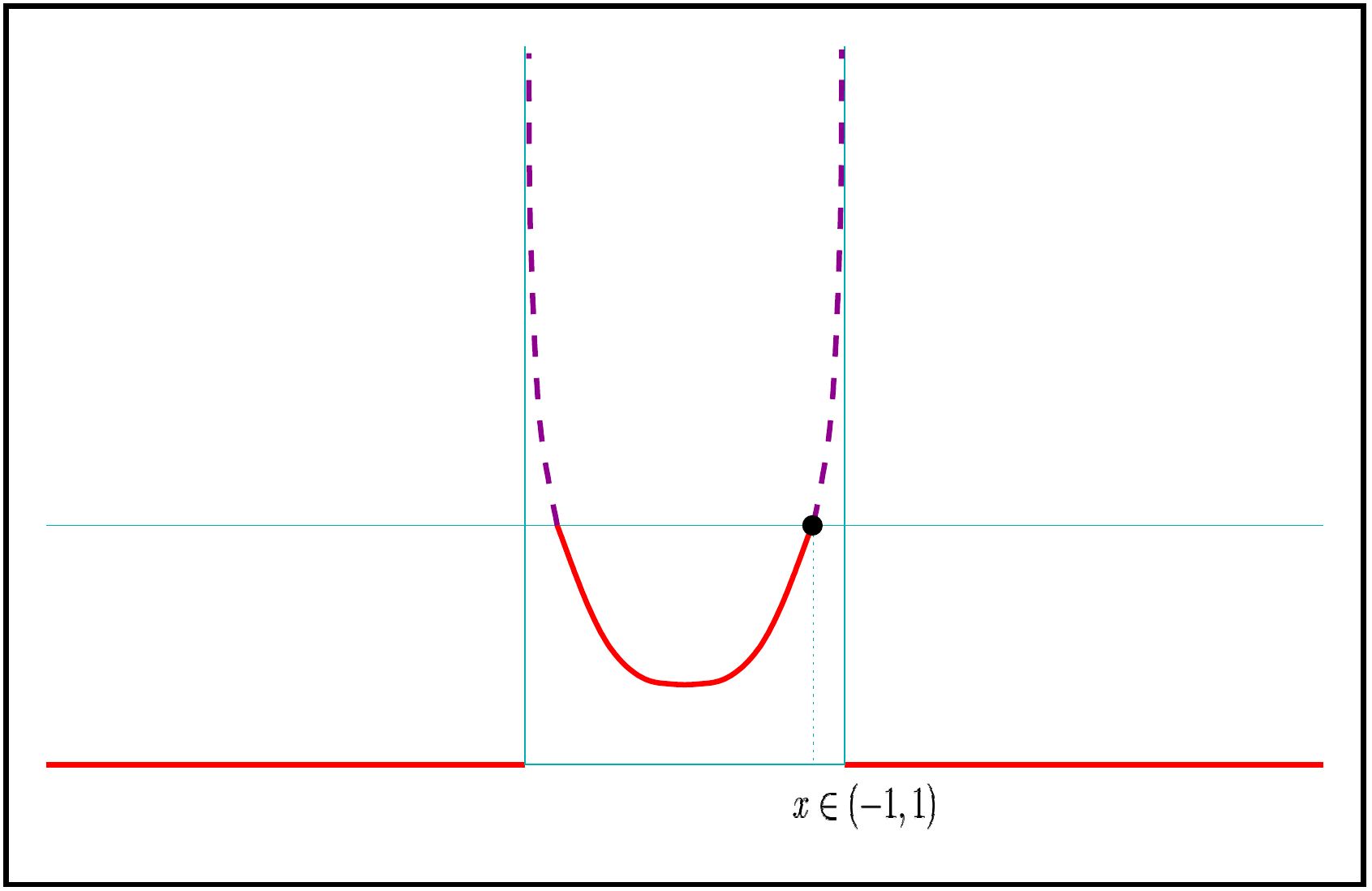}
    \caption{\it {{The function~$u_{-1/2}$ and the cancellation
occurring in~\eqref{CERC-HARM}.}}}
    \label{WI0we8484588585K3}
\end{figure}

A geometric interpretation of~\eqref{CERC-HARM}
is depicted in Figure~\ref{WI0we8484588585K3}
where a point~$x\in(-1,1)$ is selected and the graph of~$u_{-1/2}$
above the value~$u_{-1/2}(x)$ is drawn with a ``dashed curve''
(while a ``solid curve'' represents the graph of~$u_{-1/2}$
below the value~$u_{-1/2}(x)$): then, when computing the fractional Laplacian
at~$x$, the values coming from the dashed curve,
compared with~$u_{-1/2}(x)$, provide an opposite sign with
respect to the values coming from the solid curve.
The ``miracle'' occurring in~\eqref{CERC-HARM}
is that these two contributions with opposite sign
perfectly compensate and cancel each other, for any~$x\in(-1,1)$.
\medskip

More generally,
in every smooth bounded domain $\Omega\subset\R^n$ it is possible to build 
$s$-harmonic functions exploding at $\partial\Omega$ at the same rate as
dist$(\cdot,\partial\Omega)^{s-1}$. A phenomenon of this sort was spotted in
\cite{MR2985500}, and see~\cite{MR3393247} for the explicit explosion rate.
See~\cite{MR3393247} also
for a justification of the boundary behavior,
as well as the study of Dirichlet problems prescribing a singular boundary trace.
\medskip

Concerning this feature of explosive solutions
at the boundary, it is interesting to point out
a simple analogy with the classical Laplacian. Indeed, in view of~\eqref{HHA},
if~$s\in(0,1)$ and we take the function~$\R\ni
x\mapsto x_+^s$, we know that
it is~$s$-harmonic in~$(0,+\infty)$
and it vanishes on the boundary
(namely, the origin), and these features have a clear classical analogue
for~$s=1$.
Then, since for all~$s\in(0,1]$
the derivative of~$x_+^s$ is~$x_+^{s-1}$,
up to multiplicative constants,
we have that the latter is $s$-harmonic in~$(0,+\infty)$
and it
blows-up at the origin when~$s\in(0,1)$
(conversely, when $s=1$ one can do the same
computations but the resulting function is simply the
characteristic function of~$(0,+\infty)$
so no explosive effect arises). 

Similar computations
can be done in the unit ball instead of~$(0,+\infty)$,
and one simply gets functions that are bounded up
to the boundary when~$ s=1$, or explosive when~$s\in(0,1)$
(further details in
Appendices~\ref{CERC-HARM:A} and~\ref{CERC-HARM:A:2}).
}\end{difference}

\begin{difference}[Decay at infinity]{\rm The Gaussian~$e^{-|x|^2}$
reproduces the classical heat kernel.
That is, the solution of the heat equation
with initial datum concentrated at the origin, when considered at time~$t=1/4$,
produces the Gaussian (of course, the choice~$t=1/4$ is only for convenience,
any time~$t$ can be reduced to unit time by scaling the equation).

The fast decay prescribed by the Gaussian is special for the classical case
and the fractional case exhibits power law decays at infinity.
More precisely, let us consider
the heat equation
with initial datum concentrated at the origin, that is
\begin{equation}\label{DHE} \left\{\begin{matrix}
\partial_t u(x,t)=-(-\Delta)^s u(x,t) \;{\mbox{ for }}(x,t)\in \R^n\times(0,+\infty),\\
u(x,0)=\delta_0,
\end{matrix}\right.\end{equation}
and set
\begin{equation}\label{DHE:DGA}
{\mathcal{G}}_s(x)=u(x,1).\end{equation} By taking the
Fourier Transform of~\eqref{DHE}
in the~$x$ variable (and possibly neglecting normalization constants) one finds that
\begin{equation*} \left\{\begin{matrix}
\partial_t \hat u=-|\xi|^{2s}\hat u \;{\mbox{ in }} \R^n\times(0,+\infty),\\
u(\xi,0)=1,
\end{matrix}\right.\end{equation*}
hence
\begin{equation}\label{L12L12}
\hat u=e^{-|\xi|^{2s}t},\end{equation}
and consequently
\begin{equation}\label{SGA} {\mathcal{G}}_s(x)={\mathcal{F}}^{-1} (e^{-|\xi|^{2s}}),\end{equation}
being~${\mathcal{F}}^{-1}$ the anti-Fourier Transform of the Fourier Transform~${\mathcal{F}}$.
When~$s=1$, and neglecting the normalizing constants,
the expression in~\eqref{SGA} reduces to the Gaussian
(since the Gaussian is the Fourier Transform of itself).
On the other hand, as far as we know, there is no simple explicit representation of the fractional
heat kernel in~\eqref{SGA}, except in the ``miraculous'' case~$s=1/2$, in which~\eqref{SGA}
provides the explicit representation
\begin{equation}\label{GAUSS:s} 
{\mathcal{G}}_{1/2}(x)=\frac{\const}{\big(1+|x|^2\big)^{\frac{n+1}2}}.\end{equation}
See Appendix~\ref{GAUSS:sA} for a proof of~\eqref{GAUSS:s}
using Fourier methods and Appendix~\ref{GAUSS:FAC}
for a proof based on extension methods.

We stress that, differently from the classical case, the heat kernel~${\mathcal{G}}_{1/2}$
decays only with a power law. This is in fact a general feature of the fractional case,
since, for any~$s\in(0,1)$, it holds that 
\begin{equation}\label{LIMGA} \lim_{|x|\to+\infty} |x|^{n+2s} {\mathcal{G}}_s (x)=\const\end{equation}
and, for~$|x|\ge1$ and~$s\in(0,1)$, the heat kernel~${\mathcal{G}}_s (x)$
is bounded from below and from above by~$\frac{\const}{|x|^{n+2s}}$.

We refer to~\cite{MR1744782} for a detailed discussion on the fractional heat kernel.
See also~\cite{MR3211862} for more information on the fractional heat equation.
For precise asymptotics on fractional heat kernels,
see~\cite{zbMATH02598269,
MR0119247, MR1297844, 2018arXiv180311435D}.
\medskip

The decay of the heat kernel is also related to the associated distribution
in probability theory: as we will see in Section~\ref{HAJ103456},
the heat kernel represents the probability density of finding
a particle at a given point after a unit of time;
the motion of such particle is driven by a random walk in the classical
case and by a random process with long jumps in the fractional case
and, as a counterpart, the fractional probability distribution exhibits a ``long tail'',
in contrast with the rapidly decreasing classical one.\medskip

Another situation in which the classical case provides exponentially fast decaying solutions
while the fractional case exhibits polynomial tails is given by the Allen-Cahn equation
(see e.g. Section~1.1 in~\cite{MR2528756}
for a simple description of this equation also in view of phase coexistence models).
For concreteness, one can consider the one-dimensional equation
\begin{equation}\label{OEBN:s}
\left\{\begin{matrix}
(-\Delta)^s u =u-u^3\; {\mbox{ in }}\R,\\
\\
\dot{u}>0,\\
u(0)=0,\\
\displaystyle\lim_{t\to\pm\infty} u(t)=\pm1.
\end{matrix}
\right.\end{equation}
For~$s=1$, the system in~\eqref{OEBN:s} reduces to the pendulum-like system
\begin{equation}\label{OEBN:1}
\left\{\begin{matrix}
-\ddot u =u-u^3 \;{\mbox{ in }}\R,\\
\\
\dot{u}>0,\\
u(0)=0,\\
\displaystyle\lim_{t\to\pm\infty} u(t)=\pm1.
\end{matrix}
\right.\end{equation}
The solution of~\eqref{OEBN:1} is explicit and it has the form
\begin{equation}\label{11:OEBN:11}
u(t):=
\tanh\frac{t}{\sqrt{2}},
\end{equation}
as one can easily check. Also, by inspection, we see that such solution satisfies
\begin{equation}\label{GIAL}\begin{split}
&|u(t)-1|\le \const \exp({-\const\,t})\quad{\mbox{ for any }}t\ge1\\
{\mbox{and }}\;&|u(t)+1|\le \const \exp({-\const\,|t|})\quad{\mbox{ for any }}t\le-1.
\end{split}\end{equation}
Conversely, 
to the best
of our knowledge,
the solution of~\eqref{OEBN:s} has no simple explicit expression.
Also, remarkably, the solution of~\eqref{OEBN:s} decays to the equilibria~$\pm1$ only polynomially fast.
Namely, as proved in Theorem~2 of~\cite{MR3081641},
we have that the solution of~\eqref{OEBN:s} satisfies
\begin{equation}\label{SORGIA}\begin{split}
&|u(t)-1|\le\frac{\const}{t^{2s}}\quad{\mbox{ for any }}t\ge1\\
{\mbox{and }}\;&|u(t)+1|\le\frac{\const}{|t|^{2s}}\quad{\mbox{ for any }}t\le-1,
\end{split}\end{equation}
and the estimates in~\eqref{SORGIA} are optimal, namely it also holds that
\begin{equation}\label{GIA}\begin{split}
&|u(t)-1|\ge\frac{\const}{t^{2s}}\quad{\mbox{ for any }}t\ge1\\
{\mbox{and }}\;&|u(t)+1|\ge\frac{\const}{|t|^{2s}}\quad{\mbox{ for any }}t\le-1.
\end{split}\end{equation}
See Appendix~\ref{GIAGIA} for a proof of~\eqref{GIA}.
In particular, \eqref{GIA} says that
solutions of fractional Allen-Cahn equations such as the one in~\eqref{OEBN:s}
do not satisfy the exponential decay in~\eqref{GIAL} which is fulfilled
in the classical case.\medskip

The estimate in~\eqref{GIA} can be confirmed by looking at
the solution of the very similar equation
\begin{equation}\label{OEBN:11}
\left\{\begin{matrix}
(-\Delta)^s u =
\displaystyle\frac1{\pi}
\sin(\pi u)
\;{\mbox{ in }}\R,\\
\\
\dot{u}>0,\\
u(0)=0,\\
\displaystyle\lim_{t\to\pm\infty} u(t)=\pm1.
\end{matrix}
\right.\end{equation}
Though a simple expression of the solution of~\eqref{OEBN:11}
is not available in general, the ``miraculous'' case~$s=1/2$ possesses
an explicit solution, given by
\begin{equation}\label{78:87}
u(t):=\frac2\pi\,\arctan t.
\end{equation}
That~\eqref{78:87} is a solution of~\eqref{OEBN:11} when~$s=1/2$ is proved in
Appendix~\ref{0987poiutyui19}.
Another proof of this fact using~\eqref{GAUSS:s} is given in Appendix~\ref{APP:CC}.
\medskip

The reader should not be misled by the similar typographic forms of~\eqref{11:OEBN:11}
and~\eqref{78:87}, which represent two very different behaviors at infinity:
indeed
$$ \lim_{t\to+\infty} t\,\left(1- \frac2\pi\,\arctan t\right)=\frac2\pi,$$
and the function in~\eqref{78:87} satisfies the slow decay in~\eqref{GIA} (with~$s=1/2$)
and not the exponentially fast one in~\eqref{GIAL}.
\medskip

Equations like the one in~\eqref{OEBN:s} naturally arise, for instance, in long-range phase coexistence models
and in models arising in atom dislocation in crystals,
see e.g.~\cites{MR1442163, MR3296170}.
\medskip

A similar slow decay also occurs in the study of 
fractional Schr\"odinger operators, see e.g.~\cite{MR1054115}
and Lemma~C.1 in~\cite{MR3530361}.
For instance, the solution of
\begin{equation}
\label{L717}
(-\Delta)^s \Gamma+\Gamma=\delta_0\quad{\mbox{ in }}\delta_0\end{equation}
satisfies, for any~$|x|\ge1$,
$$ \Gamma(x)\simeq \frac{\const}{|x|^{n+2s}}.$$
A heuristic motivation for a bound of this type can be ``guessed'' from~\eqref{L717}
by thinking that, for large~$|x|$, the function~$\Gamma$
should decay more or less like~$(-\Delta)^s\Gamma$, which has ``typically''
the power law decay described in~\eqref{DECAT:2}.\medskip

If one wishes to keep arguing in this heuristic way, also the decays in~\eqref{LIMGA}
and~\eqref{GIA} may be seen as coming from an interplay between the right and the left
side of the equation, in the light of the decay of the fractional Laplace operator discussed
in~\eqref{DECAT:2}. For instance, to heuristically justify~\eqref{LIMGA},
one may think that the solution of the fractional heat equation
which starts from a Dirac's Delta, after a unit of time
(or an ``infinitesimal unit'' of time, if one prefers) has produced some bump, whose
fractional Laplacian, in view of~\eqref{DECAT:2}, may decay at infinity like~$\frac1{|x|^{n+2s}}$.
Since the time derivative of the solution has to be equal to that, the solution itself,
in this unit of time, gets ``pushed up'' by an amount like~$\frac1{|x|^{n+2s}}$
with respect to the initial datum, thus justifying~\eqref{LIMGA}.

A similar justification for~\eqref{GIA} may seem more tricky, since the decay
in~\eqref{GIA} is only of the type~$\frac1{|t|^{2s}}$ instead of~$\frac1{|t|^{1+2s}}$,
as the analysis in~\eqref{DECAT:2} would suggest. But to understand the problem,
it is useful to consider the derivative of the solution~$v:=\dot u$ and deduce from~\eqref{OEBN:s}
that
\begin{equation}\label{PUB} (-\Delta)^s v=
(-\Delta)^s \dot u =\dot u-3u^2\dot u=(1-3u^2) v.\end{equation}
That is, for large~$|t|$, the term~$1-3u^2$ gets close to~$1-3=-2$ and so the profile
at infinity may locally resemble the one driven by the equation~$(-\Delta)^s v=-2v$.
In this range, $v$ has to balance its fractional Laplacian, which is expected to decay like~$\frac1{|t|^{1+2s}}$,
in view of~\eqref{DECAT:2}. Then, since~$u$ is the primitive of~$v$,
one may expect that its behavior at infinity is related to the primitive of~$\frac1{|t|^{1+2s}}$,
and so to~$\frac1{|t|^{2s}}$, which is indeed the correct answer given by~\eqref{GIA}.

We are not attempting here to make these heuristic considerations rigorous,
but perhaps these kinds of comments may be useful in understanding why the behavior
of nonlocal equations is different from that of classical equations and to give at least
a partial justification of the delicate quantitative aspects involved in a rigorous
quantitative analysis (in any case, ideas like these are
rigorously exploited for instance in Appendix~\ref{GIAGIA}).

See also~\cite{MR3461371} for decay estimates
of ground states of a nonlinear nonlocal problem.\medskip

We also mention that other very interesting differences in the decay of solutions
arise in the study of different models for fractional porous medium equations,
see e.g.~\cites{MR2737788, MR2773189, MR3082241}.
}\end{difference}

\begin{difference}[Finiteness versus infiniteness
of the mean squared displacement]{\rm The mean squared displacement
is a useful notion to measure 
the ``speed of a diffusion process'',
or more precisely the
portion of the space that gets ``invaded'' at a given time by
the spreading of the diffusive quantity which is concentrated at
a point source at the initial time.
In a formula, if~$u(x,t)$ is the fundamental solution of the diffusion equation
related to the diffusion operator~$\LLL$, namely
\begin{equation}\label{F HEAT} \left\{
\begin{matrix}
\partial_t u = \LLL u &{\mbox{ for any $x\in\R^n$ and $t>0$}},\\
u(\cdot,0)=\delta_0(\cdot),
\end{matrix}
\right.\end{equation}
being~$\delta_0$ the Dirac's Delta,
one can define the mean squared displacement relative to
the diffusion process~$\LLL$ as
the ``second moment'' of~$u$ in the space variables, that is
\begin{equation}\label{MSD:MSD:1} {\rm MSD}_{\LLL}(t):=\int_{\R^n} |x|^2\, u(x,t)\,dx.\end{equation}
For the classical heat equation, by Fourier Transform one sees that,
when~$\LLL=\Delta$, the fundamental solution of~\eqref{F HEAT}
is given by the classical heat kernel
$$ u(x,t)= {\frac {1}{(4\pi t)^{n/2}}}e^{-\frac{|x|^{2}}{4t}},
$$
and therefore\footnote{See Appendix~A in~\cite{zbMATH05712796}
for a very nice explanation of the dimensional analysis and for a throughout discussion
of its role in detecting fundamental solutions.}
in such case, the substitution~$y:=\frac{x}{2\sqrt{t}}$
gives that
\begin{equation}\label{LINE}
{\rm MSD}_{\Delta}(t)=
\int_{\R^n}{\frac {|x|^2}{(4\pi t)^{n/2}}}e^{-\frac{|x|^{2}}{4t}}\,dx
=
\int_{\R^n}{\frac {
4 t |y|^2 }{ \pi^{n/2} }}e^{-|y|^2}\,dy = Ct,\end{equation}
for some~$C>0$. This says that the
mean squared displacement of the classical heat equation is finite,
and linear in the time variable.

On the other hand, in the fractional case in which~$\LLL=
-(-\Delta)^s$, by~\eqref{L12L12}
the fractional heat kernel is endowed with the scaling
property
$$ u(x,t)= \frac1{t^{\frac{n}{2s}}} {\mathcal{G}}_s 
\left( \frac{x}{t^{\frac1{2s}} }\right),$$
with~${\mathcal{G}}_s$ being as in~\eqref{DHE} and~\eqref{DHE:DGA}.
Consequently, in this case,
the substitution~$y:=
\frac{x}{t^{\frac1{2s}} }$ gives that
\begin{equation}\label{0o2r} {\rm MSD}_{-(-\Delta)^s}(t)=
\int_{\R^n} |x|^2\,
\frac1{t^{\frac{n}{2s}}} {\mathcal{G}}_s 
\left( \frac{x}{t^{\frac1{2s}} }\right)\,dx
= t^{\frac1{s}}\,
\int_{\R^n} |y|^2\, {\mathcal{G}}_s 
(y)\,dy.\end{equation}
Now, from \eqref{LIMGA},
we know that
$$ \int_{\R^n} |y|^2\, {\mathcal{G}}_s 
(y)\,dy=+\infty$$
and therefore we infer from~\eqref{0o2r}
that
\begin{equation}\label{810ehqwgoy8r3ryfefi34757575757} {\rm MSD}_{-(-\Delta)^s}(t)=+\infty.\end{equation}
This computation shows that, when~$s\in(0,1)$,
the diffusion process induced by~$-(-\Delta)^s$ does not
possess a finite mean squared displacement, in contrast with the classical
case in~\eqref{LINE}.
}\end{difference}

Other important differences
between the classical and fractional cases
arise in the study of nonlocal minimal
surfaces and in related fields:
just to list a few features, differently than in the classical
case, nonlocal minimal surfaces typically
``stick'' at the boundary, see~\cites{MR3516886, MR3596708, POICLA},
the gradient bounds of nonlocal minimal graphs are different
than in the classical case, see~\cite{DUCA},
nonlocal catenoids grow linearly and nonlocal
stable cones arise
in lower dimension, see~\cites{LAWSON, 2018arXiv181112141C},
stable surfaces of vanishing nonlocal mean curvature possess
uniform perimeter bounds, see Corollary~1.8 in~\cite{BVSUP},
the nonlocal mean curvature flow develops singularity
also in the plane, see~\cite{SIN},
its fattening phenomena are different, see~\cite{FATTE},
and the self-shrinking solutions are also different, see~\cite{2018arXiv181201847C},
and genuinely nonlocal phase transitions
present stronger rigidity properties than in the classical case,
see e.g. Theorem~1.2 in~\cite{owduir375957973}
and~\cite{FOFISEA}. Furthermore,
from the probabilistic viewpoint,
recurrence and transiency in long-jump stochastic processes
are different from the case of classical random walks, see e.g.~\cite{AFFILI}
and the references therein.
\medskip

We would like to conclude this list of differences with one similarity,
which seems to be not very well-known.
There is indeed a ``nonlocal representation'' for the classical Laplacian
in terms of a singular kernel. It reads as
\begin{equation}\label{LAP NONLOC}
-\Delta u(x) = \const\int_{\R^n}\frac{u(x+2y)+u(x-2y)-4u(x+y)-4u(x-y)+6u(x)}{{|y|}^{n+2}}\;dy.
\end{equation}
This one is somehow very close to~\eqref{FRAC LAP:2} with one important modification:
the difference operator in the numerator of the integrand has been increased in order,
in such a way that it is able to compensate the singularity of the kernel in $0$.
We include in Appendix~\ref{APP LAP NONLOC} a computation proving~\eqref{LAP NONLOC}
when $u$ is $C^{2,\alpha}$ around~$x$.
For a complete proof, involving Fourier transform techniques
and providing the explicit value of the constant, we refer to~\cite{AJS1}.

\subsection{The regional (or censored) fractional Laplacian}\label{DE:002}

A variant of the fractional Laplacian in~\eqref{FRAC LAP} consists in restricting
the domain of integration to a subset of~$\R^n$.
In this direction, an interesting operator 
is defined by the following singular integral:
\begin{equation} \label{REGION}
(-\Delta)^s_\Omega u(x) := \PV \int_\Omega \frac{u(x)-u(y)}{|x-y|^{n+2s}}\,dy.\end{equation}
We remark that when~$\Omega:=\R^n$ the regional fractional Laplacian in~\eqref{REGION}
boils down to the standard fractional Laplacian in~\eqref{FRAC LAP}.

In spite of the apparent similarity, the regional fractional Laplacian and the fractional
Laplacian are structurally two different operators.
For instance,
concerning Difference~\ref{DIALL}, we mention that solutions of regional fractional Laplace equations do not possess the same rich structure of those of fractional Laplace equations, and indeed 
\begin{equation}\label{NO:REGO}
\begin{split}&
{\mbox{it is not true that for any~$\e>0$ and any (bounded, smooth)
function~$\bar u$,}}\\&{\mbox{we can find~$v_\e$ such that}}\\ &
\left\{\begin{matrix}\|\bar u-v_\e\|_{C^2(B_1)}\le\e\\
{\mbox{and }} (-\Delta)^s_{\Omega} v_\e=0{ \mbox{ in }}B_1.\end{matrix}\right.
\end{split}\end{equation}
A proof of this observation will be given in Appendix~\ref{NO:REGOA}.\medskip

Interestingly, the regional fractional Laplacian turns out to be useful also in a possible setting
of Neumann-type conditions in the nonlocal case, as presented\footnote{Some colleagues pointed out to us that the use of~$R$ and~$r$
in some steps of formula~(5.5) of~\cite{MR3651008}
are inadequate.
We take this opportunity to amend such a flaw, presenting
a short proof of~(5.5) of~\cite{MR3651008}. Given~$\e>0$,
we notice that
\begin{eqnarray*}
&& I_1:=\iint_{{\Omega\times(\R^n\setminus\Omega)}\atop{ \{ |x-y|\ge\e\}}}
\frac{|u(x)-u(y)|^2}{|x-y|^{n+2s}}\,dx\,dy\le
\iint_{{\Omega\times(\R^n\setminus B_\e)}}
\frac{4\,\|u\|_{L^\infty(\R^n)}^2
\,dx\,d\zeta
}{|\zeta|^{n+2s}}\le \frac{
\const}{s\,\e^{2s}},
\end{eqnarray*}
where the constants are also allowed to depend on~$\Omega$ and~$u$.
Furthermore, if we define~$\Omega_\e$ to be the set of all the points
in~$\Omega$ with distance less than~$\e$ from~$\partial\Omega$,
the regularity of~$\partial\Omega$ implies that
the measure of~$\Omega_\e$ is bounded by~$\const\e$, and therefore
\begin{eqnarray*}
&& I_2:=\iint_{{\Omega\times(\R^n\setminus\Omega)}\atop{ \{ |x-y|<\e\}}}
\frac{|u(x)-u(y)|^2}{|x-y|^{n+2s}}\,dx\,dy\le
\iint_{{\Omega_\e\times B_\e(x)}}
\frac{4\,\|u\|_{C^1(\R^n)}^2\,|x-y|^2
\,dx\,dy
}{|x-y|^{n+2s}}\le
\int_{{B_\e}}
\frac{
\const\e \,d\zeta    
}{|\zeta|^{n+2s-2}}\le\frac{\const\e^{3-2s}}{1-s}.
\end{eqnarray*}
These observations imply that
\begin{eqnarray*}
&& \lim_{s\nearrow1}(1-s)\,\iint_{{\Omega\times(\R^n\setminus\Omega)}
} \frac{|u(x)-u(y)|^2}{|x-y|^{n+2s}}\,dx\,dy\le
\lim_{s\nearrow1}(1-s)
\left( \frac{\const}{s\,\e^{2s}}+\frac{\const\e^{3-2s}}{1-s}\right)
=\const\e.\end{eqnarray*}
Taking~$\e$ as small as we wish, we obtain
formula~(5.5) in~\cite{MR3651008}.}
in~\cite{MR3651008}.
Related to this, we mention that it is possible to obtain a 
regional-type operator starting from the classical Laplacian
coupled with Neumann boundary conditions (details about it will be given
in formula~\eqref{XAB} below).

\subsection{The spectral fractional Laplacian}\label{DE:003}

Another natural fractional operator arises
in taking fractional powers of the eigenvalues.
For this, we write
\begin{equation}\label{EXPA}
u(x,t)=\sum_{k=0}^{+\infty} u_k(t)\, \phi_k(x),\end{equation}
where~$\phi_k$ is the eigenfunction corresponding to the $k$th eigenvalue of the Dirichlet Laplacian,
namely
$$ \left\{
\begin{matrix}
-\Delta \phi_k=\lambda_k \phi_k {\mbox{ in }}\Omega\\
\phi_k\in H^1_0(\Omega).
\end{matrix}
\right.$$
with~$0<\lambda_0<\lambda_1\le\lambda_2\le\dots$. We normalize the sequence~$\phi_k$
to make it an orthonormal basis of~$L^2(\Omega)$
(see e.g. page~335
in~\cite{MR1625845}).
In this setting, we define
\begin{equation}\label{90ALA:123dfghjO0333}
(-\Delta)^s_{D,\Omega} u(x): = \sum_{k=0}^{+\infty} {\lambda_k^s}\,u_k(t)\,\phi_k(x).\end{equation}
We refer to~\cite{MR2754080} for extension methods for this type of operator.
Furthermore, other types of fractional operators can be defined in terms of different boundary conditions:
for instance, a spectral decomposition with respect to the eigenfunctions of the Laplacians
with Neumann boundary data naturally leads to an operator~$(-\Delta)^s_{N,\Omega}$
(and such operator also have applications in biology, see e.g.~\cite{MR3082317} and~\cite{SOA}).\medskip

It is also interesting to
observe that the spectral fractional Laplacian
with Neumann boundary conditions can also be
written in terms of a regional operator with a singular kernel.
Namely,
given an open and bounded set~$\Omega\subset\R^n$,
denoting by ${\Delta_{N,\Omega}}$ the Laplacian operator 
coupled with Neumann boundary conditions on $\partial\Omega$,
we let ${\{(\mu_j,\psi_j)\}}_{j\in\N}$ the pairs made up
of eigenvalues and eigenfunctions of ${-\Delta_{N,\Omega}}$, that is
$$ \left\{
\begin{matrix}
-\Delta \psi_j=\mu_j \psi_j {\mbox{ in }}\Omega \\
\partial_\nu\psi_j=0 \text{ on }\partial\Omega \\
\psi_j\in H^1(\Omega).
\end{matrix}
\right.$$
with~$0=\mu_0<\mu_1\le\mu_2\le\mu_3\le\dots$.

We define the following operator by making use of a spectral decomposition
\begin{equation}\label{spec-neu}
{(-\Delta)^s_{N,\Omega}}
\ :=\ \sum_{j=0}^{+\infty}\mu_j^s\hat u_j\:\psi_j,
\qquad \hat u_j=\int_\Omega u\psi_j,\ u\in C^\infty_0(\Omega).
\end{equation}
Comparing with~\eqref{90ALA:123dfghjO0333},
we can consider~$(-\Delta)^s_{N,\Omega}$ a spectral fractional
Laplacian with respect to classical Neumann data.
In this setting, the operator~$(-\Delta)^s_{N,\Omega}$ is
also an integrodifferential operator of regional type,
in the sense that one can write
\begin{equation}\label{XAB}
(-\Delta)^s_{N,\Omega}u(x)=\PV\int_\Omega\big(u(x)-u(y)\big)\,K(x,y)\,dy
,
\end{equation}
for a kernel~$K(x,y)$ which is comparable to~$\frac1{|x-y|^{n+2s}}$.
We refer to Appendix~\ref{09wdy8gfuiigsr3489trfhhdghhghghg}
for a proof of this.
\medskip

Interestingly, the fractional Laplacian and the spectral fractional Laplacian coincide, up
to a constant, for periodic functions, or functions defined on the flat torus, namely
\begin{equation}\label{TORUS}
{\mbox{if $u(x+k)=u(x)$ for any~$x\in\R^n$ and~$k\in\Z^n$, then~$(-\Delta)^s_{D,\Omega} u(x)=\const(-\Delta)^su(x)$.}}
\end{equation}
See e.g. Appendix~\ref{TORUS:A} for a proof of this fact.\medskip

On the other hand, striking differences between the fractional Laplacian
and the spectral fractional Laplacian hold true, see e.g.~\cites{MR3233760, MR3246044}.

Interestingly, it is not true that all functions are $s$-harmonic with respect to the
spectral fractional Laplacian, up to a small error, that is
\begin{equation}\label{BARU:NO}\begin{split}&
{\mbox{it is not true that for any~$\e>0$ and any (bounded, smooth)
function~$\bar u$,}}\\&{\mbox{we can find~$v_\e$ such that}}\\ &
\left\{\begin{matrix}\|\bar u-v_\e\|_{C^2(B_1)}\le\e\\
{\mbox{and }} (-\Delta)^s_{D,\Omega} v_\e=0{ \mbox{ in }}B_1.\end{matrix}\right.
\end{split}\end{equation}
A proof of this will be given in Appendix~\ref{BARU:NOA}.
The reader can easily compare~\eqref{BARU:NO}
with the setting for the fractional Laplacian discussed in Difference~\ref{DIALL}.\medskip

Remarkably, in spite of these differences, the spectral fractional Laplacian
can also be written as an integrodifferential operator of the form
\begin{equation}\label{LAKPA} \PV \int_{\Omega} \big(u(x)-u(y)\big)\,K(x,y)\,dy
+\beta(x)u(x), \end{equation}
for a suitable kernel~$K$ and potential~$\beta$, see 
Lemma~38 in~\cite{MR3610940} or Lemma~10.1 in~\cite{2016arXiv161009881B}.
This can be proved with analogous computations to those performed in the 
case of the regional fractional Laplacian in the previous paragraph.

\subsection{Fractional time derivatives}\label{w9eurgh}

The operators described in Sections in~\ref{DE:001}, \ref{DE:002} and~\ref{DE:003}
are often used in the mathematical description of anomalous types of diffusion
(i.e. diffusive processes which produce important differences with respect to the classical
heat equation, as we will
discuss in Section~\ref{PRO:SE}): the main role of
such nonlocal operators is usually
to produce a different behavior of the diffusion process with respect to the space variables.

Other types of anomalous diffusions arise from non-standard behaviors with respect
to the time variable. These aspects are often the mathematical counterpart
of memory effects. As a prototype example, we recall the notion of Caputo
fractional derivative,
which, for any~$t>0$ (and up to normalizing factors that we omit for simplicity) is given by
\begin{equation}\label{0we0wy83fguguqgferighierhg0011} \partial^s_{C,t} u(t):=
\int_0^t \frac{\dot u (\tau) }{(t-\tau)^{s}} \, d\tau.\end{equation}
We point out that, for regular enough functions~$u$,
\begin{equation}\label{PRE:CA} 
\begin{split}\partial^s_{C,t} u(t)\,&=
\int_0^t \frac{\dot u (\tau) }{(t-\tau)^{s}} \, d\tau\\
&=\int_0^t \left(\frac{d}{d\tau}\frac{\big(u (\tau) -u(t)\big)}{(t-\tau)^{s}} 
-s\,\frac{\big(u (\tau) -u(t)\big)}{(t-\tau)^{1+s}}
\right)\, d\tau\\
&=\frac{u (t) -u(0)}{t^{s}}-\lim_{\tau\to t} \frac{u (t) -u(\tau)}{(t-\tau)^{s}} 
-s\,\int_0^t \frac{\big(u (\tau) -u(t)\big)}{(t-\tau)^{1+s}}\, d\tau\\
&=\frac{u (t) -u(0)}{t^{s}}-\dot u(t)\,\lim_{\tau\to t}(t-\tau)^{1-s} 
-s\,\int_0^t \frac{u (\tau) -u(t)}{(t-\tau)^{1+s}}\, d\tau\\
&= 
\frac{u (t) -u(0)}{t^{s}}+s\,\int_0^t \frac{u (t) -u(\tau)}{(t-\tau)^{1+s}}\, d\tau.
\end{split}\end{equation}
Though in principle this expression takes into account only the values of~$u(t)$
for~$t\ge0$, hence~$u$ does not need to be defined for negative times,
as pointed out e.g. in Section~2 of~\cite{MR3488533},
it may be also convenient to constantly extend~$u$ in~$(-\infty,0)$.
Hence, we take the convention for which~$u(t)=u(0)$ for any~$t\le0$.
With this extension, one has that, for any~$t>0$,
\begin{eqnarray*}
&& s\,\int_{-\infty}^0 \frac{ u (t) -u(\tau) }{(t-\tau)^{1+s}}\, d\tau
=s\,\int_{-\infty}^0 \frac{u (t)-u(0)}{(t-\tau)^{1+s}}\, d\tau=\frac{u(t)-u(0)}{t^s}.
\end{eqnarray*}
Hence, one can write~\eqref{PRE:CA} as
\begin{equation}\label{PRE:CAq2} 
\partial^s_{C,t} u(t)=
s\,\int_{-\infty}^t \frac{u (t) -u(\tau)}{(t-\tau)^{1+s}}\, d\tau.
\end{equation}
This type of formulas also relates the Caputo derivative to the so-called Marchaud derivative,
see e.g.~\cite{MR1347689}. 

In the literature, one can also consider higher order
Caputo derivatives, see e.g.~\cites{MR2090004, MR1829592}
and the references therein.

Also, it is useful to consider the Caputo derivative
in light of the (unilateral) Laplace Transform
(see e.g. Chapter~2.8
in~\cite{MR1658022}, and~\cite{iqo0eiKAL})
\begin{equation}\label{L7849576AP:L}
{\mathcal{L}} u(\omega):=
\int_{0}^{+\infty } u(t)\,e^{-\omega t}\,dt.\end{equation}
With this notation, up to dimensional constants, one can write 
(for a smooth function with exponential
control at infinity) that
\begin{equation}\label{LAP:L}
{\mathcal{L}} (\partial^s_{C,t} u)(\omega)=\omega^s
{\mathcal{L}} u(\omega)-\omega^{s-1} u(0)
,\end{equation}
see Appendix~\ref{LAP:LA} for a proof.
\medskip

In this way, one can also link equations driven by the Caputo derivative
to the so-called Volterra integral equations: namely one can invert the expression~$
\partial^s_{C,t} u=f$ by 
\begin{equation}\label{LAP:L2}
u(t)=u(0)+C\int_0^t \frac{f(\tau)}{(t-\tau)^{1-s}}\,d\tau
,\end{equation}
for some normalization constant~$C>0$,
see Appendix~\ref{LAP:LA2} for a proof.
\medskip

It is also worth mentioning that the Caputo derivative of order~$s$
of a power gives, up to normalizing constants, the ``power minus~$s$'': more precisely,
by~\eqref{0we0wy83fguguqgferighierhg0011}
and using the substitution~$\vartheta:=\tau/t$, we see that, for any~$r>0$,
$$
\partial^s_{C,t} (t^r)=r\,
\int_0^t \frac{\tau^{r-1} }{(t-\tau)^{s}} \, d\tau=r\,t^{r-s}\,
\int_0^1 \frac{\vartheta^{r-1} }{(1-\vartheta)^{s}} \, d\vartheta= Ct^{r-s},
$$
for some~$C>0$.\medskip

Moreover, in relation to the comments on page~\pageref{F HEAT},
we have that
\begin{equation}\label{MSDES}
\begin{split}&{\mbox{the mean squared displacement
related to the diffusion operator}}\\
& \left\{
\begin{matrix}
\partial^s_{C,t} u = \Delta u &{\mbox{ for any $x\in\R^n$ and $t>0$}},\\
u(\cdot,0)=\delta_0(\cdot),
\end{matrix}
\right. \\
&{\mbox{is finite and proportional to~$t^s$}}.\end{split}
\end{equation}
See Appendix~\ref{MSDES:A} for a proof of this.\medskip

The Caputo derivatives describes a process ``with memory'', in the sense
that it ``remembers the past'', though ``old events count less than recent ones''. We sketch a memory effect of Caputo type in Appendix~\ref{MEMORY}.

Due to its memory effect, operators related to Caputo derivatives have found several applications
in which the basic parameters of a physical system change in time, in view of the evolution
of the system itself: for instance, in studying flows in porous media, when
time goes, the fluid may either ``obstruct'' the holes of the medium, 
thus slowing down the diffusion, or ``clean'' the holes,
thus making
the diffusion faster, and the Caputo derivative may be a convenient approach to describe
such modification in time of the diffusion coefficient, see~\cite{MR2379269}.
\medskip

Other applications of Caputo derivatives occur in biology and neurosciences,
since the network of neurons exhibit time-fractional diffusion, also in view of
their highly ramified structure, see e.g.~\cite{PuPOw2e3we}
and the references therein.

We also refer to~\cites{MR3469920, MR3177769, 2017arXiv170608241V} and to the references
therein for further discussions on different types of anomalous diffusions.

\section{A more general point of view: the ``master equation''}

The operators discussed in Sections~\ref{DE:001},
\ref{DE:002}, \ref{DE:003}
and~\ref{w9eurgh} can be framed into a more general setting,
that is that of the ``master equation'', see e.g.~\cite{MR3329847}.

Master equations describe the evolution of a quantity
in terms of averages in space and time of the quantity itself.
For concreteness one can consider a quantity~$u=u(x,t)$
and describe its evolution by
an equation of the kind
\begin{equation*}
c\partial_t u(x,t)= Lu(x,t)+f\big(x,t,u(x,t)\big)
\end{equation*}
for some~$c\in\R$ and a forcing term~$f$,
and the operator~$L$ has the integral form
\begin{equation}\label{MASTER} Lu(x,t):=\iint_{\R^n\times(0,+\infty)}
\big( u(x,t)-u(x-y,t-\tau)\big)\,
{\mathcal{K}}(x,t,y,\tau)\,d\mu(x,\tau),\end{equation}
for a suitable measure~$\mu$ (with the integral
possibly taken in the principal value sense, which is
omitted here for simplicity; also one can consider even more general operators by
taking actions
different than translations and more general ambient spaces).

Though the form of such operator
is very general, one can also consider simplifying structural
assumptions. For instance, one can take~$\mu$ to be the space-time
Lebesgue measure over~$\R^n\times(0,+\infty)$, namely
\begin{equation*}
d\mu(x,\tau)=dx\,d\tau.
\end{equation*}
Another common simplifying assumption is to assume that
the kernel is induced by an uncorrelated effect of the space and time
variables, with the product structure
\begin{equation*}
{\mathcal{K}}(x,t,y,\tau)={\mathcal{K}}_{\rm space}(x,y)\,{\mathcal{K}}_{\rm time}(t,\tau).
\end{equation*}
The fractional Laplacian of Section~\ref{DE:001}
is a particular case of this setting (for functions depending on the space variable), with the choice,
up to normalizing constants,
$${\mathcal{K}}_{\rm space}(x,y):=\frac{1}{|y|^{n+2s}}.$$
More generally, for~$\Omega\subseteq\R^n$, the regional fractional
Laplacian in Section~\ref{DE:002} comes from the choice
$${\mathcal{K}}_{\rm space}(x,y):=\frac{\chi_\Omega(x-y)}{|y|^{n+2s}}.$$
Finally, in view of~\eqref{PRE:CAq2}, for time-dependent functions,
the choice
$${\mathcal{K}}_{\rm time}(t,\tau):=\frac{\chi_{(-\infty,t)}(\tau)}{|\tau|^{1+s}}.$$
produces the Caputo
derivative discussed in
Section~\ref{w9eurgh}.\medskip

We recall that one of 
the fundamental structural
differences in partial differential equations
consists in the distinction between operators ``in divergence form'',
such as
\begin{equation}\label{DIU8234j:1}-\sum_{i,j=1}^n
\frac{\partial}{\partial x_i}\left( a_{ij}(x) \frac{\partial u}{\partial x_j}(x)\right)
\end{equation}
and those ``in non-divergence form'', such as
\begin{equation}\label{DIU8234j:2}-\sum_{i,j=1}^n
a_{ij}(x) \frac{\partial^2 u}{\partial x_i\partial x_j}(x).
\end{equation}
This structural difference can also be recovered from the master equation.
Indeed, if we consider a (say, for the sake of concreteness,
strictly positive, bounded and smooth)
matrix function~$M:\R^n\to{\rm Mat }(n\times n)$, we can take into account the master spatial operator induced by the kernel
\begin{equation}\label{KERVA23}
{\mathcal{K}}_{\rm space}(x,y):=\frac{1-s}{|M(x-y,y)\,y|^{n+2s}},\end{equation}
that is, in the notation of~\eqref{MASTER},
\begin{equation}\label{DIU8234j:3}
(1-s)\,\int_{\R^n} \frac{ u(x)-u(x-y) }{|M(x-y,y)\, y|^{n+2s}}\,dy.
\end{equation}
Then, up to a normalizing constant, if
\begin{equation}\label{DIU8234j:IPOT}
M(x-y,y)=M(x,-y),\end{equation}
then
\begin{equation}\label{DIV:EQ}
 \begin{split}
&{\mbox{the limit as~$s\nearrow1$ of the operator in~\eqref{DIU8234j:3}}}\\
&{\mbox{recovers the classical divergence form operator in~\eqref{DIU8234j:1},}}\\
&{\mbox{with }} a_{ij}(x):= \const\int_{S^{n-1}} \frac{\omega_i\, \omega_j}{
|M(x,0)\,\omega|^{n+2} }\,d{\mathcal{H}}^{n-1}_{\omega}.
 \end{split}
\end{equation}
A proof of this will be given in Appendix~\ref{DIV:FORM}.

It is interesting to observe that condition~\eqref{DIU8234j:IPOT} says that,
if we set~$z:=x-y$, then 
\begin{equation}
\label{BL33}
M(z,x-z)=M(x,z-x)
\end{equation}
and so
the kernel in~\eqref{KERVA23} is invariant by exchanging~$x$ and~$z$.
This invariance naturally leads to a (possibly formal) energy functional
of the form
\begin{equation}\label{E:2 en}
\frac{1-s}{2}\,\iint_{\R^n\times\R^n} \frac{ \big( u(x)-u(z) \big)^2}{|M(z,x-z)\,(x-z)|^{n+2s}}\,dx\,dz.\end{equation}
We point out that condition~\eqref{BL33} translates, roughly speaking, into
the fact that the energy density in~\eqref{E:2 en} ``charges
the variable~$x$ as much as the variable~$z$''.\medskip

The study of the energy functional in~\eqref{E:2 en}
also drives to a natural quasilinear generalization,
in which the fractional energy takes the form
\[
\int_{\R^n} \frac{ \Phi\big( u(x)-u(z) \big)}{|M(z,x-z)\,(x-z)|^{n+2s}}\,dx\,dz,\]
for a suitable~$\Phi$, see e.g.~\cites{MR3339179, MR3456825}
and the references therein for further details on quasilinear nonlocal operators.
See also~\cite{MR3177769} and the references therein for other type of nonlinear fractional equations.
\medskip

Another case of interest (see e.g.~\cite{MR1918242})
is the one in which one considers the master equation driven by the spatial kernel
$${\mathcal{K}}_{\rm space}(x,y):=\frac{1-s}{|M(x,y)\, y|^{n+2s}}\,dy,$$
that is, in the notation of~\eqref{MASTER},
\begin{equation}\label{DIU8234j:4}
(1-s)\,\int_{\R^n} \frac{
u(x)-u(x-y)
}{|M(x,y)\,y|^{n+2s}}\,dy.
\end{equation}
Then, up to a normalizing constant, 
if
\begin{equation}\label{DIU8234j:IPOT:2}
M(x,y)=M(x,-y),\end{equation} 
then
\begin{equation}\label{NONDIV:EQ}
 \begin{split}
&{\mbox{the limit as~$s\nearrow1$ of the operator in~\eqref{DIU8234j:4}}}\\
&{\mbox{recovers the classical non-divergence form operator in~\eqref{DIU8234j:2},}}
\\&{\mbox{with }} a_{ij}(x):= \const
\int_{ S^{n-1}}\frac{
\omega_i \omega_j}{ |M(x,0)\,\omega|^{n+2} }\,d{\mathcal{H}}^{n-1}_{\omega}
 \end{split}
\end{equation}
A proof of this will be given in Appendix~\ref{NONDIV:FORM}.
\medskip

We recall that nonlocal linear operators in non-divergence form
can also be useful in the definition of fully nonlinear nonlocal operators,
by taking appropriate infima and suprema of combinations of linear operators,
see e.g.~\cite{MR3148110} and the references therein for further discussions about this
topic (which is also related to 
stochastic games).\medskip

We also remark that
understanding the role of the affine transformations of the spaces
on suitable nonlocal operators (as done for instance
in~\eqref{DIU8234j:4} and~\eqref{DIU8234j:4}) often permits a deeper
analysis of the problem in nonlinear settings too, see e.g. the
very elegant way in which a fractional Monge-Amp\`ere equation is introduced in~\cite{MR3479063}
by considering
the infimum of fractional linear operators
corresponding to all affine transformations of determinant one of a given multiple of the fractional
Laplacian.
\medskip

As a general comment, we also think that an interesting consequence
of the considerations given in this section is that classical, local equations
can also be seen as a limit approximation of more general master equations.\medskip

We mention that there are also many other interesting kernels, both in space and time,
which can be taken into account in integral equations. Though we focused here mostly
on the case of singular kernels, there are several important problems that focus
on ``nice'' (e.g. integrable) kernels, see e.g.~\cites{MR2722295, MR2958346, MR3641640}
and the references therein.

As a technical comment let us point out that, in a sense, the nice kernels may have
computational advantages, but may provide loss of compactness and loss of regularity
issues: roughly speaking, convolutions with smooth kernel are always smooth,
thus any smoothness information on a convolved function gives little information
on the smoothness of the original function -- viceversa, if the convolution of an ``object'' with a singular
kernel is smooth, then it means that the original object has a ``good order of vanishing at the origin''.
When the original object is built by the difference of a function and its translation, such
vanishing implies some control of the oscillation of the function, hence
opening a door towards a regularity result.

\section{Probabilistic motivations}\label{PRO:SE}

We provide here some elementary, and somewhat
heuristic, motivations for the operators described in Section~\ref{SOME}
in view of
probability and statistics applications. The treatment of this section is mostly colloquial
and not to be taken at a strictly rigorous level (in particular, all functions
are taken to be smooth, some uniformity problems are neglected, convergence
is taken for granted, etc.). See e.g.~\cite{MR1409607} for rigorous
explanations linking
pseudo-differential operators and Markov/L\'evy processes.
See also~\cites{MR1406564, MR2345912, MR2512800, MR2584076, MR2963050}
for other perspectives and links between probability and fractional calculus
and~\cite{K201717} for a complete
survey
on jump processes and
their connection to nonlocal operators.

The probabilistic approach to study nonlocal effects
and the analysis of distributions with polynomial tails
are also some of the cornerstones
of the application of mathematical theories to finance,
see e.g.~\cites{PARETO, MANDEL}, and
models with jump process
for prices
have been proposed in~\cite{COX}.

\subsection{The heat equation and the classical Laplacian}\label{SE:HE}

The prototype of parabolic equations is the heat equation
\begin{equation}\label{HC}
\partial_t u(x,t)=c\,\Delta u(x,t)
\end{equation}
for some~$c>0$.
The solution~$u$ may represent, for instance, a temperature, and the foundation
of~\eqref{HC} lies on two basic assumptions:
\begin{itemize}
\item the variation of $u$ in a given region~$U\subset\R^n$ is due to the flow of some
quantity~$v:\R^n\to\R^n$ through~$U$,
\item $v$ is produced by the local variation of~$u$.
\end{itemize}
The first ansatz can be written as
\begin{equation}\label{BH10} \partial_t \int_U u(y,t)\,dy = \int_{\partial U} v(y,t)\cdot
\nu(y)\,d{\mathcal{H}}^{n-1}_y,\end{equation}
where~$\nu$ denotes the exterior normal vector of~$U$ and~${\mathcal{H}}^{n-1}$
is the standard $(n-1)$-dimensional surface Hausdorff measure.

The second ansatz can be written as~$v=c\nabla u$, which combined with~\eqref{BH10} 
and the Divergence Theorem gives that
$$ \partial_t \int_U u(y,t)\,dy = c \int_{\partial U} \nabla u(y,t)\cdot
\nu(y)\,d{\mathcal{H}}^{n-1}_y= c \int_U {\rm div}\,(\nabla u(y,t))\,dy=
c\int_U \Delta u(y,t)\,dy.$$
Since~$U$ is arbitrary, this gives~\eqref{HC}.\medskip

Let us recall a probabilistic interpretation of~\eqref{HC}. 
The idea is that~\eqref{HC} follows by taking suitable limits of a
discrete ``random walk''. For this, we take a small space scale~$h>0$
and a time step
\begin{equation}\label{KAL}
\tau=h^2.
\end{equation}
We consider the random motion of a particle in the lattice~$h\Z^n$, as follows.
At each time step, the particle can move in any coordinate direction with equal probability.
That is, a particle located at~$h\bar k\in h\Z^n$ at time~$t$
is moved to one of the~$2n$ points~$h\bar k\pm he_1$, $\dots$, $h\bar k\pm h e_n$
with equal probability (here, as usual, $e_j$ denotes the $j$th element of the
standard Euclidean basis of~$\R^n$).

We now look at the expectation to find the particle at a point~$x\in h\Z^n$ at time~$t\in\tau\N$.
For this, we
denote by~$u(x,t)$ the probability density of such expectation.
That is, the probability for the particle of lying in the spatial region~$B_r(x)$
at time~$t$ is, for small~$r$, comparable with
$$ \int_{B_r(x)} u(y,t)\,dy.$$
Then, the probability of finding a particle at
the point~$x\in h\Z^n$
at time~$t+\tau$ is the sum of the probabilities of finding the particle at a closest neighborhood of~$x$
at time~$t$, times the probability of jumping from this site to~$x$. That is,
\begin{equation}\label{AL}
u(x,t+\tau)=\frac{1}{2n} \sum_{j=1}^n \Big(u(x+h e_j)+u(x-he_j)\Big).\end{equation}
Also,
\begin{eqnarray*}
&& u(x+h e_j)+u(x-he_j)-2u(x,t)\\&=&
\left( u(x,t)+h\nabla u(x,t)\cdot e_j+\frac{h^2 \,D^2 u(x,t)\,e_j\cdot e_j}{2}\right)
+
\left( u(x,t)-h\nabla u(x,t)\cdot e_j+\frac{h^2 \,D^2 u(x,t)\,e_j\cdot e_j}{2}\right)\\&&\qquad-2u(x,t)+O(h^3)
\\&=&
h^2\,\partial^2_{x_j}u(x,t)+O(h^3).
\end{eqnarray*}
Thus, subtracting~$u(x,t)$ to both sides
in~\eqref{AL}, dividing by~$\tau$, recalling~\eqref{KAL},
and taking the limit (and neglecting any possible
regularity issue), we formally find that
\begin{eqnarray*}
\partial_t u(x,t) &=&\lim_{\tau\searrow0}\frac{u(x,t+\tau)-u(x,t)}{\tau}\\&=&
\lim_{h\searrow0}
\frac{1}{2n} \sum_{j=1}^n \frac{u(x+h e_j)+u(x-he_j)-2u(x,t)}{h^2}\\
&=& \lim_{h\searrow0}
\frac{1}{2n} \sum_{j=1}^n\partial^2_{x_j}u(x,t)+O(h)\\
&=&\frac1{2n} \Delta u(x,t),
\end{eqnarray*}
which is~\eqref{HC}.

\subsection{The fractional Laplacian and the regional fractional Laplacian}\label{HAJ103456}

Now we consider an open set~$\Omega\subseteq\R^n$
and a discrete random process in~$h\Z^n$ which can be roughly speaking described
in this way. The space parameter~$h>0$  is linked to the time step
\begin{equation}\label{PROh} \tau:=h^{2s}.\end{equation}
A particle starts its journey from a given point~$h\bar k\in\Omega$ of the lattice~$h\Z^n$ and, at each time
step~$\tau$, it can reach any other point of the lattice~$hk$, with~$k\ne\bar k$,
with probability
\begin{equation}\label{PROC} 
P_h(\bar k,k):=
\frac{\chi_\Omega(h\bar k)\,\chi_\Omega(hk)
}{C\,|k-\bar k|^{n+2s}},\end{equation}
then the process continues following the same law.
Notice that the above probability density
does not allow the process to leave the domain~$\Omega$,
since~$P_h$ vanishes in the complement of~$\Omega$
(in jargon, this process is called ``censored'').

In~\eqref{PROC}, the constant~$C>0$ is needed to
normalize to total probability and is defined by 
$$ C:=
\sum_{{k\in\Z^n\setminus\{0\}}}\frac{
1}{|k|^{n+2s}}.$$
We let
\begin{equation}
\label{ONE} c_h(\bar k):=
\sum_{k\in\Z^n\setminus\{0\}} P_h(\bar k,k)=
\sum_{k\in\Z^n\setminus\{0\}} P_h(k,\bar k)
\end{equation}
and
$$p_k(\bar k):=1-c_h(\bar k).$$
Notice that, 
for any~$\bar k\in\Z^n$, it holds that
\begin{equation}\label{36}
c_h(\bar k)\le
\sum_{k\in\Z^n\setminus\{\bar k\}}
\frac{1
}{C\,|k-\bar k|^{n+2s}}=1,\end{equation}
hence, for a fixed~$h>0$ and~$\bar k\in\Z^n$, this
aggregate probability does not equal to~$1$:
this means that there is a remaining probability~$p_h(\bar k)\ge0$
for which 
the particle does not move
(in principle, such
probability is small when so is~$h$,
but, for a bounded domain~$\Omega$,
it is not negligible with respect to the time
step, hence it must be taken into account in the analysis of the process in the general setting
that we present here).

We define $u(x,t)$ to be the probability density for the particle to lie at the point~$x\in\Omega\cap(h\Z^n)$
at time~$t\in\tau\N$.
We show that, for small space and time scale, the function~$u$ is well described by
the evolution of the nonlocal heat equation
\begin{equation}\label{HEAT}
\partial_t u(x,t) = -c\,(-\Delta)^s_\Omega u(x,t)
\qquad{\mbox{ in }}\Omega,\end{equation}
for some normalization constant~$c>0$.
To check this, up to a translation, we suppose that~$x=0\in\Omega$
and we set~$c_h:=c_h(0)$ and~$p_h:=p_h(0)$.
We observe that
the probability of being at~$0$ at time~$t+\tau$ is the sum of
the probabilities of being somewhere else, say at~$hk\in h\Z^n$, at time~$t$, times
the probability of jumping from~$hk$ to the origin,
plus the probability of staying put: that is
\begin{eqnarray*}
u(0,t+\tau)&=&\sum_{k\in\Z^n\setminus\{0\}} u(hk,t)\,P_h(k,0)+u(0,t)\,p_h\\&=&
\sum_{k\in\Z^n\setminus\{0\}} u(hk,t)\,P_h(k,0)+
(1-c_h)\,u(0,t).
\end{eqnarray*}
Thus, recalling~\eqref{ONE},
\begin{eqnarray*}
u(0,t+\tau)- u(0,t)&=&
\sum_{k\in\Z^n\setminus\{0\}} u(hk,t)\,P_h(k,0)-c_h\,u(0,t)
\\&=&
\sum_{k\in\Z^n\setminus\{0\}} \Big( u(hk,t)-u(0,t)\Big)\,P_h(k,0)\\&=&
\sum_{k\in\Z^n\setminus\{0\}} \Big( u(hk,t)-u(0,t)\Big)\,
\frac{\chi_\Omega(hk)
}{C \,|k|^{n+2s}}\\&=& \frac{h^{n+2s}}{C}
\sum_{k\in\Z^n\setminus\{0\}} \Big( u(hk,t)-u(0,t)\Big)\,
\frac{\chi_\Omega(hk)
}{|hk|^{n+2s}}.
\end{eqnarray*}
So, we divide by~$\tau$ and, in view of~\eqref{PROh}, we find that
\[ C\;\frac{u(0,t+\tau)- u(0,t)}{\tau}=
h^{n}\;
\sum_{k\in\Z^n\setminus\{0\}} \Big( u(hk,t)-u(0,t)\Big)\,
\frac{\chi_\Omega(hk)
}{|hk|^{n+2s}}.
\]
We write this identity changing~$k$ to~$-k$ and we sum up: in this way, we obtain that
\begin{equation}\label{PR:0}\begin{split}& 2C\;\frac{u(0,t+\tau)-
u(0,t)}{\tau}\\ =\;&
h^{n}\;
\sum_{k\in\Z^n\setminus\{0\}} 
\frac{\big( u(hk,t)-u(0,t)\big)\,\chi_\Omega(hk)+
\big( u(-hk,t)-u(0,t)\big)\,\chi_\Omega(-hk)
}{|hk|^{n+2s}}.\end{split}
\end{equation}
Now, for small~${y}$, if~$u$ is smooth enough,
\begin{eqnarray*}
&& \Big| \big( u(y,t)-u(0,t)\big)\,\chi_\Omega({y})+
\big( u(-{y},t)-u(0,t)\big)\,\chi_\Omega(-{y})\Big|\\
&=&\Big| \big( u({y},t)-u(0,t)\big)+
\big( u(-{y},t)-u(0,t)\big)\Big|\\
&=&\Big| \big( \nabla u(0,t){y}+O(|{y}|^2)\big)+
\big( -\nabla u(0,t){y}+O(|{y}|^2)\big)\Big|\\
&=& O(|{y}|^2)
\end{eqnarray*}
and therefore, if we write
$$ \varphi(y):=\frac{\big( u(y,t)-u(0,t)\big)\,\chi_\Omega(y)+
\big( u(-y,t)-u(0,t)\big)\,\chi_\Omega(-y)
}{|y|^{n+2s}},$$
we (formally) have that
\begin{equation}\label{00:ord}
\varphi(y)=O(|y|^{2-n-2s})
\end{equation}
for small~$|y|$.

Now, we fix~$\delta>0$
and use the Riemann sum representation of an integral
to write (for a bounded Riemann integrable function~$\varphi:\R^n\setminus B_\delta\to\R$),
\begin{equation}\label{00:ord1}
\int_{\R^n\setminus B_\delta}\varphi(y)\,dy= \lim_{h\searrow0}
h^n \sum_{k\in\Z^n} \varphi(hk) \,\chi_{\R^n\setminus B_\delta}(hk)=
\lim_{h\searrow0}
h^n \sum_{{k\in\Z^n}\atop{k\ne0}} \varphi(hk) \,\chi_{\R^n\setminus B_\delta}(hk)
.\end{equation}
If, in addition, \eqref{00:ord} is satisfied, one has that, for small~$\delta$,
$$ \int_{B_\delta}\varphi(y)\,dy=O(\delta^{2-2s}).$$
{F}rom this and~\eqref{00:ord1} we have that
\begin{equation}\label{00:ord2}\begin{split}
\int_{\R^n}\varphi(y)\,dy\;&=\,O(\delta^{2-2s})+
\lim_{h\searrow0}
h^n \sum_{{k\in\Z^n}\atop{k\ne0}} \varphi(hk) \,\chi_{\R^n\setminus B_\delta}(hk)
\\ &=\,O(\delta^{2-2s})+
\lim_{h\searrow0}
h^n \sum_{{k\in\Z^n}\atop{k\ne0}} \varphi(hk) 
+
h^n \sum_{{k\in\Z^n}\atop{k\ne0}} \varphi(hk) \,\big(\chi_{\R^n\setminus B_\delta}(hk)-1\big).
\end{split}\end{equation}
Also, in view of~\eqref{00:ord},
\begin{eqnarray*}
&& \left|
h^n \sum_{{k\in\Z^n}\atop{k\ne0}} \varphi(hk) \,\big(\chi_{\R^n\setminus B_\delta}(hk)-1\big)\right|
=
\left|
h^n \sum_{{k\in\Z^n}\atop{0<h|k|<\delta}} \varphi(hk)\right|\\
&&\qquad\le\const h^n \sum_{{k\in\Z^n}\atop{0<|k|<\delta/h}} |hk|^{2-n-2s}=
\const h^{2-2s} \sum_{{k\in\Z^n}\atop{0<|k|<\delta/h}} \frac{|k|^{1-s}}{|k|^{n+s-1}}
\\&&\qquad\le\const h^{2-2s} \,\left(\frac{\delta}h\right)^{1-s}\,
\sum_{{k\in\Z^n}\atop{1\le|k|<\delta/h}} \frac1{|k|^{n+s-1}}
\le \const h^{2-2s} \,\left(\frac{\delta}h\right)^{1-s}\,\left(\frac{\delta}h\right)^{1-s}=\const\delta^{2-2s}.
\end{eqnarray*}
Hence, \eqref{00:ord2} boils down to
\begin{equation*}
\int_{\R^n}\varphi(y)\,dy=\,O(\delta^{2-2s})+
\lim_{h\searrow0}
h^n \sum_{{k\in\Z^n}\atop{k\ne0}} \varphi(hk) 
\end{equation*}
and so, taking~$\delta$ arbitrarily small,
\begin{equation*}
\int_{\R^n}\varphi(y)\,dy=\,
\lim_{h\searrow0}
h^n \sum_{{k\in\Z^n}\atop{k\ne0}} \varphi(hk) .
\end{equation*}
Therefore, recalling~\eqref{PR:0}, 
\begin{eqnarray*}
2C\partial_t u(0,t)&=& \lim_{h\searrow0} 2C\;
\frac{u(0,t+\tau)- u(0,t)}{\tau}
\\&=&
\lim_{h\searrow0}
h^{n}\;
\sum_{k\in\Z^n\setminus\{0\}} 
\frac{\big( u(hk,t)-u(0,t)\big)\,\chi_\Omega(hk)+
\big( u(-hk,t)-u(0,t)\big)\,\chi_\Omega(-hk)
}{|hk|^{n+2s}}\\
&=&\lim_{h\searrow0}
h^n \sum_{{k\in\Z^n}\atop{k\ne0}} \varphi(hk)\\&=&
\int_{\R^n}\varphi(y)\,dy\\&=&
\int_{\R^n}
\frac{\big( u(y,t)-u(0,t)\big)\,\chi_\Omega(y)+
\big( u(-y,t)-u(0,t)\big)\,\chi_\Omega(-y)
}{|y|^{n+2s}}\\&=& -2(-\Delta)^s_\Omega u(x,0).
\end{eqnarray*}
This confirms~\eqref{HEAT}.

As a final comment, in view of these calculations and those of Section~\ref{SE:HE},
we may compare the classical random walk, which leads to the classical
heat equation, and the long-jump random walk which leads to the nonlocal
heat equation and relate such jumps to an ``infinitely fast''
diffusion, in the light of
the computations of the associated mean squared displacements
(recall~\eqref{LINE} and~\eqref{810ehqwgoy8r3ryfefi34757575757}).

\subsection{The spectral fractional Laplacian}

Now, we briefly discuss a heuristic motivation for the fractional heat equation
run by the spectral fractional Laplacian, that is
\begin{equation}\label{SPHEAT}
\partial_t u(x,t) = -c\,(-\Delta)^s_{D,\Omega} u(x,t)
\qquad{\mbox{ in }}\Omega,\end{equation}
for some normalization constant~$c>0$.
To this end, we consider a bounded and smooth set~$\Omega\subset\R^n$
and we define a random motion of a ``distribution of particles'' in~$\Omega$.
For any~$x\in\Omega$ and~$t\ge0$, the function~$u(x,t)$ denotes the ``number of particles''
present at the point~$x$ at the time~$t$. No particles lie outside~$\Omega$
and we write~$u$ as a suitable superposition of eigenfunctions~$\{\phi_k\}_{k\ge1}$ of the Laplacian
with Dirichlet boundary data
(this is a reasonable assumption, given that such
eigenfunctions provide a basis of~$L^2(\Omega)$, see e.g. page~335
in~\cite{MR1625845}). In this way, we write
$$ u(x,t)=\sum_{k=1}^{+\infty} u_k(t)\,\phi_k(x).$$
Namely, in the notation in~\eqref{EXPA},
the evolution of the particle distribution~$u$ is 
defined on each spectral component~$u_k$ and it is
taken to follow
a ``classical'' random walk, but the space/time scale is supposed to depend on~$k$ as well:
namely, spectral components relative to high frequencies will move slower
than the ones relative to low frequencies (namely, the time step is taken to be longer if
the frequency is higher).

More precisely, for any~$k\in\N$, we suppose that each of the~$u_k$ particles
of the $k$th spectral component undergo a classical
random walk in a lattice~$h_k\Z^d$, as described in Section~\ref{SE:HE},
but with time step
\begin{equation}\label{TAUK}\tau_k:=\lambda_k^{1-s}\,h_k^2.\end{equation}
We suppose that~$h_k$ and~$\tau_k$ are ``small space and time increments''.
Namely, after a time step~$\tau_k$, each of these $u_k(t)\,\phi_k(x)$ particles
will move, with equal probability~$\frac1{2n}$, to one of the points~$x\pm h_k e_1$, $\dots$, $x\pm h_k e_n$
(for simplicity, we are imaging here~$u_k$ to be positive; the case of negative~$u_k$
represents a ``lack of particles'', which is supposed to diffuse with the same law).
Hence, the number of particles at time~$t+\tau_k$ which correspond to the $k$th frequency of
the spectrum and which lie at the point~$x\in\Omega$ is equal to the sum
of the number of the particles at time~$t$ which lie somewhere else times the probability
of jumping to~$x$ in this time step, that is, in formula,
\begin{equation} \label{20394u}
u_k(t+\tau_k) \,\phi_k(x)=\frac1{2n} \sum_{j=1}^n u_k(t)\, \Big( \phi_k(x+h_ke_j)+\phi_k(x-h_ke_j)\Big).\end{equation}
Moreover,
\begin{eqnarray*}
&& \phi_k(x+h_ke_j)+\phi_k(x-h_ke_j)-2\phi_k(x)\\
&=&
\left( \phi_k(x)+h_k\nabla\phi_k(x)\cdot e_j
+\frac{h_k^2\,D^2\phi_k(x) e_j\cdot e_j}{2}\right)+
\left( \phi_k(x)-h_k\nabla\phi_k(x)\cdot e_j
+\frac{h_k^2\,D^2\phi_k(x) e_j\cdot e_j}{2}\right)\\&&\qquad
-2\phi_k(x)+O(h_k^3)\\&=&
h_k^2\,\partial^2_{x_j}\phi_k(x)+O(h_k^3).
\end{eqnarray*}
Consequently, from this and~\eqref{20394u},
\begin{eqnarray*}
\Big(u_k(t+\tau_k)-u_k(t)\Big) \,\phi_k(x)
&=&\frac1{2n} \sum_{j=1}^n u_k(t)\, \Big( \phi_k(x+h_ke_j)+\phi_k(x-h_ke_j)-2\phi_k(x)\Big)\\
&=&\frac{h_k^2}{2n}\,u_k(t)\,\Delta \phi_k(x)+O(h_k^3)\\&=&-\frac{\lambda_k\,h_k^2}{2n}\,u_k(t)\,\phi_k(x)+O(h_k^3).
\end{eqnarray*}
Hence, with a formal computation, dividing by~$\tau_k$, using~\eqref{TAUK}
and sending~$h_k$, $\tau_k\searrow0$ (for a fixed~$k$), we obtain
$$ \partial_t u_k(t)=\lim_{\tau_k\searrow0}\frac{
\big(u_k(t+\tau_k)-u_k(t)\big) \,\phi_k(x)}{\tau_k}=
\lim_{h_k\searrow0}
-\frac{\lambda_k^s}{2n}\,u_k(t)\,\phi_k(x)+O(h_k)=
-\frac{\lambda_k^s}{2n}\,u_k(t)\,\phi_k(x).$$
Hence, from~\eqref{EXPA} (and neglecting converge issues
in~$k$), we have
\begin{equation*}
\partial_t u(x,t)=\sum_{k=0}^{+\infty} \partial_t u_k(t)\, \phi_k(x)
=
-\sum_{k=0}^{+\infty}\frac{\lambda_k^s}{2n}\,u_k(t)\,\phi_k(x)
,\end{equation*}
that is~\eqref{SPHEAT}.

\subsection{Fractional time derivatives}\label{90woo22o3jejf}

\begin{figure}
    \centering
    \includegraphics[width=11cm]{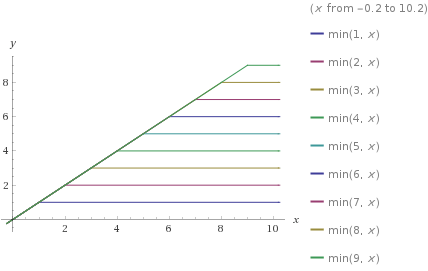}
    \caption{\it {{The motions~$u_k$ described in Section~\ref{90woo22o3jejf}
when the velocity field~$f$ is constant.}}}
    \label{STOP}
\end{figure}

We consider a model in which a bunch of people
is supposed to move along the real line
(say, starting at the origin)
with some given velocity~$f$, which depends on time. We consider the case in which
the environment surrounding the moving
people is ``tricky'', and some of them risk to get stuck
for some time, and they are able
to ``exit the trap'' only by overcoming their
past velocity. Concretely, we fix a function~$\varphi:[0,+\infty)\to[0,+\infty)$
with
\begin{equation} \label{CPHI}
C_\varphi:=\sum_{k=1}^{+\infty} \varphi(k)<+\infty.\end{equation}
Then we define
$$ p_k:=\frac{\varphi(k)}{C_\varphi}$$
and we notice that
$$ \sum_{k=1}^{+\infty} p_k=\frac{1}{C_\varphi}\,\sum_{k=1}^{+\infty} \varphi(k)=1.$$
Then, we denote by~$u(t)$ the position of the ``generic person'' at time~$t$, with~$u(0)=0$.
We suppose that some people, say a proportion~$p_1$
of the total population, 
move with the prescribed velocity for a unit of time,
after which
their velocity is the difference between the prescribed
velocity at that time and the one at the preceding time
with respect to the time unit. 
In formulas, this says
that there is a proportion~$p_1$ of the total people who travels 
with velocity
\[
\dot u_1(t) := \left\{
\begin{matrix}
f(t) & {\mbox{ if }} t\in [0,1],\\
f(t)-f(t-1) 
& {\mbox{ if }} t> 1.
\end{matrix}
\right. \]
After integrating, we thus obtain that
there is a proportion~$p_1$ of the total people whose position
is described by the function
\begin{eqnarray*}
u_1(t) &=& \left\{
\begin{matrix}\displaystyle\int_0^t f(\vartheta)\,d\vartheta & {\mbox{ if }} t\in [0,1],\\
\displaystyle\int_0^1 f(\vartheta)\,d\vartheta+\int_1^t\big(
f(\vartheta)-f(\vartheta-1)\big)\,d\vartheta 
& {\mbox{ if }} t> 1,
\end{matrix}
\right. \\
&=& \left\{
\begin{matrix}\displaystyle\int_0^t f(\vartheta)\,d\vartheta & {\mbox{ if }} t\in [0,1],\\
\displaystyle\int_0^t f(\vartheta)\,d\vartheta-
\int_1^t f(\vartheta-1)\,d\vartheta 
& {\mbox{ if }} t> 1,
\end{matrix}
\right.
\\
&=& \left\{
\begin{matrix}\displaystyle\int_0^t f(\vartheta)\,d\vartheta & {\mbox{ if }} t\in [0,1],\\
\displaystyle\int_0^t f(\vartheta)\,d\vartheta-
\int_0^{t-1} f(\vartheta)\,d\vartheta 
& {\mbox{ if }} t> 1,
\end{matrix}
\right.
\\ &=& \int_{(t-1)_+}^t f(\vartheta)\,d\vartheta.
\end{eqnarray*}
For instance, if~$f$ is constant, then
the position~$u_1$ grows linearly for a unit of time and then remains put
(this would correspond to consider ``stopping times''
in the motion, see Figure~\ref{STOP}).

Similarly, a proportion~$p_2$
of the total population 
evolves with prescribed velocity~$f$ for two units of time,
after which
its velocity becomes the difference between the prescribed
velocity at that time and the one at the preceding time
with respect to two time units, namely
\[
\dot u_2(t) := \left\{
\begin{matrix}
f(t) & {\mbox{ if }} t\in [0,2],\\
f(t)-f(t-2) 
& {\mbox{ if }} t> 2.
\end{matrix}
\right. \]
In this case, an integration gives that
there is a proportion~$p_2$ of the total people whose position
is described by the function
\begin{eqnarray*}
u_2(t) &=& \left\{
\begin{matrix}\displaystyle\int_0^t f(\vartheta)\,d\vartheta & {\mbox{ if }} t\in [0,2],\\
\displaystyle\int_0^2 f(\vartheta)\,d\vartheta+\int_2^t\big(
f(\vartheta)-f(\vartheta-2)\big)\,d\vartheta 
& {\mbox{ if }} t> 2,
\end{matrix}
\right. \\
&=& \left\{
\begin{matrix}\displaystyle\int_0^t f(\vartheta)\,d\vartheta & {\mbox{ if }} t\in [0,2],\\
\displaystyle\int_0^t f(\vartheta)\,d\vartheta-
\int_2^t f(\vartheta-2)\,d\vartheta 
& {\mbox{ if }} t> 2.
\end{matrix}
\right.
\\
&=& \left\{
\begin{matrix}\displaystyle\int_0^t f(\vartheta)\,d\vartheta & {\mbox{ if }} t\in [0,2],\\
\displaystyle\int_0^t f(\vartheta)\,d\vartheta-
\int_0^{t-2} f(\vartheta)\,d\vartheta 
& {\mbox{ if }} t> 2.
\end{matrix}
\right.
\\ &=& \int_{(t-2)_+}^t f(\vartheta)\,d\vartheta.
\end{eqnarray*}
Repeating this argument, we suppose that for each~$k\in\N$
we have a proportion~$p_k$ of the people that move initially
with the prescribed velocity~$f$, but,
after~$k$ units of time,
get their velocity
changed into the difference of the actual velocity field
and that of $k$ units of time before (which is indeed a ``memory effect'').
In this way, we have that a proportion~$p_k$ of the total population
moves with law of motion given by
\begin{eqnarray*}
u_k(t) = \int_{(t-k)_+}^t f(\vartheta)\,d\vartheta.\end{eqnarray*}
The average position of the moving population is then given by
\begin{equation}\label{02345:1234523ws3tgv46777124345}
u(t):= \sum_{k=1}^{+\infty} p_k\,u_k(t)=\frac{1}{C_\varphi}\,
\sum_{k=1}^{+\infty} \varphi(k)\, \int_{(t-k)_+}^t f(\vartheta)\,d\vartheta.\end{equation}
We now specialize the computation above for the case
$$ \varphi(x):=x^{s-2},$$
with~$s\in(0,1)$. Notice that the quantity in~\eqref{CPHI}
is finite in this case, and we can denote it simply by~$C_s$.
In addition, we will consider long time asymptotics in~$t$
and introduce a small time increment~$h$ which
is inversely proportional to~$t$, namely
\begin{equation*}
h:=\frac1t.\end{equation*}
In this way, recalling that the motion was supposed to
start at the origin (i.e., $u(0)=0$)
and using the substitution~$\eta:=\vartheta/t$, we can write~\eqref{02345:1234523ws3tgv46777124345}
as
\begin{equation}\label{02345:1234523ws3tgv46777124345:x}
\begin{split}
u(t)-u(0)\, &=\frac{1}{C_s}\,
\sum_{k=1}^{+\infty} k^{s-2}\, \int_{(t-k)_+}^t f(\vartheta)\,d\vartheta\\
&=\frac{t^s \, h}{C_s}\,
\sum_{k=1}^{+\infty} (hk)^{s-2}\, \int_{(1-kh)_+}^1 f(t\eta)\,d\eta\\
&\simeq\frac{t^s}{C_s}\, 
\int_0^{+\infty}\left[
\lambda^{s-2}\, \int_{(1-\lambda)_+}^1 f(t\eta)\,d\eta\right]\,d\lambda
,\end{split}\end{equation}
where we have recognized a Riemann sum in the last line.

We also point out that the conditions
$$ \lambda\in(0,+\infty) \;{\mbox{ and }}\;
0<\xi<\min\{1,\lambda\}$$
are equivalent to
$$ 0<\xi< 1\;{\mbox{ and }}\;\lambda\in(\xi,+\infty),$$
and, furthermore,
$$ 1-(1-\lambda)_+=
1-\max\{0,\,1-\lambda\}=
\min\{1-0,\,1-(1-\lambda)\}=
\min\{1,\lambda\}.$$
Therefore we use the substitution~$\xi:=1-\eta$
and we exchange the order of integrations, to deduce
from~\eqref{02345:1234523ws3tgv46777124345:x} that
\begin{equation*}
\begin{split}
u(t)-u(0)\, &=\frac{t^s}{C_s}\,\int_0^{+\infty}\left[
\int_0^{\min\{1,\lambda\}} \lambda^{s-2}\, f(t-t\xi)\,d\xi\right]\,d\lambda\\
&=\frac{t^s}{C_s}\,
\int_0^{1}\left[
\int_{\xi}^{+\infty} \lambda^{s-2}\, f(t-t\xi)\,d\lambda\right]\,d\xi
\\&=\frac{t^s}{C_s\,(1-s)}\,
\int_0^{1}\xi^{s-1}\, f(t-t\xi)\,d\xi
.\end{split}\end{equation*}
The substitution~$\tau:=t\xi$ then gives
$$ u(t)-u(0)=
\frac{1}{C_s\,(1-s)}\,
\int_0^{t} \tau^{s-1}\, f(t-\tau)\,d\tau,$$
which, comparing with~\eqref{LAP:L2}
and possibly redefining constants, gives that~$
\partial^s_{C,t} u=f$.

Of course, one can also take into account the case in which the velocity
field~$f$ is induced by a classical diffusion in space, i.e.~$f=\Delta u$,
and in this case one obtains the time fractional diffusive equation~$\partial^s_{C,t} u=\Delta u$.

\subsection{Fractional time diffusion arising from
heterogeneous media}

\begin{figure}
    \centering
    \includegraphics[width=10cm]{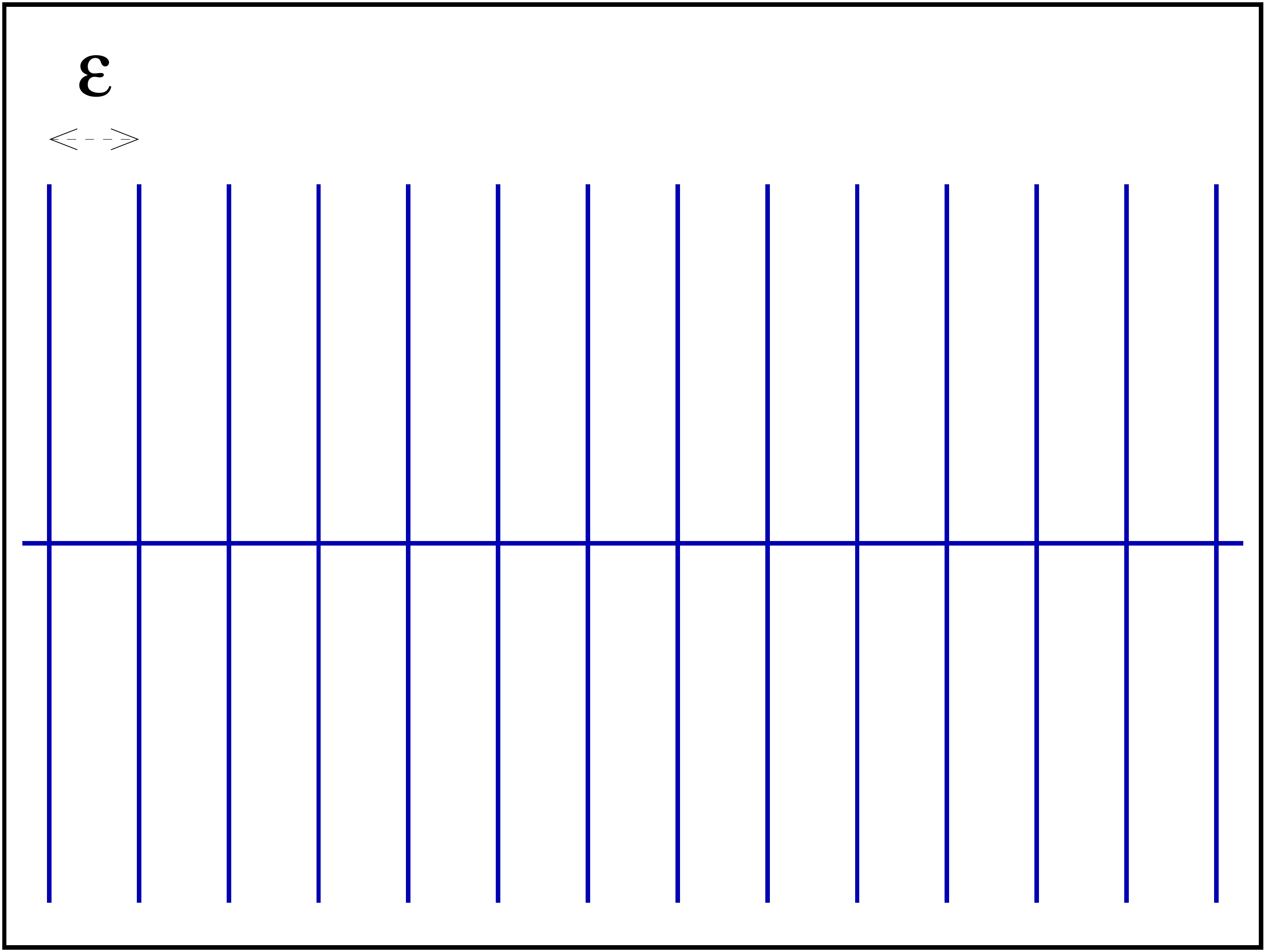}
    \caption{\it {{The comb~${\mathcal{C}}_\varepsilon$.}}}
    \label{PcombO}
\end{figure}

A very interesting phenomenon
to observe is that the geometry of the diffusion medium can naturally
transform classical diffusion into an anomalous one.
This feature can be very well understood by an elegant model,
introduced in~\cite{comb}
(see also~\cite{MR3856678} and the references therein
for an exhaustive account of the research in this direction)
consisting in random walks on a ``comb'', that we briefly
reproduce here for the facility of the reader.
Given~$\varepsilon>0$, the comb may be considered as a transmission medium
that is the union of a ``backbone''~${\mathcal{B}}:=\R\times\{0\}$
with the ``fingers''~${\mathcal{P}}_k:=\{\varepsilon k\}\times\R$, namely
$$ {\mathcal{C}}_\varepsilon:={\mathcal{B}}\cup
\left( \bigcup_{k\in\Z} {\mathcal{P}}_k\right),$$
see Figure~\ref{PcombO}.

We suppose that a particle experiences a random walk on the comb,
starting at the origin,
with some given
horizontal and vertical speeds.
In the limit, this random walk can be modeled by the diffusive equation
along the comb~${\mathcal{C}}_\varepsilon$
\begin{equation}\label{COMB0}
\left\{
\begin{matrix}
u_t = d_1 \delta_0(y)\,u_{xx}+d_2\,
\varepsilon\,\displaystyle \sum_{k\in\Z}\delta_0(\varepsilon k)u_{yy},\\
u(x,y,0)=\delta_0(x)\,\delta_0(y),
\end{matrix}
\right.
\end{equation}
with~$d_1$, $d_2>0$.
The case~$d_1=d_2$ corresponds to equal horizontal and vertical speeds
of the random walk (and this case is already quite interesting).
Also, in the limit as~$\varepsilon\searrow0$, we can consider the Riemann
sum approximation
$$ \int_\R f(y)\,dy\simeq \varepsilon\,\displaystyle \sum_{k\in\Z}f(\varepsilon k),$$
and~${\mathcal{C}}_\varepsilon$ tends to cover the whole of~$\R^2$
when~$\varepsilon$ gets small.
Accordingly, at least at a formal level, as the fingers of the comb become thicker
and thicker, we can think that
$$ 1=
\int_\R \delta_0(y)\,dy\simeq \varepsilon\,\displaystyle \sum_{k\in\Z}\delta_0(\varepsilon k),$$
and reduce~\eqref{COMB0} to
the diffusive equation
in~$\R^2$ given by
\begin{equation}\label{COMB1}
\left\{
\begin{matrix}
u_t = d_1 \delta_0(y)\,u_{xx}+d_2 u_{yy},\\
u(x,y,0)=\delta_0(x)\,\delta_0(y).
\end{matrix}
\right.
\end{equation}
The very interesting feature of~\eqref{COMB1}
is that it naturally induces a fractional time  
diffusion along the backbone.
The quantity that experiences this fractional diffusion
is the total diffusive mass at a point of the backbone. Namely, one sets
\begin{equation}\label{COMBDE}
U(x,t):=\int_\R u(x,y,t)\,dy,
\end{equation}
and we claim that
\begin{equation}\label{MAG}
\partial^{1/2}_{C,t} \,U(x,t)=\frac{d_1}{2\sqrt{d_2}}\,
\Delta U(x,t)\qquad{\mbox{for all }}\,
(x,t)\in\R\times(0,+\infty).
\end{equation}
Equation~\eqref{MAG} reveals the very relevant phenomenon
that a diffusion governed by the Caputo derivative may naturally
arise from classical diffusion, only in view of the particular geometry
of the domain.

To check~\eqref{MAG},
we first point out that
\begin{equation}\label{COMBDE2}
\hat U(\xi,0):=\int_\R \hat u(\xi,y,0)\,dy=\int_\R\delta_0(y)\,dy=1.
\end{equation}
Then, we observe that, if~$a$, $b\in\C$, and
\begin{equation}\label{abdis}
\begin{split}
&g(y):= b\,e^{-a|y|}\qquad{\mbox{ for any }}y\in\R,\\
&{\mbox{then we have that }}\quad g''(y)=a^2 g(y)-2ab\,\delta_0(y).
\end{split}
\end{equation}
To check this let~$\varphi\in C^\infty_0(\R)$. Then, integrating twice by parts,
\begin{eqnarray*}
&& \frac1b\,\int_\R \big(g(y)\varphi''(y)-a^2 g(y)\varphi(y)\big)\,dy\\
&=& \int_0^{+\infty} e^{-ay}\varphi''(y)\,dy+\int_{-\infty}^0 
e^{ay}\varphi''(y)\,dy
-a^2\int_\R e^{-a|y|}\varphi(y)\,dy\\
&=& a\int_0^{+\infty} e^{-ay}\varphi'(y)\,dy-a\int_{-\infty}^0 e^{ay}\varphi'(y)\,dy
-a^2\int_\R e^{-a|y|}\varphi(y)\,dy\\
&=& -2a\varphi(0)
+a^2\int_0^{+\infty} e^{-ay}\varphi(y)\,dy+a^2\int_{-\infty}^0 e^{ay}
\varphi'(y)\,dy
-a^2\int_\R e^{-a|y|}\varphi(y)\,dy\\
&=& -2a\varphi(0),
\end{eqnarray*}
thus proving~\eqref{abdis}.

We also remark that, in the notation of~\eqref{abdis},
we have that~$\delta_0(y)g(y)=\delta_0(y)g(0)=b\delta_0(y)$,
and so, for every~$c\in\R$,
\begin{equation}\label{gabc}
g''(y)=a^2 g(y)-b(2a+c)\delta_0(y)+c\delta_0(y)g(y).\end{equation}
Now, taking the Fourier Transform of~\eqref{COMB1} in the variable~$x$, using the notation~$\hat u(\xi,y,t)$
for the Fourier Transform of~$u(x,y,t)$, and possibly neglecting normalization constants,
we get
\begin{equation}\label{COMB2}
\left\{
\begin{matrix}
\hat u_t = -d_1 |\xi|^2 \delta_0(y)\,\hat u+d_2 \hat u_{yy},\\
\hat u(\xi,y,0)=\delta_0(y).
\end{matrix}
\right.
\end{equation}
Now, we take the Laplace Transform of~\eqref{COMB2} in the variable~$t$,
using the notation~$w(\xi,y,\omega)$ for the Laplace Transform of~$\hat u(\xi,y,t)$, namely~$
w(\xi,y,\omega):={\mathcal{L}} \hat u(\xi,y,\omega)$. In this way, recalling that
$$ {\mathcal{L}}(\dot f)=\omega{\mathcal{L}}f(\omega)-f(0),$$
and therefore
$$ {\mathcal{L}}(\hat u_t)(\xi,y,\omega)=
\omega{\mathcal{L}}\hat u(\xi,y,\omega)-\hat u(\xi,y,0)
=\omega w(\xi,y,\omega)-\delta_0(y),$$
we deduce from~\eqref{COMB2} that
\begin{equation}\label{AVVAjc}
\omega w-\delta_0(y)
= -d_1 |\xi|^2 \delta_0(y)\,w+d_2 w_{yy}
.\end{equation}
That is, setting
$$ a(\omega):=\left( {\frac{\omega}{d_2}}\right)^{1/2},\qquad
b(\xi,\omega):=\frac{1}{\left(4d_2\omega\right)^{1/2}+d_1\,|\xi|^2}
\qquad{\mbox{and}}\qquad
c(\xi):=\frac{d_1\,|\xi|^2}{d_2},$$
we see that
$$ b(2a+c)=
\frac{2
\left(\displaystyle\frac{\omega}{d_2}\right)^{1/2}+\displaystyle\frac{d_1\,|\xi|^2}{d_2}
}{\left( {4d_2\omega}\right)^{1/2}+d_1\,|\xi|^2}=\frac1{d_2},$$
and hence
we can write~\eqref{AVVAjc} as
$$ w_{yy}=\frac{\omega}{d_2}\, w-\frac{1}{d_2}\delta_0(y)+
\frac{d_1 |\xi|^2}{d_2} \delta_0(y)\,w
=a^2 w-b(2a+c)\delta_0(y)+c\delta_0(y)w.$$
In light of~\eqref{gabc}, we know that this equation is solved by taking~$w=g$, that is
$$ {\mathcal{L}} \hat u(\xi,y,\omega)=w(\xi,y,\omega)=
b(\xi,\omega) e^{-a(\omega)|y|}.$$
As a consequence, by~\eqref{COMBDE},
\begin{eqnarray*}
&&{\mathcal{L}} \hat U(\xi,t)=\int_\R {\mathcal{L}}\hat u(\xi,y,\tau)\,dy
=\int_\R b(\xi,\omega) e^{-a(\omega)|y|}\,dy=
\frac{2b(\xi,\omega)}{a}.
\end{eqnarray*}
This and~\eqref{COMBDE2} give that
\begin{eqnarray*}
\left( \left(4d_2\omega\right)^{1/2}+d_1\,|\xi|^2\right){\mathcal{L}} \hat U(\xi,t)=
\frac{{\mathcal{L}} \hat U(\xi,t)}{b(\xi,\omega)}=\frac2{a}=
\left( {\frac{4d_2}{\omega}}\right)^{1/2}=
2\sqrt{d_2}\,{\omega}^{-1/2}\hat U(\xi,0),
\end{eqnarray*}
that is
\begin{eqnarray*}
\omega^{1/2} {\mathcal{L}} \hat U(\xi,t)-{\omega}^{-1/2}\hat U(\xi,0)=
-\frac{d_1}{2\sqrt{d_2}}\,|\xi|^2\,{\mathcal{L}} \hat U(\xi,t).
\end{eqnarray*}
Transforming back and recalling \eqref{LAP:L},
we obtain~\eqref{MAG}, as desired.

\section{All functions are locally $s$-caloric (up to a small error):
proof of~\eqref{BARU:CALORIC}}\label{CALORIC}

We let~$(x,t)\in\R\times\R$ and
consider the operator~$\LLL:=\partial_t+(-\Delta)^s_x$.
One defines
\begin{equation*}
{\mathcal{V}}:= \big\{ h:\R\times\R\to\R
{\mbox{ s.t. $\LLL h=0$ in some neighborhood of the origin in~$\R^2$}} \big\},\end{equation*}
and for any~$J\in\N$, we define
\begin{equation*} {\mathcal{V}}_J := \Big\{ \big(
\partial^\alpha h(0,0)\big)_{{\alpha=(\alpha_x,\alpha_t)
\in\N\times\N}\atop{\alpha_x+\alpha_t\in[0,J]}}
{\mbox{ with }}h\in {\mathcal{V}}\Big\}.\end{equation*}
Notice that~${\mathcal{V}}_J$ is a linear subspace of~$\R^{N+1}$, for some~$N\in\N$.
The core of the proof is to establish the maximal span condition
\begin{equation}\label{TUTT:P}
{\mathcal{V}}_J=\R^{N+1}.
\end{equation}
To this end, we argue for a contradiction
and we suppose that~${\mathcal{V}}_J$ is a linear subspace strictly
smaller than~$\R^{N+1}$: hence, there exists
\begin{equation}\label{TUTT2:P}\nu=(\nu_0,\dots,\nu_N)\in S^{N}\end{equation}
such that
\begin{equation}\label{TUTT3:P}
{\mathcal{V}}_J\subseteq \left\{ X=(X_0,\dots,X_N)\in\R^{J+1} {\mbox{ s.t. }}\nu\cdot X=0
\right\} .
\end{equation}
One considers~$\phi$ to be the first eigenfunctions of~$(-\Delta)^s$ in~$(-1,1)$ with Dirichlet data,
normalized to have unit norm in~$L^2(\R)$. Accordingly,
$$ \left\{
\begin{matrix}
(-\Delta)^s \phi(x)=\lambda\,\phi(x) & {\mbox{ for any }}x\in(-1,1),\cr
\phi(x)=0 &{\mbox{ for any $x$ outside $(-1,1)$,}}
\end{matrix}
\right.$$
for some~$\lambda>0$.

In view of the boundary properties discussed in 
Difference~\ref{009}, one can prove that
\begin{equation} \label{phi e}
\partial^\ell\phi(-1+\delta) = 
\const \delta^{s-\ell}(1+o(1)),\end{equation}
with~$o(1)$ infinitesimal as~$\delta\searrow0$.
So, fixed~$\e$, ${\tau}>0$, we define
$$ h_{\e,{\tau}}(x,t):= e^{-{\tau}\, t} \,\phi\left( -1+\e +\frac{{\tau}^{\frac1{2s}}\,x}{\lambda^{\frac1{2s}} }\right).$$
This function is smooth for
any~$x$ in a small neighborhood of the origin and any~$t\in\R$,
and, in this domain,
\begin{eqnarray*}
\LLL h_{\e,{\tau}}(x,t)&=&
\partial_t\left(e^{-{\tau}\, t} \,\phi\left( 
-1+\e +\frac{{\tau}^{\frac1{2s}} \,x}{\lambda^{\frac1{2s}} }
\right)\right)
+(-\Delta)^s_x\left( e^{-{\tau}\, t} \,\phi\left( -1+ \frac{{\tau}^{\frac1{2s}}
\,x}{\lambda^{\frac1{2s}} }\right)\right)\\&=&
-{\tau}\, e^{-{\tau}\, t} \,\phi\left(
-1+\e +\frac{{\tau}^{\frac1{2s}} \,x}{\lambda^{\frac1{2s}} }
\right)
+\frac{{\tau}\,e^{-\tau\,t}}{{\lambda}} (-\Delta)^s\phi
\left( -1+ \e +\frac{{\tau}^{\frac1{2s}}\,x}{\lambda^{\frac1{2s}} }\right)
\\&=&
-{\tau}\,e^{-{\tau}\, t} \,\phi\left( -1+\e + 
\frac{{\tau}^{\frac1{2s}} \,x}{\lambda^{\frac1{2s}} }\right)
+{\tau}\, e^{-{\tau}\, t} \,\phi\left( -1+ \e +\frac{{\tau}^{\frac1{2s}}\,x}{\lambda^{\frac1{2s}} }\right)
\\ &=&0.
\end{eqnarray*}
This says that~$h_{\e,{\tau}}\in{\mathcal{V}}$ and therefore
$$\big(
\partial^\alpha h_{\e,{\tau}}0,0)\big)_{{\alpha=(\alpha_x,\alpha_t)
\in\N\times\N}\atop{\alpha_x+\alpha_t\in[0,J]}}
\in {\mathcal{V}}_J.$$
This, together with~\eqref{TUTT3:P}, implies that, for any fixed and positive~$\tau$ and~$y$,
\begin{equation*}
\begin{split}&
0 = \sum_{{\alpha=(\alpha_x,\alpha_t)
\in\N\times\N}\atop{\alpha_x+\alpha_t\in[0,J]}} \nu_\alpha\,\partial^\alpha h_{\e,{\tau}}(0,0)
=\sum_{{(\alpha_x,\alpha_t)
\in\N\times\N}\atop{\alpha_x+\alpha_t\in[0,J]}} \nu_{(\alpha_x,\alpha_t)}
\,\partial^{\alpha_t}_t \partial^{\alpha_x}_x h_{\e,{\tau}}(0,0)
\\ &\qquad=\left.\sum_{{(\alpha_x,\alpha_t)
\in\N\times\N}\atop{\alpha_x+\alpha_t\in[0,J]}} \nu_{(\alpha_x,\alpha_t)}
(-{\tau})^{\alpha_t} \,
\left( \frac{{\tau}^{\frac1{2s}}}{\lambda^{\frac1{2s}} }\right)^{\alpha_x}\,
e^{-{\tau}\, t} \;
\partial^{\alpha_x} 
\phi\left( -1+\e +\frac{{\tau}^{\frac1{2s}}\,x}{\lambda^{\frac1{2s}} }\right)\right|_{(x,t)=(0,0)}\\
&\qquad=\sum_{{(\alpha_x,\alpha_t)
\in\N\times\N}\atop{\alpha_x+\alpha_t\in[0,J]}} \nu_{(\alpha_x,\alpha_t)}
\frac{(-1)^{\alpha_t}}{\lambda^{\frac{\alpha_x}{2s}} }\,{\tau}^{\alpha_t+\frac{\alpha_x}{2s}}\;
\partial^{\alpha_x} 
\phi\left( -1+\e \right).\end{split}\end{equation*}
Hence, fixed $\tau>0$, this identity and~\eqref{phi e} yield that
\begin{equation}\label{ZZER:1} 0=\sum_{{(\alpha_x,\alpha_t)
\in\N\times\N}\atop{\alpha_x+\alpha_t\in[0,J]}} \nu_{(\alpha_x,\alpha_t)}
\frac{(-1)^{\alpha_t}}{\lambda^{\frac{\alpha_x}{2s}} }\,{\tau}^{\alpha_t+\frac{\alpha_x}{2s}}\;
\e^{s-\alpha_x}
(1+o(1)), \end{equation}
with~$o(1)$ infinitesimal as~$\e\searrow0$.

We now take~$\bar\alpha_x$ be the largest integer~$\alpha_x$
for which there exists an integer~$\alpha_t$ such that~$\bar\alpha_x+\alpha_t\in[0,J]$
and~$\nu_{(\bar\alpha_x,\alpha_t)}\ne0$. Notice that this definition is well-posed, since
not all the~$\nu_{(\alpha_x,\alpha_t)}$
can vanish, due to~\eqref{TUTT2:P}. Then, \eqref{ZZER:1} becomes
\begin{equation}\label{ZZER:2} 0=\sum_{{{(\alpha_x,\alpha_t)
\in\N\times\N}\atop{\alpha_x+\alpha_t\in[0,J]}}\atop{\alpha_x\le\bar\alpha_x}} \nu_{(\alpha_x,\alpha_t)}
\frac{(-1)^{\alpha_t}}{\lambda^{\frac{\alpha_x}{2s}} }\,{\tau}^{\alpha_t+\frac{\alpha_x}{2s}}\;
\e^{s-\alpha_x}
(1+o(1)), \end{equation}
since the other coefficients vanish by definition of~$\bar\alpha_x$.

Thus, we multiply~\eqref{ZZER:2} by~$\e^{\bar\alpha_x-s}
{\tau}^{-\frac{\bar\alpha_x}{2s}}$ and we take the limit as~$\e\searrow0$:
in this way, we obtain that
\begin{eqnarray*}
0&=&\lim_{\e\searrow0}\sum_{{{(\alpha_x,\alpha_t)
\in\N\times\N}\atop{\alpha_x+\alpha_t\in[0,J]}}\atop{\alpha_x\le\bar\alpha_x}} \nu_{(\alpha_x,\alpha_t)}
\frac{(-1)^{\alpha_t}}{\lambda^{\frac{\alpha_x}{2s}} }\,{\tau}^{\alpha_t
+\frac{\alpha_x}{2s}
-\frac{\bar\alpha_x}{2s}
}\;
\e^{\bar\alpha_x-\alpha_x}
(1+o(1))
\\&=&\sum_{{{\alpha_t
\in\N}\atop{\bar\alpha_x+\alpha_t\in[0,J]}}} \nu_{(\bar\alpha_x,\alpha_t)}
\frac{(-1)^{\alpha_t}}{\lambda^{\frac{\bar\alpha_x}{2s}} }\,{\tau}^{\alpha_t}.\end{eqnarray*}
Since this is valid for any~$\tau>0$, 
by the
Identity Principle for Polynomials
we obtain that
$$ \nu_{(\bar\alpha_x,\alpha_t)}
\frac{(-1)^{\alpha_t}}{\lambda^{\frac{\bar\alpha_x}{2s}} }=0,$$
and thus~$\nu_{(\bar\alpha_x,\alpha_t)}=0$,
for any integer~$\alpha_t$ for which~$\bar\alpha_x+\alpha_t\in[0,J]$.
But this is in contradiction with the definition of~$\bar\alpha_x$
and so we have completed the proof of~\eqref{TUTT:P}.

{F}rom this maximal span property, the proof of~\eqref{BARU:CALORIC}
follows by scaling 
(arguing as done, for instance, in~\cite{HB}).


\begin{appendix}
\addcontentsline{toc}{section}{Appendices}
\vspace{9\bigskipamount}

\addtocontents{toc}{\protect\setcounter{tocdepth}{0}}

{\setstretch{2.0}
\begin{quotation}{{\wesa\large
$\langle\langle$The {\underline{longest appendix}} measured~26cm (10.24in) when it was removed from~72-year-old
Safranco August (Croatia) during an autopsy at the Ljudevit Jurak University Department 
of Pathology, Zagreb, Croatia, on~26 August 2006.$\rangle\rangle$}}
\end{quotation}
\begin{flushright}
(Source: {\tt http://www.guinnessworldrecords.com/world-records/largest-appendix-removed})\bigskip
\end{flushright}
}

\section{Confirmation of~\eqref{EX:APP:F}}\label{EX:APP}

We write~$\Delta_x$ to denote the Laplacian in the coordinates~$x\in\R^n$.
In this way, the total Laplacian in the variables~$(x,y)\in\R^n\times(0,+\infty)$
can be written as
\begin{equation}\label{09uyfds0iuyfd9uytdg}\Delta=\Delta_x+\partial^2_y.\end{equation}
Given a (smooth and bounded, in the light of
footnote~\ref{0qeiuru9467890000} on page~\pageref{0qeiuru9467890000})~$u:\R^n\to\R$, we take~$U:=E_u$ be (smooth
and bounded) as in~\eqref{EX:APP:E}.

We also consider the operator
\begin{equation}\label{9wefgu45647832}
Lu(x):=-\partial_y E_u(x,0)\end{equation}
and we take~$V(x,y):=-\partial_y U(x,y)$.
Notice that~$\Delta V=-\partial_y \Delta U=0$ in~$\R^n\times(0,+\infty)$
and~$V(x,0)=Lu(x)$ for any~$x\in\R^n$. In this sense, $V$ is the harmonic extension of~$Lu$
and so we can write~$V=E_{Lu}$ and so, in the notation of~\eqref{9wefgu45647832},
and recalling~\eqref{EX:APP:E} and~\eqref{09uyfds0iuyfd9uytdg},
we have
\begin{eqnarray*}&& L (Lu)(x) = -\partial_y E_{Lu}(x,0) =-\partial_y V(x,0)=
\partial_y^2 U(x,0)\\ &&\qquad=
\Delta U(x,0)-\Delta_x U(x,0)=-\Delta_x U(x,0)=-\Delta u(x).
\end{eqnarray*}
This gives that~$L^2=-\Delta$, which is consistent with~$L=(-\Delta)^{1/2}$,
thanks to~\eqref{QUASDA}.

\section{Proof of \eqref{DECAT:2}}\label{DECAT:3}

Let~$u\in{\mathcal{S}}$. By~\eqref{DECATS}, we can write
\begin{equation}\label{BOUn2}
\sup _{x\in \R^{n}}(1+|x|^{n})\left|u(x)\right|
+
\sup _{x\in \R^{n}}(1+|x|^{n+2 })\left|D^{2}u(x)\right|\le\const.\end{equation}
Fixed~$x\in\R^n$ (with~$|x|$ to be taken large), recalling the notation in~\eqref{2PO},
we consider the map~$y\mapsto\delta_u(x,y)$ and we observe that
\begin{eqnarray*}
&&\delta_u(x,0)=0,\\
&& \nabla_y \delta_u(x,y)=\nabla u(x+y)-\nabla u(x-y),\\{\mbox{and }}
&& D^2_y \delta_u(x,y)=D^2 u(x+y)+D^2 u(x-y).\end{eqnarray*}
Hence, if~$|Y|\le |x|/2$ we have that~$|x\pm Y|\ge |x|-|Y|\ge |x|/2$, and thus
$$ \big| D^2_y \delta_u(x,Y)\big|\le 2\sup_{|\zeta|\ge|x|/2} \big|D^2 u(\zeta)\big|
\le 2\sup_{|\zeta|\ge|x|/2}\frac{(2|\zeta|)^{n+2}\big|D^2 u(\zeta)\big|}{|x|^{n+2}}
\le \frac{\const}{|x|^{n+2}},$$
thanks to~\eqref{BOUn2}.

Therefore, a second order Taylor expansion of~$\delta_u$ in the variable~$y$ gives that,
if~$|y|\le|x|/2$,
\begin{eqnarray*}&&
\big| \delta_u(x,y)\big| \le\sup_{|Y|\le|x|/2}
\left| \delta_u(x,0)+\nabla \delta_u(x,0)\cdot y+\frac{D^2 \delta_u(x,Y)\,y\cdot y}{2}\right|\\
&&\qquad=\sup_{|Y|\le|x|/2}
\left| \frac{D^2 \delta_u(x,Y)\,y\cdot y}{2}\right|\le\frac{\const\,|y|^2}{|x|^{n+2}}.
\end{eqnarray*}
Consequently,
\begin{equation}\label{9765eolkhr56789oytrdcv}
\left| \int_{B_{|x|/2}} \frac{\delta_u(x,y)}{|y|^{n+2s}}\,dy\right|\le
\frac{\const}{|x|^{n+2}}\int_{B_{|x|/2}}\frac{|y|^2}{|y|^{n+2s}}\,dy
\le\frac{\const\,|x|^{2-2s}}{|x|^{n+2}}=\frac{\const}{|x|^{n+2s}}.\end{equation}
Moreover, by~\eqref{BOUn2},
\begin{eqnarray*}
&& \left| \int_{\R^n\setminus B_{|x|/2}} \frac{\delta_u(x,y)}{|y|^{n+2s}}\,dy\right|\\
&\le& \int_{\R^n\setminus B_{|x|/2}} \frac{|u(x+y)|}{|y|^{n+2s}}\,dy
+\int_{\R^n\setminus B_{|x|/2}} \frac{|u(x-y)|}{|y|^{n+2s}}\,dy
+2\int_{\R^n\setminus B_{|x|/2}} \frac{|u(x)|}{|y|^{n+2s}}\,dy\\
&\le& \int_{\R^n\setminus B_{|x|/2}} \frac{|u(x+y)|}{(|x|/2)^{n+2s}}\,dy
+\int_{\R^n\setminus B_{|x|/2}} \frac{|u(x-y)|}{(|x|/2)^{n+2s}}\,dy
+\frac{\const\,|u(x)|}{|x|^{2s}}\\
&\le& \frac{\const}{|x|^{n+2s}} \int_{\R^n} |u(\zeta)|\,d\zeta+
\frac{\const\,|u(x)|}{|x|^{2s}}\\
&\le& \frac{\const}{|x|^{n+2s}}.
\end{eqnarray*}
This and~\eqref{9765eolkhr56789oytrdcv}, recalling~\eqref{2PO}, establish~\eqref{DECAT:2}.

\section{Proof of~\eqref{BOU:REG}}\label{BOU:REG:S}

Let~$M:=\frac1{2n}\,\left(1+\sup_{B_1}|f|\right)$ and~$v(x):= M(1-|x|^2)-u(x)$. 
Notice that~$v=0$ along~$\partial B_1$ and
$$ \Delta v = -2n M -\Delta u \le -M-f\le -M+\sup_{B_1}|f|\le0$$
in~$B_1$. Consequently, $v\ge0$ in~$B_1$, which gives that~$u(x)\le M(1-|x|^2)$.

Arguing similarly, by looking at~$\tilde v(x):= M(1-|x|^2)+u(x)$, one sees that~$-u(x)
\le M(1-|x|^2)$. Accordingly, we have that
$$ |u(x)|\le M(1-|x|^2)\le M(1+|x|)(1-|x|)\le 2M(1-|x|).$$
This proves~\eqref{BOU:REG}.

\section{Proof of \eqref{CERC}}\label{CERC:A}

\begin{figure}
    \centering
    \includegraphics[width=11cm]{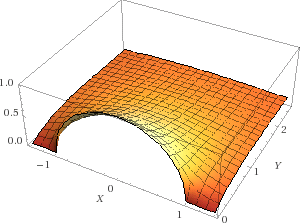}
    \caption{\it {{Harmonic extension in
the halfplane of the function~$\R\ni x\mapsto(1-x^2)_+^{1/2}$.}}}
    \label{MK12}
\end{figure}

The idea of the proof is described in Figure~\ref{MK12}.
The trace of the function in Figure~\ref{MK12} is exactly
the function~$u_{1/2}$ in~\eqref{CERCu12}.
The function plotted in Figure~\ref{MK12}
is the harmonic extension of~$u_{1/2}$ in the halfplane
(like an elastic membrane pinned at the halfcircumference along the trace).
Our objective is to show that the normal derivative of such extended function
along the trace is constant, and so we can
make use of the extension method
in~\eqref{EX:APP:E}
and~\eqref{EX:APP:F} to obtain~\eqref{CERC}.

\begin{figure}
    \centering
    \includegraphics[width=9cm]{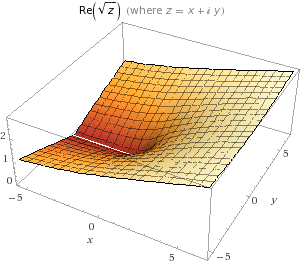}
    \includegraphics[width=9cm]{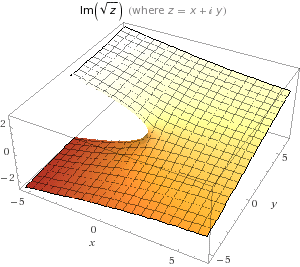}
    \caption{\it {{Real and imaginary part of the complex principal square root.}}}
    \label{RIE}
\end{figure}

In further detail, we use complex coordinates,
identifying~$(x,y)\in\R\times(0,+\infty)$
with~$z:=x+iy\in \C$ with~$\Im(z)>0$.
Also, as customary, we define the principal square root
in the cut complex plane
$$\C_\star:=\{ 
z=re^{{i\varphi }}{\text{ with }}r>0 {\mbox{ and }} -\pi <\varphi < \pi \}$$
by defining, for any~$z=re^{{i\varphi }}\in\C_\star$,
\begin{equation} \label{RAD}
\SURD({z}) := \sqrt{r} \, e^{i \varphi / 2},\end{equation}
see Figure~\ref{RIE} (for typographical convenience, we distinguish between
the complex and the real square root, by using the symbols~${\SURD}(\cdot)$ and~$\sqrt{\cdot}$
respectively).

The principal square root function is
defined using the nonpositive
real axis as a ``branch cut''
and
\begin{equation}\label{INC}
\big( \SURD({z}) \big)^2= r \, e^{i \varphi }=z.
\end{equation}
Moreover,
\begin{eqnarray}\label{9eydiuewgfi3278X197345}
&{\mbox{the function~$\SURD$ is holomorphic in~$\C_\star$}}\\
\label{DERIVATIVE}
& {\mbox{and }}\;\partial_z\SURD({z})=\frac1{2\SURD(z)}.
\end{eqnarray}
To check these facts, we take~$z\in\C_\star$: since~$\C_\star$ is open,
we have that~$z+w\in\C_\star$ for any~$w\in\C\setminus\{0\}$
with small module. Consequently, by~\eqref{INC}, we obtain that
\begin{eqnarray*}
&& w=(z+w)-z=\big( \SURD({z+w}) \big)^2-\big( \SURD({z}) \big)^2\\
&&\qquad=
\Big(
\SURD({z+w})+\SURD({z})\Big)
\Big(
\SURD({z+w})-\SURD({z})\Big).
\end{eqnarray*}
Dividing by~$w$ and taking the limit, we thus find that
\begin{equation}\label{HJA:001234985}\begin{split}
1\;&=\lim_{w\to0}\Big(
\SURD({z+w})+\SURD({z})\Big)\,\frac{
\SURD({z+w})-\SURD({z})}{w}\\
&= 2\SURD({z})\;\lim_{w\to0}\frac{
\SURD({z+w})-\SURD({z})}{w}
\end{split}\end{equation}
Since~$\C_\star\subseteq\C\setminus\{0\}$, we have that~$z\ne0$, and thus~$\SURD(z)\ne0$.
As a result, we can divide~\eqref{HJA:001234985} by~$2\SURD({z})$ and conclude that
$$ \lim_{w\to0}\frac{
\SURD({z+w})-\SURD({z})}{w}=\frac1{ 2\SURD({z}) },$$
which establishes, at the same time, both~\eqref{9eydiuewgfi3278X197345}
and~\eqref{DERIVATIVE}, as desired.

We also remark that
\begin{equation}\label{inclus}
{\mbox{if~$z\in\C$ with~$\Im (z)>0$, then~$1-z^2\in\C_\star$.}}
\end{equation}
To check this, if~$z=x+iy$ with~$y>0$, we observe that
\begin{equation} \label{1men}1-z^2 =1-(x+iy)^2 = 1-x^2+y^2-2ixy.\end{equation}
Hence, if~$1-z^2$ lies on the real axis, we have that~$xy=0$,
and so~$x=0$. Then, the real part of~$1-z^2$ in this case is equal
to~$1+y^2$ which is strictly positive. This
proves~\eqref{inclus}.

Thanks to~\eqref{inclus}, for any~$z\in\C$ with~$\Im (z)>0$ we can define
the function~$\SURD(1-z^2)$.
{F}rom~\eqref{1men}, we can write
\begin{equation*}\begin{split}
&
1-z^2 = r(x,y)\, e^{{i\varphi(x,y) }},\\{\mbox{where }}\;&
r(x,y)=\big(  (1-x^2+y^2)^2+4x^2y^2 \big)^{1/2},\\&
r(x,y)\,\cos\varphi(x,y)=1-x^2+y^2\\
{\mbox{and }}&r(x,y)\,\sin\varphi(x,y)=2xy.
\end{split}\end{equation*}
Notice that
\begin{equation*}
\lim_{y\searrow0} r(x,y)=
\big(  (1-x^2)^2\big)^{1/2}=|1-x^2|.\end{equation*}
As a consequence,
\begin{eqnarray*}&&
|1-x^2|\,\lim_{y\searrow0} \cos\varphi(x,y)=
\lim_{y\searrow0} r(x,y)\,\cos\varphi(x,y)=\lim_{y\searrow0}(1-x^2+y^2)=1-x^2\\
{\mbox{and }}&&
|1-x^2|\,\lim_{y\searrow0} \sin\varphi(x,y)=\lim_{y\searrow0}
r(x,y)\,\sin\varphi(x,y)=\lim_{y\searrow0} 2xy=0.
\end{eqnarray*}
This says that, if~$x^2>1$ then
\begin{eqnarray*}&&
\lim_{y\searrow0} \cos\varphi(x,y)=-1\\
{\mbox{and }}&&
\lim_{y\searrow0} \sin\varphi(x,y)=0,
\end{eqnarray*}
while if~$x^2<1$ then
\begin{eqnarray*}&&
\lim_{y\searrow0} \cos\varphi(x,y)=1\\
{\mbox{and }}&&
\lim_{y\searrow0} \sin\varphi(x,y)=0.
\end{eqnarray*}
On this account, we deduce that
\begin{equation}\label{0iug9uy7678766666} \lim_{y\searrow0} \varphi(x,y)=\left\{
\begin{matrix}
\pi & {\mbox{ if }}x^2>1,\\
0& {\mbox{ if }}x^2<1
\end{matrix}\right.\end{equation}
and therefore, recalling~\eqref{RAD},
\begin{equation}\label{J11KLA}\begin{split}& \lim_{y\searrow0} \SURD(1-z^2)=
\lim_{y\searrow0}\sqrt{r(x,y)} \, e^{i \varphi (x,y)/ 2}
=\left\{
\begin{matrix}
\sqrt{|1-x^2|} \;e^{i\pi/2} & {\mbox{ if }}x^2>1,\\
\sqrt{|1-x^2|} \;e^{i0}& {\mbox{ if }}x^2<1,\\
0 & {\mbox{ if }}x^2=1
\end{matrix}\right.\\&\qquad\qquad=\left\{
\begin{matrix}
i\,\sqrt{|1-x^2|} & {\mbox{ if }}x^2>1,\\
\sqrt{|1-x^2|} & {\mbox{ if }}x^2<1,\\
0 & {\mbox{ if }}x^2=1.
\end{matrix}\right.
\end{split}\end{equation}
This implies that
\begin{equation}\label{9e8dufeiyt4857}
\begin{split}
\lim_{y\searrow0} \Re\Big(\SURD(1-z^2)\Big)\;&=
\left\{
\begin{matrix}
0 & {\mbox{ if }}x^2\geq1,\\
\sqrt{|1-x^2|} & {\mbox{ if }}x^2<1
\end{matrix}\right.\\
&=(1-x^2)^{1/2}_+.
\end{split}\end{equation}
Now we define
$$ z=x+iy\mapsto \Re\Big( \SURD(1-z^2)+i z\Big)=:U_{1/2}(x,y).$$
The function~$U_{1/2}$ is the harmonic extension of~$u_{1/2}$
in the halfplane, as plotted in Figure~\ref{MK12}. Indeed,
from~\eqref{9e8dufeiyt4857},
\begin{eqnarray*} \lim_{y\searrow0}
U_{1/2}(x,y)=
\lim_{y\searrow0}\Re\Big( \SURD(1-z^2)+i x-y\Big)=(1-x^2)^{1/2}_+=u_{1/2}(x).
\end{eqnarray*}
Furthermore, from~\eqref{9eydiuewgfi3278X197345},
we have that~$U_{1/2}$ is the real part of a holomorphic function
in the halfplane and so it is harmonic.

These considerations give that~$U_{1/2}$ solves
the harmonic extension problem
in~\eqref{EX:APP:E}, hence, in the light of~\eqref{EX:APP:F},
\begin{equation}\label{0987098709870987}
\begin{split}& (-\Delta)^{1/2} u_{1/2}(x)
=\lim_{y\searrow0}
-\partial_y U_{1/2}(x,y)=\lim_{y\searrow0}
-\Re\Big(\partial_y \SURD(1-z^2)+i \partial_y z\Big)\\&\qquad\qquad=
\lim_{y\searrow0}
-\Re\Big(\partial_y \SURD(1-z^2)-1\Big)=
1-\lim_{y\searrow0}\Re\Big(\partial_y \SURD(1-z^2)\Big).
\end{split}\end{equation}
Now, recalling~\eqref{DERIVATIVE}, we see that, for any~$x\in(-1,1)$
and small~$y>0$,
\begin{equation}
\label{023943548767676}\partial_y\SURD(1-z^2)=\partial_z\SURD(1-z^2)\;\partial_y z=
\frac1{2\SURD(1-z^2)}\;\partial_z(1-z^2)\;\partial_y (x+iy)=
-\frac{iz}{\SURD(1-z^2)}.
\end{equation}
We stress that 
the latter denominator does not vanish
when~$x\in(-1,1)$
and~$y>0$ is small.
So, using that~$\Re(ZW)=\Re Z\Re W-\Im Z\Im W$ for any~$Z$, $W\in\C$,
we obtain that
\begin{equation}\label{09iu345jhgdfg}\begin{split}&
y=
\Re\big(-i(x+iy)\big)=
\Re(-iz)= \Re\Big(\SURD(1-z^2)\;\partial_y\SURD(1-z^2)\Big)\\
&\qquad=
\Re\Big(\SURD(1-z^2)\Big)\;\Re\Big(\partial_y\SURD(1-z^2)\Big)
-
\Im\Big(\SURD(1-z^2)\Big)\;\Im\Big(\partial_y\SURD(1-z^2)\Big)
.\end{split}\end{equation}
{F}rom~\eqref{J11KLA}, for any~$x\in(-1,1)$
we have that
$$ \lim_{y\searrow0}
\Im\Big(\SURD(1-z^2)\Big)=\Im\Big( \sqrt{|1-x^2|}\Big)=0.$$
This and the fact that~$\partial_y\SURD(1-z^2)$
is bounded (in view of~\eqref{023943548767676}) give that,
for any~$x\in(-1,1)$,
\begin{equation*}
\lim_{y\searrow0}
\Im\Big(\SURD(1-z^2)\Big)\;\Im\Big(\partial_y\SURD(1-z^2)\Big)
=0.
\end{equation*}
This, \eqref{J11KLA} and~\eqref{09iu345jhgdfg}
imply that, for any~$x\in(-1,1)$,
\begin{eqnarray*}&&
0=\lim_{y\searrow0} y=
\lim_{y\searrow0}
\Re\Big(\SURD(1-z^2)\Big)\;\Re\Big(\partial_y\SURD(1-z^2)\Big)
-
\Im\Big(\SURD(1-z^2)\Big)\;\Im\Big(\partial_y\SURD(1-z^2)\Big)
\\&&\qquad=\Re\Big( \sqrt{|1-x^2|}\Big)\;
\lim_{y\searrow0}
\Re\Big(\partial_y\SURD(1-z^2)\Big)+0\\
&&\qquad=\sqrt{|1-x^2|}\;
\lim_{y\searrow0}\Re\Big(\partial_y\SURD(1-z^2)\Big)
\end{eqnarray*}
and therefore
\begin{equation}\label{NU99} \lim_{y\searrow0}\Re\Big(\partial_y\SURD(1-z^2)\Big)=0.\end{equation}
Plugging this information into~\eqref{0987098709870987},
we conclude the proof of~\eqref{CERC}, as desired.

\section{Deducing~\eqref{AHHA} from~\eqref{HHA}
using a space inversion}\label{JSE}

\begin{figure}
    \centering
    \includegraphics[height=5.6cm]{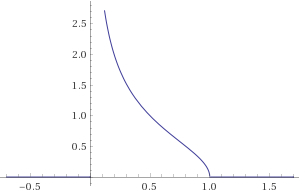} \quad \includegraphics[height=5.6cm]{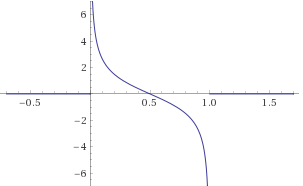}
    \caption{\it {{The functions~$w_{1/2}$ and~$W^\star_{1/2}$.}}}
    \label{F284050050202020TO}
\end{figure}

{F}rom~\eqref{HHA}, up to a translation, we know that
\begin{equation}\label{274835768dsp39057yfh578}
{\mbox{the function~$\R\ni x\mapsto v_s(x):=(x-1)_+^s$ is $s$-harmonic in~$(1,+\infty)$.}}
\end{equation}
We let~$w_s$ be the space inversion of~$v_s$ induced by the Kelvin transform
in the fractional setting, namely
$$ w_s(x):=|x|^{2s-1}\, v_s\left(\frac{x}{|x|^2}\right)
=|x|^{2s-1}\,\left( \frac{x}{|x|^2}-1\right)_+^s=
\left\{\begin{matrix}
x^{s-1}(1-x)^s & {\mbox{ if }}x\in(0,1),\\
0 & {\mbox{ otherwise}}.
\end{matrix}
\right.
.$$
By \eqref{274835768dsp39057yfh578},
see Corollary~2.3 in~\cite{MR2964681},
it follows that~$w_s(x)$ is $s$-harmonic in~$(0,1)$.
Consequently, the function
$$ w^\star_s(x):=w_s(1-x)
=\left\{\begin{matrix}
x^{s}(1-x)^{s-1} & {\mbox{ if }}x\in(0,1),\\
0 & {\mbox{ otherwise}}.
\end{matrix}
\right.
$$
is also $s$-harmonic in~$(0,1)$.
We thereby conclude that the function
\begin{eqnarray*}
&& W^\star_s(x):= w_s(x)-w^\star_s(x)
=\left\{\begin{matrix}
x^{s-1}(1-x)^s-x^{s}(1-x)^{s-1} & {\mbox{ if }}x\in(0,1),
\\0 & {\mbox{ otherwise}}.
\end{matrix}
\right.
\end{eqnarray*}
is also $s$-harmonic in~$(0,1)$. See Figure~\ref{F284050050202020TO} for a picture of~$w_s$
and~$W^\star_s$ when~$s=1/2$.
Let now
$$ U_s(x):=x_+^s(1-x)_+^s
=\left\{\begin{matrix}
x^{s}(1-x)^{s} & {\mbox{ if }}x\in(0,1),
\\0 & {\mbox{ otherwise}}.
\end{matrix}
\right.$$
and notice that~$U_s$ is the primitive of~$s W^\star_s$. Since the latter function is $s$-harmonic
in~$(0,1)$, after an integration we thereby deduce that~$(-\Delta)^s U_s=\const$ in~$(0,1)$.
This and the fact that
$$U_s\left(\frac{x+1}2\right)=2^{-s}\,u_s(x)$$
imply~\eqref{AHHA}.

\section{Proof of~\eqref{BOU:REG:LIP}}\label{BOU:REG:LIP:A}

Fixed~$y\in\R^n\setminus\{0\}$ we let~${\mathcal{R}}^y$ be a rotation which
sends~$\frac{y}{|y|}$ into the vector~$e_1=(1,0,\dots,0)$, that is
\begin{equation}\label{3434} \sum_{k=1}^n {\mathcal{R}}_{ik}^y\, y_k =|y|\delta_{i1},
\end{equation}
for any~$i\in\{1,\dots,n\}$. We also denote by
$$ K(y):=\frac{y}{|y|^2}$$
the so-called Kelvin Transform. We recall that
for any~$i$, $j\in\{1,\dots,n\}$,
$$ \partial_{y_i} K_j(y)=\frac{\delta_{ij}}{|y|^2}-\frac{2y_iy_j}{|y|^4}$$
and so, by~\eqref{3434},
$$ \Big( {\mathcal{R}}^y\;(DK(y))\;({\mathcal{R}}^y)^{-1}\Big)_{ij}=
\sum_{k,h=1}^n {\mathcal{R}}^y_{ik}
\,\partial_{y_k} K_h(y)\,
{\mathcal{R}}^y_{jh}=\frac{\delta_{ij}}{|y|^2}
-\frac{2\delta_{i1}\delta_{j1}}{|y|^2}.$$
This says that~$ {\mathcal{R}}^y\;(DK(y))\;({\mathcal{R}}^y)^{-1}$
is a diagonal\footnote{{F}rom the geometric point of view, one can 
also take radial coordinates,
compute the derivatives of~$K$ along the unit sphere and use scaling.}
matrix, with first entry equal to~$-\frac{1}{|y|^2}$
and the others equal to~$\frac{1}{|y|^2}$.

As a result,
\begin{equation}\label{MESK}
\big|\det(DK(y))\big|=\Big| \det\big(
{\mathcal{R}}^y\;(DK(y))\;({\mathcal{R}}^y)^{-1}\big)\Big|=\frac{1}{|y|^{2n}}.
\end{equation}
The Kelvin Transform is also useful to write the Green function of the ball~$B_1$,
see e.g. formula~(41) on p.~40 
and Theorem~13 on p.~35 of~\cite{MR1625845}. Namely, we take~$n\ge3$
for simplicity, and we write
\begin{eqnarray*}
&& G(x,y):=
\const\left( \frac{1}{|y-x|^{n-2}}-\frac{1}{\big|\,|x|(y-K(x))\big|^{n-2}}\right)\\
&&\qquad\qquad=
\const\left( \frac{1}{|x-y|^{n-2}}-\frac{1}{\big|\,|y|(x-K(y))\big|^{n-2}}\right)=G(y,x)
\end{eqnarray*}
and, for a suitable choice of the constant, for any~$x\in B_1$
we can write the solution of~\eqref{LA:E10}
in the form
$$ u(x)=\int_{B_1}f(y)\,G(x,y)\,dy.$$
see e.g. page~35 in~\cite{MR1625845}.

On this account, we have that, for any~$x\in B_1$,
\begin{eqnarray*}
|\nabla u(x)|&\le& \int_{B_1}|f(y)|\,|\partial_x G(x,y)|\,dy
\\ &\le& \const\sup_{B_1}|f|
\int_{B_1}\left( 
\frac{1}{|x-y|^{n-1}}+\frac{1}{|y|^{n-2}\big|x-K(y)\big|^{n-1}}
\right)\,dy\\
&\le& 
\const\sup_{B_1}|f|
\left( \int_{B_2}
\frac{d\zeta}{|\zeta|^{n-1}} +
\int_{\R^n\setminus B_1}
\frac{d\eta}{|\eta|^{n+2}\big|x-\eta\big|^{n-1}}
\right)
\\&\le& 
\const\sup_{B_1}|f|
\left( 1+
\int_{B_2\setminus B_1}
\frac{d\eta}{|x-\eta|^{n-1}}
+
\int_{\R^n\setminus B_2}
\frac{d\eta}{|\eta|^{n+2}}
\right)
\\&\le&
\const\sup_{B_1}|f|.
\end{eqnarray*}
Notice that here we have used the transformations~$\zeta:=x-y$ and~$\eta:=K(y)$,
exploiting also~\eqref{MESK}.
The claim in~\eqref{BOU:REG:LIP}
is thus established.


\section{Proof of~\eqref{CERC-HARM} and probabilistic insights}\label{CERC-HARM:A}

We give a proof of~\eqref{CERC-HARM} by taking a derivative
of~\eqref{CERC}. To this aim,
we claim\footnote{The difficulty in proving~\eqref{DER91-9dz}
is that the function~$u_{1/2}$ is not differentiable
at~$\pm1$ and the derivative taken inside the integral
might produce a singularity (in fact,
formula~\eqref{DER91-9dz}
exactly says that such derivative can be performed with no harm
inside the integral).
The reader who is already familiar
with the basics of functional analysis can prove~\eqref{DER91-9dz}
by using the theory of
absolutely continuous functions, see e.g. Theorem~8.21
in~\cite{MR0210528}. We provide
here a direct proof, available to everybody.}
that
\begin{equation}\label{DER91-9dz}\begin{split}&
\frac{d}{dx} \int_\R \frac{u_{1/2}(x+y)+u_{1/2}(x-y)-2u_{1/2}(x)}{|y|^{2}}\,dy\\ =\;&
-\int_\R
\frac{(x+y) u_{-1/2}(x+y)+(x-y) u_{-1/2}(x-y)-2x u_{-1/2}(x)}{|y|^{2}}
\,dy.\end{split}
\end{equation}
To this end, we fix~$x\in(-1,1)$ and~$h\in\R$.
We define
$$ \ell_x:=\min\{ |x-1|,\,|x+1|\}>0.$$
In the sequel,
we will take~$|h|$ as small as we wish in order to compute incremental
quotients, hence we can assume that
\begin{equation}\label{P0ow2992924rf}
|h|<\frac{\ell_x}4.\end{equation}
We also define
\begin{equation}\label{9rfyw7rf203948ks}
I_x(h) := \Big\{
y\in\R {\mbox{ s.t. }} \min\{ |(x+y)-1|,\,|(x-y)-1|,\,
|(x+y)+1|,\,|(x-y)+1|
\}\le 2|h|
\Big\}.
\end{equation}
Since~$I_x(h)\subseteq(x-1-2|h|,x-1+2|h|)\cup
(x+1-2|h|,x+1+2|h|)
\cup(1-x-2|h|,1-x+2|h|)\cup(-1-x-2|h|,-1-x+2|h|)$, we have that
\begin{equation}\label{0-er9858585859gigi}
{\mbox{the measure of $I_x$ is less than $\const|h|$.}}
\end{equation}
Furthermore,
\begin{equation}\label{01019238397757548483993}
I_x(h) \subseteq \left\{
y\in\R {\mbox{ s.t. }} |y|\ge \frac{\ell_x}2
\right\}. 
\end{equation}
To check this, let~$y\in I_x(h)$. Then, by~\eqref{9rfyw7rf203948ks},
there exist~$\sigma_{1,x,y}$,
$\sigma_{2,x,y}\in\{-1,1\}$ such that
$$|x+\sigma_{1,x,y} y
+\sigma_{2,x,y}|\le2|h|$$
and therefore
$$ |y|=|\sigma_{1,x,y} y| \ge
|x+\sigma_{2,x,y}| -
|x+\sigma_{1,x,y} y+\sigma_{2,x,y}|\ge\ell_x-2|h|\ge\frac{\ell_x}2,$$
where the last inequality is a consequence of~\eqref{P0ow2992924rf},
and this establishes~\eqref{01019238397757548483993}.

Now, we introduce the following notation for the incremental quotient
$$ Q_h(x,y):=\frac{\big(
u_{1/2}(x+y+h)+u_{1/2}(x-y+h)-2u_{1/2}(x+h)\big)-
\big(u_{1/2}(x+y)+u_{1/2}(x-y)-2u_{1/2}(x)\big)
}{h}$$
and we observe that, since~$u_{1/2}$ is globally
H\"older continuous with
exponent~$1/2$, it holds that
\begin{equation*}\begin{split}
|Q_h(x,y)|\;&\le
\frac{\big|
u_{1/2}(x+y+h)-u_{1/2}(x+y)\big|+\big|
u_{1/2}(x-y+h)-u_{1/2}(x-y)\big|+2\,\big|
u_{1/2}(x+h)-u_{1/2}(x)\big|
}{|h|}\\&\le\frac{\const |h|^{1/2}}{|h|}\\&=\frac{\const}{|h|^{1/2}},
\end{split}\end{equation*}
for any~$x$, $y\in\R$. Consequently, recalling~\eqref{0-er9858585859gigi}
and~\eqref{01019238397757548483993},
we conclude that
\begin{equation}\label{67034l1234567sevtyucoles02w945685r9fovk}
\lim_{h\to0} \left| \int_{I_x(h)} \frac{Q_h(x,y)}{|y|^{2}}\,dy\right|
\le \lim_{h\to0}
\int_{I_x(h)}
\frac{\const}{|h|^{1/2}\,\ell_x^{2}}\,dy\le
\lim_{h\to0}
\frac{\const |h|}{|h|^{1/2}\,\ell_x^{2}}=0.
\end{equation}
Now we take derivatives of~$u_{1/2}$. For this,
we observe that, for any~$\xi\in(-1,1)$,
$$ u_{1/2}'(\xi)=-\xi(1-\xi^2)^{-1/2}=-\xi u_{-1/2}(\xi).$$
Since the values outside~$(-1,1)$ are trivial, this implies that
\begin{equation}\label{92eir4950010020200220485}
u_{1/2}'(\xi)=-\xi u_{-1/2}(\xi)\quad{\mbox{ for any }}
\xi\in\R\setminus\{-1,1\}.
\end{equation}
Now, by~\eqref{9rfyw7rf203948ks}, we know that if~$y\in\R\setminus I_x(h)$
we have that~$x+y+t\in\R\setminus\{-1,1\}$ for all~$t\in\R$ with~$|t|<|h|$
and therefore we can exploit~\eqref{92eir4950010020200220485} and find that
\[ \lim_{h\to0} \frac{u_{1/2}(x+y+h)-u_{1/2}(x+y)}{h}=
-(x+y) u_{-1/2}(x+y).\]
Similar arguments show that, for any~$y\in\R\setminus I_x(h)$,
\begin{eqnarray*}
&& \lim_{h\to0} \frac{u_{1/2}(x-y+h)-u_{1/2}(x-y)}{h}=
-(x-y) u_{-1/2}(x-y)\\ {\mbox{and }}&&
\lim_{h\to0} \frac{u_{1/2}(x+h)-u_{1/2}(x)}{h}=
-x u_{-1/2}(x).
\end{eqnarray*}
Consequently, for any~$y\in\R\setminus I_x(h)$,
\begin{equation}\label{012394949494-1023949949-1023co}
\lim_{h\to0} \frac{Q_h(x,y)}{|y|^{2}}=
-\frac{(x+y) u_{-1/2}(x+y)+(x-y) u_{-1/2}(x-y)-2x u_{-1/2}(x)}{|y|^{2}}.
\end{equation}
Now we set
\begin{eqnarray*} \Xi_h(x,y)&:=&\frac{Q_h(x,y)\,\chi_{\R\setminus
I_x(h)}(y)}{|y|^{2}}\\&=&
\frac{1}{h\,|y|^{2}}\,\Big( \big(
u_{1/2}(x+y+h)+u_{1/2}(x-y+h)-2u_{1/2}(x+h)\big)\\&&\qquad-
\big(u_{1/2}(x+y)+u_{1/2}(x-y)-2u_{1/2}(x)\big)\Big)\,\,\chi_{\R\setminus
I_x(h)}(y)
\end{eqnarray*}
and we claim that
\begin{equation}\label{01-239402-39458723edf}
|\Xi_h(x,y)|\le C_x\,\left[\chi_{(-3,3)}(y)\,\left( 
\frac{1}{|1-(x+y)^2|^{1/2}}+\frac{1}{|1-(x-y)^2|^{1/2}}\right)+\frac{ \chi_{\R\setminus
(-3,3)}(y)}{|y|^2}
\right],
\end{equation}
for a suitable~$C_x>0$, possibly depending on~$x$.
For this, we first observe that if~$|y|\ge3$ then~$|x\pm y|\ge1$
and also~$|x\pm y+h|\ge1$. This implies that if~$|y|\ge3$, then~$u_{1/2}(x\pm y)=
u_{1/2}(x\pm y+h)=0$ and therefore
$$ \Xi_h(x,y)=
\frac{1}{h\,|y|^{2}}\,\big( 2u_{1/2}(x)-2u_{1/2}(x+h)\big).$$
This and the fact that~$u_{1/2}$ is smooth in the vicinity of the fixed~$x\in(-1,1)$
imply that~\eqref{01-239402-39458723edf} holds true when~$|y|\ge3$.
Therefore, from now on, to prove~\eqref{01-239402-39458723edf}
we can suppose that
\begin{equation}\label{938eufji394875489576584956ksdjhf-1}
|y|<3.
\end{equation}
We will also distinguish two regimes,
the one in which~$|y|\le\frac{\ell_x}4$
and the one in which~$|y|>\frac{\ell_x}4$.

If~$|y|\le\frac{\ell_x}4$ and~$|t|\le h$,
we have that
$$ |(x+y+t)+1|\ge |x+1|-|y|-|t|\ge\ell_x-|y|-|h|\ge \frac{\ell_x}2,$$
due to~\eqref{P0ow2992924rf},
and similarly $|(x-y+t)-1|\ge\frac{\ell_x}2$.
This implies that 
$$ |u_{1/2}(x+y+t)+u_{1/2}(x-y+t)-2u_{1/2}(x+t)|\le C_x\,|y|^2,$$
for some~$C_x>0$ that depends on~$\ell_x$. Consequently, we find that
if~$|y|\le\frac{\ell_x}4$ then
\begin{equation}\label{82376r293847001001}
|\Xi_h(x,y)|\le\frac{\const\,C_x\,|y|^2}{|y|^{2}}=\const C_x.\end{equation}
Conversely, if~$y\in\R\setminus I_x(h)$,
with
$|y|>\frac{\ell_x}4$, then we make use of~\eqref{92eir4950010020200220485}
and~\eqref{938eufji394875489576584956ksdjhf-1} to see that
\begin{equation}\label{0-1w033058-13933939wert}
\begin{split}&
|u_{1/2}(x+y+h)-u_{1/2}(x+y)|\le \int_0^{|h|} |u_{1/2}'(x+y+\tau)|\,d\tau\\
&\qquad= \int_0^{|h|} |x+y+\tau|\, |u_{-1/2}(x+y+\tau)|\,d\tau
\le 5\,\int_0^{|h|}  |u_{-1/2}(x+y+\tau)|\,d\tau\\
&\qquad\le5\,\int_0^{|h|} \frac{d\tau}{|1-(x+y+\tau)^2|^{1/2}}.
\end{split}\end{equation}
Also, if~$y\in\R\setminus I_x(h)$ we deduce from~\eqref{9rfyw7rf203948ks}
that~$|1\pm (x+y)|\ge 2 |h|$ and therefore, if~$|\tau|\le|h|$, then
$$ |1\pm(x+y+\tau)|\ge |1\pm (x+y)|-|\tau|\ge 
|1\pm (x+y)|-|h|\ge\frac{|1\pm (x+y)|}{2}.$$
Therefore
\begin{eqnarray*} &&|1-(x+y+\tau)^2|=
|1+(x+y+\tau)|\,|1-(x+y+\tau)|\\&&\qquad\ge \frac{1}{4}\,
|1+(x+y)|\,|1-(x+y)|=\frac{1}{4}\,|1-(x+y)^2|.\end{eqnarray*}
Hence, we insert this information into~\eqref{0-1w033058-13933939wert}
and we conclude that
\begin{equation}\label{01200-2-2020292929000:1}
|u_{1/2}(x+y+h)-u_{1/2}(x+y)|\le \const
\int_0^{|h|} \frac{d\tau}{|1-(x+y)^2|^{1/2}}=\frac{\const |h|}{|1-(x+y)^2|^{1/2}}
.
\end{equation}
Similarly, one sees that
\begin{equation}\label{01200-2-2020292929000:2}
|u_{1/2}(x-y+h)-u_{1/2}(x-y)|\le \frac{\const |h|}{|1-(x-y)^2|^{1/2}}
.
\end{equation}
In view of~\eqref{01200-2-2020292929000:1} and~\eqref{01200-2-2020292929000:2},
we get that, for any~$y\in\R\setminus I_x(h)$
with
$|y|>\frac{\ell_x}4$, 
\begin{eqnarray*} |\Xi_h(x,y)|&\le& \frac{1}{h\,|y|^2}\,
\left( \const|h|+\frac{\const |h|}{|1-(x+y)^2|^{1/2}}+
\frac{\const |h|}{|1-(x-y)^2|^{1/2}}\right)\\&\le&
\frac{\const}{\ell_x^2}\,
\left( 1+\frac{1}{|1-(x+y)^2|^{1/2}}+
\frac{1}{|1-(x-y)^2|^{1/2}}\right).\end{eqnarray*}
Combining this with~\eqref{82376r293847001001}, we obtain~\eqref{01-239402-39458723edf}, up to renaming constants.

Now, we point out that the right hand side of~\eqref{01-239402-39458723edf}
belongs to~$L^1(\R)$. Accordingly,
using~\eqref{01-239402-39458723edf} and 
the Dominated Convergence Theorem, and recalling also~\eqref{92eir4950010020200220485}, it follows that
\begin{eqnarray*}
&& \lim_{h\to0} \int_{\R\setminus I_x(h)} \frac{1}{h|y|^{2}}\,\Big(
\big( u_{1/2}(x+y+h)+u_{1/2}(x-y+h)-2u_{1/2}(x+h) \big)\\
&&\qquad-
\big( u_{1/2}(x+y)+u_{1/2}(x-y)-2u_{1/2}(x) \big)\Big)
\,dy \\ &&\qquad\qquad=
\lim_{h\to0} \int_{\R} \Xi_h(x,y)\,dy\\
&&\qquad\qquad=\int_\R \lim_{h\to0} \Xi_h(x,y)\,dy=\int_\R
\frac{ u_{1/2}'(x+y)+u_{1/2}'(x-y)-2u_{1/2}'(x) }{|y|^{2}}\,dy\\
&&\qquad\qquad=-\int_\R
\frac{(x+y) u_{-1/2}(x+y)+(x-y) u_{-1/2}(x-y)-2x u_{-1/2}(x)}{|y|^{2}}
\,dy,
\end{eqnarray*}
where the last identity is a consequence of~\eqref{012394949494-1023949949-1023co}.

{F}rom this and~\eqref{67034l1234567sevtyucoles02w945685r9fovk},
the claim in~\eqref{DER91-9dz} follows, as desired.

Now, we rewrite~\eqref{DER91-9dz} as
\begin{equation}\label{DER91-9dz-BIS}\begin{split}&
\frac{d}{dx} \int_\R \frac{u_{1/2}(x+y)+u_{1/2}(x-y)-2u_{1/2}(x)}{|y|^{2}}\,dy\\ =\;&
-{\mathcal{J}}(x)-x\int_\R
\frac{ u_{-1/2}(x+y)+ u_{-1/2}(x-y)-2 u_{-1/2}(x)}{|y|^{2}}
\,dy\\
&\qquad{\mbox{ where }}
{\mathcal{J}}(x):=
\int_\R
\frac{y \,\big(u_{-1/2}(x+y)- u_{-1/2}(x-y)\big)}{|y|^{2}}
\,dy=
\int_\R
\frac{u_{-1/2}(x+y)- u_{-1/2}(x-y)}{y}
\,dy.
\end{split}
\end{equation}
We claim that
\begin{equation}\label{DER91-9dz-TRIS}
{\mathcal{J}}(x)=0.
\end{equation}
This follows plainly for $x=0$, since~$u_{-1/2}$ is even. Hence, from here on, 
to prove~\eqref{DER91-9dz-TRIS}
we assume without loss of generality that $x\in(0,1)$.
Moreover, by changing variable~$y\mapsto-y$,
we see that
$$ -\PV\int_\R
\frac{ u_{-1/2}(x-y)}{y}
\,dy=\PV\int_\R
\frac{ u_{-1/2}(x+y)}{y}
\,dy $$
and therefore
\begin{equation}\label{0123983583095885rdf99494}\begin{split}&
{\mathcal{J}}(x)=
2\;\PV\int_\R\frac{u_{-1/2}(x+y)}{y}\;dy\ =\ 
2\;\PV\int_{-1-x}^{1-x}\frac{dy}{y\sqrt{1-(x+y)^2}}\\ &\qquad\qquad=\ 
2\;\PV\int_{-1}^1\frac{dz}{(z-x)\sqrt{1-z^2}}.
\end{split}\end{equation}
Now, we apply the change of variable 
\[
\xi:=\frac{1-\sqrt{1-z^2}}{z},\qquad{\mbox{hence }} z=\frac{2\xi}{1+\xi^2}.
\]
We observe that when $z$ ranges in $(-1,1)$, then~$\xi$ ranges therein as well. Moreover,
\[
\sqrt{1-z^2}=1-\xi z=\frac{1-\xi^2}{1+\xi^2}
,\]
thus, by~\eqref{0123983583095885rdf99494},
\begin{align*}
& {\mathcal{J}}(x)
=2\;\PV\int_{-1}^1\frac1{(\frac{2\xi}{1+\xi^2}-x)\frac{1-\xi^2}{1+\xi^2}}
\cdot\frac{2-2\xi^2}{(1+\xi^2)^2}\;d\xi \\
&\qquad=4\;\PV\int_{-1}^1\frac{d\xi}{2\xi-x(1+\xi^2)}
=4x\;\PV\int_{-1}^1\frac{d\xi}{1-x^2-(1-x\xi)^2}.
\end{align*}
We now apply another change of variable
\[
\eta:=\frac{1-x\xi}{\sqrt{1-x^2}}
\]
which gives
\begin{equation}
\label{48546374839904032539535939666} {\mathcal{J}}(x)=
\frac{4}{\sqrt{1-x^2}}\;\PV\int_{a_-}^{a^+}\frac{d\eta}{1-\eta^2},
\end{equation}
where
\[ a_+:=\sqrt{\frac{1+x}{1-x}}\; {\mbox{ and }} \; a_-:=\sqrt{\frac{1-x}{1+x}}=\frac1{a_+}.
\]
Now we notice that
\[
\PV\int_{a_-}^{a_+}\frac{d\eta}{1-\eta^2}=
\frac12\ln\left|\frac{(1+a_+)(1-a_-)}{(1-a_+)(1+a_-)}\right|=0.
\]
Inserting this identity into~\eqref{48546374839904032539535939666},
we obtain~\eqref{DER91-9dz-TRIS}, as desired.

Then, from~\eqref{DER91-9dz-BIS}
and~\eqref{DER91-9dz-TRIS} we get that
$$ \frac{d}{dx} \int_\R \frac{u_{1/2}(x+y)+u_{1/2}(x-y)-2u_{1/2}(x)}
{|y|^{2}}\,dy=
-x\int_\R
\frac{ u_{-1/2}(x+y)+ u_{-1/2}(x-y)-2 u_{-1/2}(x)}{|y|^{2}}
\,dy$$
that is
$$ \frac{d}{dx}(-\Delta)^{1/2} u_{1/2}=-x\,(-\Delta)^{1/2} u_{-1/2}\qquad{\mbox{in }}(-1,1).$$
{F}rom this and~\eqref{CERC} we infer that~$x\,(-\Delta)^{1/2} u_{-1/2}=0$
and so~$(-\Delta)^{1/2} u_{-1/2}=0$ in~$(-1,1)$.

These consideration establish~\eqref{CERC-HARM}, as desired. Now, we give
a brief probabilistic insight on it.
In probability - and in stochastic calculus - 
a measurable function $f:\R^n\to\R$ is said to be harmonic
in an open set $D\subset\R^n$ if, for any~$D_1\subset D$ and any $x\in D_1$,
\begin{equation}\label{87rte9wugidvqwifdweheri3i5575757}
\begin{split}
& f(x)\ =\ \mathbb E_x\left[f(W_{\tau_{D_1}})\right], \\
& \qquad\text{where $W_t$ is a Brownian motion and $\tau_{D_1}$
is the first exit time from $D_1$, namely} \\
& \qquad\tau=\inf\{t>0:W_t\not\in D_1\}.
\end{split}
\end{equation}
Notice that, since $W_t$ has (a.s.) continuous trajectories,
then (a.s.) $W_{\tau_{D_1}}\in\partial D_1$. This notion of harmonicity
coincides with the analytic one.

If one considers a L\'evy-type process $X_t$ in place of the 
Brownian motion, the definition of harmonicity (with respect
to this other process) can be given in the very same way.
When $X_t$ is an isotropic $(2s)$-stable process,
the definition amounts to having zero fractional Laplacian
${(-\Delta)}^s$ at every point of $D$
and replace~\eqref{87rte9wugidvqwifdweheri3i5575757} by
\[
f(x)\ =\ \mathbb E_x[f(X_{\tau_{D_1}})],
\qquad {\mbox{for any }}D_1\subseteq D.
\]
In this identity, we can consider a sequence
$\{D_j:D_j\subset D,j\in\N\}$, with~$D_j\nearrow D$, and equality
\begin{equation}\label{sing-sharm}
f(x)\ =\ \mathbb E_x[f(X_{\tau_{D_j}})],\qquad \text{for any }j\in\N.
\end{equation}
When $f=0$ in $\R^n\setminus D$, 
the right-hand side of \eqref{sing-sharm} can be not~$0$
(since~$X_{\tau_{D_j}}$ may also end up in~$D\setminus D_j$),
and this leaves the possibility of finding~$f$ which satisfies~\eqref{sing-sharm}
without vanish identically (an example of this phenomenon is
exactly given by the function~$u_{-1/2}$ in~\eqref{CERC-HARM}).

It is interesting to observe that if~$f$ vanishes outside~$D$
and does not vanish identically, then, the only possibility to
satisfy~\eqref{sing-sharm} is that~$f$ diverges along~$\partial D$.
Indeed, if~$|f|\le\kappa$,
since $f(X_{\tau_{D_j}})\neq 0$ only when $x\in D\setminus D_j$
and $|D\setminus D_j|\searrow 0$ as $j\to\infty$,
we would have that
$$ \lim_{j\to+\infty}\mathbb E_x[f(X_{\tau_{D_j}})]\le\lim_{j\to+\infty}\const
\kappa \,|D\setminus D_j|=0,$$
and \eqref{sing-sharm} would imply that~$f$ must
vanish identically.

Of course, the function~$u_{-1/2}$ in~\eqref{1092373957349tufh39597306434676949767666794tfoe}
embodies exactly this singular boundary behavior.

\section{Another proof of~\eqref{CERC-HARM}}\label{CERC-HARM:A:2}

Here we give a different proof of~\eqref{CERC-HARM}
by using complex analysis and extension methods.
We use the principal complex square root introduced in~\eqref{INC}
and, for any~$x\in\R$ and~$y>0$ we define
$$ U_{-1/2}(x,y):=\Re \left( \frac1{\SURD{1-z^2}}\right)
,$$ 
where~$z:=x+iy$.

\begin{figure}
    \centering
    \includegraphics[width=11cm]{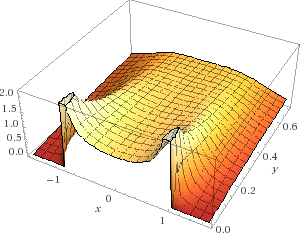}
    \caption{\it {{Harmonic extension in
the halfplane of the function~$\R\ni x\mapsto (1-x^2)^{-1/2}_+$.}}}
    \label{SPIKES}
\end{figure}

The function~$U_{-1/2}$ is plotted in Figure~\ref{SPIKES}.
We recall that the function~$U_{-1/2}$
is well-defined, thanks to~\eqref{inclus}. Also, the denominator
never vanishes when~$y>0$ and so~$U_{-1/2}$
is harmonic in the halfplane, being the real part of a holomorphic function
in such domain.

Furthermore, in light of~\eqref{J11KLA},
we have that
\[
\lim_{y\searrow0} \frac1{\SURD{1-z^2}} \,=\,
\left\{
\begin{matrix}
-\displaystyle\frac{i}{\sqrt{|1-x^2|} } & {\mbox{ if }}x^2>1,\\
\displaystyle\frac{1}{\sqrt{|1-x^2|}} & {\mbox{ if }}x^2<1,\\
+\infty & {\mbox{ if }}x^2=1,
\end{matrix}\right. \]
and therefore
\begin{eqnarray*}&& \lim_{y\searrow0} U_{-1/2}(x,y)
=\Re \left(\lim_{y\searrow0}
\frac1{\SURD{1-z^2}}\right)=
\left\{
\begin{matrix}
0 & {\mbox{ if }}x^2>1,\\
\displaystyle\frac{1}{\sqrt{|1-x^2|}} & {\mbox{ if }}x^2<1,\\
+\infty & {\mbox{ if }}x^2=1,
\end{matrix}\right.
\\ &&\qquad=(1-x^2)^{-1/2}_+=u_{-1/2}(x).\end{eqnarray*}
This gives that~$U_{-1/2}$ is the harmonic extension of~$u_{-1/2}$
to the halfplane. Therefore, by~\eqref{EX:APP:E}
\eqref{EX:APP:F} and~\eqref{NU99}, for any~$x\in(-1,1)$ we have that
\begin{eqnarray*}
-(-\Delta)^{-1/2} u_{1/2}(x)
&=&\lim_{y\searrow0} \partial_y U_{-1/2}(x,y)
\\&=&
\lim_{y\searrow0} \partial_y \left(
\Re \left( \frac1{\SURD{1-z^2}}\right)\right) \\
&=&
\lim_{y\searrow0} \Re\left(
\partial_y \left( \frac1{\SURD{1-z^2}}\right)\right) \\&=&
-\lim_{y\searrow0} \Re\left(
\left( \frac1{\SURD{1-z^2}}\right)^2
\partial_y \left( \SURD{1-z^2}\right)\right)\\&=&
-\frac{1}{1-x^2}\lim_{y\searrow0} \Re\left(
\partial_y \left( \SURD{1-z^2}\right)\right)\\&=&
0,
\end{eqnarray*}
that is~\eqref{CERC-HARM}.

\section{Proof of~\eqref{GAUSS:s} (based on Fourier methods)}\label{GAUSS:sA} 

When~$n=1$, we use~\eqref{SGA} to find that\footnote{As a historical remark,
we recall that~$e^{-|\xi|}$ is sometimes called the ``Abel Kernel''
and its Fourier Transform the ``Poisson Kernel'',
which in dimension~$1$ reduces to the ``Cauchy-Lorentz,
or Breit-Wigner,
Distribution'' (that has also classical geometric
interpretations as the ``Witch of Agnesi'', and so many names
attached to a single function clearly demonstrate its importance
in numerous applications).} 
\begin{equation}\label{209rewifhweot894tfhkfheo485749}\begin{split}
{\mathcal{G}}_{1/2}(x)\;&= \int_{\R} e^{-|\xi|} \,e^{ix\xi}\,d\xi
=\lim_{R\to+\infty}\int_0^{R} e^{-\xi} \,e^{ix \xi}\,d\xi
+\int_{-R}^0 e^{\xi} \,e^{ix \xi}\,d\xi\\
&= \lim_{R\to+\infty}\frac{e^{R(ix-1)} -1 }{ix-1}+\frac{1-e^{-R(ix+1)}}{ix+1}
= -\frac{1 }{ix-1}+\frac{1}{ix+1}=\frac{2}{x^2+1}.
\end{split}\end{equation}
This proves~\eqref{SGA} when~$n=1$.

Let us now deal with the case~$n\ge2$.
By changing variable~$Y:=1/y$, we see that
$$ \int_0^{+\infty} e^{-\frac{|\xi|\,\left(y-\frac1y\right)^2}{2}}\,dy
=\int_0^{+\infty} e^{-\frac{|\xi|\,\left(Y-\frac1Y\right)^2}{2}}\,\frac{dY}{Y^2}.$$
Therefore, summing up the left hand side to both sides of this identity
and using
the transformation~$\eta:=y-\frac1y$,
\begin{eqnarray*}2\int_0^{+\infty} e^{-\frac{|\xi|\,\left(y-\frac1y\right)^2}{2}}\,dy
&=&\int_0^{+\infty} \left(1+\frac1{y^2}\right)\,e^{-\frac{|\xi|\,\left(y-\frac1y\right)^2}{2}}\,dy\\
&=&\const\int_0^{+\infty} e^{-\frac{|\xi|\,\eta^2}{2}}\,d\eta\\&=&\frac{\const}{\sqrt{|\xi|}}.\end{eqnarray*}
As a result,
\begin{eqnarray*}e^{-|\xi|} = \frac{\const e^{-|\xi|}\,\sqrt{|\xi|}}{\sqrt{|\xi|}}
&=&\const e^{-|\xi|}\,\sqrt{|\xi|}\int_0^{+\infty} e^{-\frac{|\xi|\,\left(y-\frac1y\right)^2}{2}}\,dy
\\
&=& \const e^{-|\xi|}\,\sqrt{|\xi|}\int_0^{+\infty} e^{-\frac{|\xi|\,\left(y^2+\frac1{y^2}-2\right)}{2}}\,dy
\\
&=& \const \sqrt{|\xi|}\int_0^{+\infty} e^{-\frac{|\xi|\,\left(y^2+\frac1{y^2}\right)}{2}}\,dy
\\
&=& \const \int_0^{+\infty} \frac{1}{\sqrt{t}}\;e^{-\frac{t}{2}}\, e^{-\frac{|\xi|^2}{2t}}\,dt,
\end{eqnarray*}
where the substitution~$t:=|\xi|\,y^2$ has been used.

Accordingly, by~\eqref{SGA}, the Gaussian Fourier transform and the change of variable~$\tau:=
t(1+|x|^2)$,
\begin{eqnarray*}
{\mathcal{G}}_{1/2}(x)&=& \int_{\R^n} e^{-|\xi|} \,e^{ix\cdot\xi}\,d\xi\\
&=& \const \iint_{\R^n\times(0,+\infty)}
\frac{1}{\sqrt{t}}\;
e^{-\frac{t}{2}}\, e^{-\frac{|\xi|^2}{2t}}\,e^{ix\cdot\xi}
\,d\xi\,dt\\
&=&
\const \int_{(0,+\infty)}
t^{\frac{n-1}{2}} e^{-\frac{t}{2}}\, e^{-\frac{ t|x|^2}2} \,dt
\\&=&
\const \int_{(0,+\infty)}
\left(\frac\tau{1+|x|^2}\right)^{\frac{n-1}{2}} e^{-\frac\tau2} \,\frac{d\tau}{1+|x|^2}\\&=&
\frac{\const}{ \big(1+|x|^2\big)^{\frac{n+1}{2}} }.
\end{eqnarray*}
This establishes~\eqref{GAUSS:s}.

\section{Another proof of~\eqref{GAUSS:s} (based on extension methods)}\label{GAUSS:FAC}

\begin{figure}
    \centering
    \includegraphics[width=11cm]{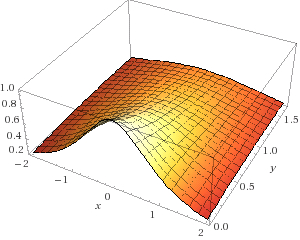}
    \caption{\it {{Harmonic extension in
the halfplane of the function~$\R\ni x\mapsto
\frac{1}{1+|x|^2}$.}}}
    \label{GAFACI}
\end{figure}

The idea is to consider the fundamental solution in the extended space and take a derivative
(the time variable acting as a translation
and, to favor the intuition, one may keep in mind
that the Poisson kernel is the normal derivative of the Green function). Interestingly, this proof
is, in a sense,  ``conceptually simpler'', and ``less technical'' than that
in Appendix~\ref{GAUSS:sA}, thus demonstrating that, at least in some cases,
when appropriately used, fractional methods may lead to cultural advantages\footnote{Let us mention another conceptual simplification of nonlocal problems:
in this setting, the integral representation often allows the formulation
of problems with minimal requirements on the functions involved
(such as measurability
and possibly minor pointwise or integral bounds). Conversely,
in the classical setting, even to just formulate a problem, one
often needs assumptions
and tools from functional analysis, comprising e.g.
Sobolev differentiability, distributions or functions of bounded
variations.}
with respect
to more classical approaches.

For this proof, we consider variables~$X:=(x,y)\in\R^n\times(0,+\infty)\subset\R^{n+1}$
and fix~$t>0$.
We let~$\Gamma$ be the fundamental solution in~$\R^{n+1}$, namely
$$ \Gamma(X):=\left\{
\begin{matrix}
-\const\log|X| &{\mbox{ if }}n=1,\cr
\displaystyle\frac{\const}{|X|^{n-1}}&{\mbox{ if }}n\ge2.
\end{matrix}
\right. $$
By construction~$\Delta\Gamma$ is the Delta Function at the origin
and so, for any~$t>0$, we have that~$\tilde\Gamma(X;t)=\tilde\Gamma(x,y;t):=\Gamma(x,y+t)$
is harmonic for~$(x,y)\in\R^n\times(0,+\infty)$.
Accordingly, the function~$U(x,y;t):=\partial_y\tilde\Gamma(x,y;t)$ is 
also harmonic for~$(x,y)\in\R^n\times(0,+\infty)$.
We remark that
$$ U(x,y;t)=\partial_y \Gamma(x,y+t)=
\displaystyle\frac{\const}{|(x,y+t)|^{n}} \partial_y \sqrt{|x|^2+(y+t)^2}
=
\displaystyle\frac{\const(y+t)}{|(x,y+t)|^{n+1}} =
\displaystyle\frac{\const(y+t)}{\big(|x|^2+(y+t)^2\big)^{\frac{n+1}{2}}}.$$
This function is plotted in Figure~\ref{GAFACI} (for the model case in the plane).
We observe that
$$ \lim_{y\searrow0}U(x,0;t)=
\displaystyle\frac{\const \;t}{\big(|x|^2+t^2\big)^{\frac{n+1}{2}}}
=
\displaystyle\frac{\const}{t^n\,\big(1+(|x|/t)^2\big)^{\frac{n+1}{2}}} =: u(x,t).$$
As a consequence, by~\eqref{EX:APP:E}
and~\eqref{EX:APP:F} (and noticing that the role played by
the variables~$y$ and~$t$ in the function~$U$ is the same),
\begin{eqnarray*}&& -(-\Delta)^{1/2} u(x,t)=
\lim_{y\searrow0} \partial_y U(x,y;t)=
\lim_{y\searrow0} \partial_y \displaystyle\frac{\const(y+t)}{\big(|x|^2+(y+t)^2\big)^{\frac{n+1}{2}}}
=
\lim_{y\searrow0} 
\partial_t \displaystyle\frac{\const(y+t)}{\big(|x|^2+(y+t)^2\big)^{\frac{n+1}{2}}}
\\&&\qquad\qquad=
\partial_t \displaystyle\frac{\const\; t}{\big(|x|^2+t^2\big)^{\frac{n+1}{2}}}
=\partial_t u(x,t).\end{eqnarray*}
This shows that~$u$ solves the fractional heat equation,
with~$u$ approaching a Delta function when~$t\searrow0$. Hence
$${\mathcal{G}}_{1/2}(x)=u(x,1)=\displaystyle\frac{\const}{\big(1+|x|^2\big)^{\frac{n+1}{2}}}
,$$
that is~\eqref{GAUSS:s}.

\section{Proof of~\eqref{GIA}}\label{GIAGIA}

First, we construct a useful barrier. Given~$A>1$, we define
$$ w(t):=\left\{
\begin{matrix}
A & {\mbox{ if }} |t|\le 1,\\
t^{-1-2s} & {\mbox{ if }} |t|>1.
\end{matrix}
\right.$$
We claim that if~$A$ is sufficiently large, then
\begin{equation}\label{BUONA}
(-\Delta)^s w(t)<-3w(t)
\quad{\mbox{ for all }}t\in\R\setminus (-3,3).
\end{equation}
To prove this, fix~$t\ge3$ (the case~$t\le-3$ being similar). Then, if~$|\xi-t|<1$, we have that
$$ \xi\ge t-1=\frac{2t}{3}+\frac{t}3-1\ge\frac{2t}{3}.$$
As a consequence, if~$|\tau-t|<1$, 
\begin{eqnarray*}
&& \big| w(t)-w(\tau)+\dot w(t)(\tau-t)\big|\le
\sup_{|\xi-t|<1} |\ddot w(\xi)|\,|t-\tau|^2\\
&&\qquad
\le\const \sup_{\xi\ge 2t/3} \xi^{-3-2s}\,|t-\tau|^2
\le\const t^{-3-2s}\,|t-\tau|^2.\end{eqnarray*}
This implies that
\begin{equation}\label{337.1}
\begin{split} &
\int_{\{|\tau-t|<1\}} \frac{w(t)-w(\tau)}{|t-\tau|^{1+2s}}\,d\tau=
\int_{\{|\tau-t|<1\}} \frac{w(t)-w(\tau)+\dot w(t)(\tau-t)}{|t-\tau|^{1+2s}}\,d\tau\\
&\qquad\le
\const t^{-3-2s}\int_{\{|\tau-t|<1\}} \frac{|t-\tau
|^2}{|t-\tau|^{1+2s}}\,d\tau=
\const t^{-3-2s}\\&\qquad
\le \const t^{-1-2s}=\const w(t).
\end{split}
\end{equation}
On the other hand,
\begin{equation}\label{337.2}
\int_{\{|\tau-t|\ge1\}\cap\{|\tau|>1\}} \frac{w(t)-w(\tau)}{|t-\tau|^{1+2s}}\,d\tau\le
\int_{\{|\tau-t|\ge1\}} \frac{w(t)}{|t-\tau|^{1+2s}}\,d\tau\le\const w(t)
.\end{equation}
In addition, if~$|\tau|\le1$ then~$|\tau-t|\ge t-\tau\ge 3-1>1$, hence
$$ \{|\tau-t|\ge1\}\cap\{|\tau|\le1\}=\{|\tau|\le1\}.$$
Accordingly,
\begin{equation}\label{337}
\int_{\{|\tau-t|\ge1\}\cap\{|\tau|\le1\}}
\frac{w(t)-w(\tau)}{|t-\tau|^{1+2s}}\,d\tau
=
\int_{\{|\tau|\le1\}} \frac{t^{-1-2s}-A}{|t-\tau|^{1+2s}}\,d\tau\le
\int_{\{|\tau|\le1\}} \frac{1-A}{|t-\tau|^{1+2s}}\,d\tau.
\end{equation}
We also observe that if~$|\tau|\le1$ then~$|t-\tau|\le t+1\le 2t$
and therefore
$$ \int_{\{|\tau|\le1\}} \frac{d\tau}{|t-\tau|^{1+2s}}\ge\frac{\const}{t^{1+2s}}=
{\const}w(t).$$
So, we plug this information into~\eqref{337},
assuming~$A>1$ and we obtain that
\begin{equation}\label{337.3}
\int_{\{|\tau-t|\ge1\}\cap\{|\tau|\le1\}} \frac{w(t)-w(\tau)}{|t-\tau|^{1+2s}}\,d\tau
\le -{(A-1)\,\const}w(t).
\end{equation}
Thus, gathering together the estimates in~\eqref{337.1}, \eqref{337.2} and~\eqref{337.3},
we conclude that
$$ \int_\R\frac{w(t)-w(\tau)}{|t-\tau|^{1+2s}}\,d\tau\le\const w(t)
-{(A-1)\,\const}w(t)\le -4 w(t)<-3w(t),$$
as long as~$A$ is sufficiently large.
This proves~\eqref{BUONA}.

Now, to prove~\eqref{GIA},
we define~$v:=\dot u>0$. {F}rom~\eqref{PUB}, we know that
\begin{equation}\label{09wfyheiogfvudigfcwerfygrgfvior}
(-\Delta)^s v=(1-3u^2) v\ge -3u^2v\ge -3v. 
\end{equation}
Given~$\e>0$, we define
$$ w_\e(t):= \frac{\iota}{A}\,w(t)-\e,\qquad{\mbox{ where }}
\iota:=\min_{t\in[-3,3]} v(t).$$
We claim that
\begin{equation}\label{SOt}
w_\e\le v. 
\end{equation}
Indeed, for large~$\e$, it holds that~$w_\e<0<v$ and so~\eqref{SOt}
is satisfied. 
In addition, for any~$\e>0$,
\begin{equation}
\label{90uwefguivrifhgvweioghreherog}
\lim_{t\to+\infty} w_\e(t)=-\e<0\le\inf_{t\in\R} v(t).\end{equation}
Suppose now that~$\e_\star>0$ produces a touching point
between~$w_{\e_\star}$ and~$v$, namely~$w_{\e_\star}\le v$
and~$w_{\e_\star}(t_\star)=v(t_\star)$ for some~$t_\star\in\R$.
Notice that, if~$|\tau|\le3$,
$$ w_{\e_\star}(\tau)\le
\frac{\iota}{A}\,\sup_{t\in\R}w(t)-\e\le
\iota-\e =\min_{t\in[-3,3]} v(t)-\e\le v(\tau)-\e<v(\tau),$$
and therefore~$|t_\star|>3$.
Accordingly, if we set~$v_\star:=v-w_{\e_\star}$, using~\eqref{BUONA}
and~\eqref{09wfyheiogfvudigfcwerfygrgfvior}, we see that
\begin{eqnarray*}
&& 0=-3 v_\star(t_\star)= -3v(t_\star)+3w_{\e_\star}(t_\star)
\le (-\Delta)^s v(t_\star)-(-\Delta)^s w_{\e_\star}(t_\star)
\\ &&\qquad\quad\qquad
=(-\Delta)^s v_\star(t_\star) =\int_\R \frac{
v_\star(t_\star)-v_\star(\tau)
}{|t_\star-\tau|^{1+2s}}\,d\tau
=-\int_\R \frac{v_\star(\tau)
}{|t_\star-\tau|^{1+2s}}\,d\tau.
\end{eqnarray*}
Since the latter integrand is nonnegative, we conclude that~$v_\star$
must vanish identically, and thus~$w_{\e_\star}$ must coincide with~$v$.
But this is in contradiction with~\eqref{90uwefguivrifhgvweioghreherog}
and so the proof of~\eqref{SOt} is complete.

Then, by sending~$\e\searrow0$ in~\eqref{SOt} we find that~$v\ge \frac\iota{A}\,w$,
and therefore, for~$t\ge1$ we have that~$\dot u(t)=
v(t)\ge \kappa t^{-1-2s}$, for all~$t>1$,
for some~$\kappa>0$.

Consequently, for any~$t>1$,
$$ 1-u(t) = \lim_{T\to+\infty} u(T)-u(t)
=\lim_{T\to+\infty} \int^T_t \dot u(\tau)\,d\tau
=\int^{+\infty}_t \dot u(\tau)\,d\tau\ge
\kappa \int^{+\infty}_t \tau^{-1-2s}\,d\tau =\frac{\kappa}{2s}\,t^{-2s},$$
and a similar estimates holds for~$1+u(t)$ when~$t<-1$.

These considerations establish~\eqref{GIA}, as desired.

\section{Proof of~\eqref{78:87}}\label{0987poiutyui19}

\begin{figure}
    \centering
    \includegraphics[width=11cm]{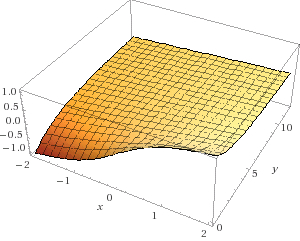}
    \caption{\it {{Harmonic extension in
the halfplane of the function~$\R\ni x\mapsto\frac2\pi\,\arctan x$.}}}
    \label{MKXX}
\end{figure}

Here we prove that~\eqref{78:87} is a solution of~\eqref{OEBN:11}.
The idea of the proof, as showed in Figure~\ref{MKXX},
is to consider the harmonic extension 
of the function~$\R\ni x\mapsto\frac2\pi\,\arctan x$ in the halfplane~$\R\times(0,+\infty)$
and use the method described in~\eqref{EX:APP:E}
and~\eqref{EX:APP:F}.

We let
\[ U(x,y):=\frac{2}{\pi}\arctan\frac{x}{y+1}.\]
The function~$U$ is depicted\footnote{In complex variables, one can
also interpret the function~$U$ in terms of the principal argument function
$$ {\rm Arg}(r e^{i\varphi})=\varphi\in(-\pi,\pi],$$
with branch cut along the nonpositive
real axis.
Notice indeed that, if~$z=x+iy$ and~$y>0$,
$$ {\rm Arg}(z+i)=\frac\pi2-\arctan\frac{x}{y+1}=\frac\pi2\left(1-U(x,y) \right).$$
This observation would also lead to~\eqref{0229}.}
in Figure~\ref{MKXX}.
Of course, it coincides with~$u$ when~$y=0$
and, for any~$x\in\R$ and~$y>0$,
\begin{equation}\label{0229}
\frac\pi2\,\Delta U(x,y)=
-\frac{2 x (1 + y)}{(x^2 + (1 + y)^2)^2}+
\frac{2 x (1 + y)}{(x^2 + (1 + y)^2)^2}=0.
\end{equation}
Hence, the setting in~\eqref{EX:APP:E} is satisfied and so, in light of~\eqref{EX:APP:F}.
we have
\begin{equation}\label{89:1Q}
(-\Delta)^{1/2} u(x)=-\lim_{y\searrow0} \partial_y U(x,y)=\frac2\pi \lim_{y\searrow0}
\frac{x}{x^2 + (1 + y)^2}
=\frac{2x}{\pi\,(x^2 + 1)}
\end{equation}
Also, by
the trigonometric Double-angle Formula, for any~$\theta\in\left(-\frac\pi2,\frac\pi2\right)$,
$$ \sin(2\theta )=2\sin \theta \cos \theta ={\frac {2\tan \theta }{\tan ^{2}\theta +1}}.$$
Hence, taking~$\theta:=\arctan x$,
$$ \sin(\pi u(x))=\sin(2\arctan x )={\frac {2x }{x^2+1}}.$$
This and~\eqref{89:1Q} show
that~\eqref{78:87} is a solution of~\eqref{OEBN:11}.

\section{Another proof of~\eqref{78:87} (based on~\eqref{GAUSS:s})}\label{APP:CC}

This proof of~\eqref{78:87} is based on the fractional heat kernel in~\eqref{GAUSS:s}.
This approach has the advantage of being quite general (see e.g. Theorem~3.1 in~\cite{MR3280032})
and also to relate the two ``miraculous'' explicit formulas~\eqref{GAUSS:s}
and~\eqref{78:87}, which are available only in the special case~$s=1/2$.

For this, we let~$P=P(x,t)$ the fundamental solution of the heat flow in~\eqref{DHE}
with~$n=1$ and~$s=1/2$. Notice that, by~\eqref{GAUSS:s}, we know that
\begin{equation}\label{GAUSS:s:INAP} 
P(x,1)={\mathcal{G}}_{1/2}(x)=\frac{c}{1+x^2},\end{equation}
with
$$ c:=\left( \int_\R \frac{dx}{1+x^2}\right)^{-1}=\frac1\pi.$$
Also, by scaling,
\begin{equation}\label{GAUSS:s:INAP:2} 
P(x,t)=t^{-1} P(t^{-1}x,1)=t^{-1}{\mathcal{G}}_{1/2}(t^{-1}x)
.\end{equation}
For any~$x\in\R$ and any~$t>0$, we define
\begin{equation} \label{GAUSS:s:INAP:3}
U(x,t):= 2\int_0^x P(\eta, t+1)\,d\eta.\end{equation}
In light of~\eqref{GAUSS:s:INAP:2}, we see that
$$ |U(x,t)|\le 2(t+1)^{-1} \int_0^x {\mathcal{G}}_{1/2}\big((t+1)^{-1}\eta\big)\,d\eta
=2 \int_0^{(t+1)^{-1}x} {\mathcal{G}}_{1/2}(\zeta)\,d\zeta,$$
which is bounded in~$\R\times[0,+\infty)$, and infinitesimal as~$t\to+\infty$ for any fixed~$x\in\R$.

Notice also that
$$ \partial^2_t P = \partial_t(\partial_t P)=\partial_t (-\Delta)^{1/2}P=
(-\Delta)^{1/2}\partial_t P=(-\Delta)^{1/2}(-\Delta)^{1/2} P=-\partial^2_x P,$$
by~\eqref{QUASDA}.
As a consequence,
\begin{equation}\label{HA:AP:013}
\begin{split}
\frac12(\partial^2_x+\partial^2_t)U(x,t)\,&=
\partial_x P(x,t+1)+
\int_0^x \partial^2_t P(\eta, t+1)\,d\eta
\\
&=
\partial_x P(x,t+1)-
\int_0^x \partial^2_x P(\eta, t+1)\,d\eta\\
&=
\partial_x P(0,t+1)\\&=0,
\end{split}\end{equation}
where the last identity follows from~\eqref{GAUSS:s:INAP:2}.

Besides, from~\eqref{GAUSS:s:INAP:2} we have that
$$ \partial_t P(x,t)= \partial_t\Big(
t^{-1}{\mathcal{G}}_{1/2}(t^{-1}x) \Big)
= -t^{-2}{\mathcal{G}}_{1/2}(t^{-1}x)
-t^{-3} x{\mathcal{G}}_{1/2}'(t^{-1}x)$$
and so
$$ -\partial_t P(x,1)={\mathcal{G}}_{1/2}(x)
+x{\mathcal{G}}_{1/2}'(x)
= \partial_x\big( x\,{\mathcal{G}}_{1/2}(x)\big).
$$
In view of this, we have that
\begin{equation}\label{HA:AP:014}
\partial_t U(x,0)= 2\int_0^x \partial_t P(\eta, 1)\,d\eta=
2\int_0^x \partial_\eta\big( \eta\,{\mathcal{G}}_{1/2}(\eta)\big)\,d\eta=
2x\,{\mathcal{G}}_{1/2}(x).
\end{equation}
Accordingly, from~\eqref{HA:AP:013}
and~\eqref{HA:AP:014}, using the extension method in~\eqref{EX:APP:E}
and~\eqref{EX:APP:F} (with the variable~$y$ called~$t$ here), we conclude that, if
\begin{equation*} u(x):=U(x,0),\end{equation*}
then
\begin{equation} \label{0923883456LA}
(-\Delta)^{1/2}u(x)=2x\,{\mathcal{G}}_{1/2}(x).\end{equation}
We remark that, by~\eqref{GAUSS:s:INAP} and~\eqref{GAUSS:s:INAP:3},
\begin{equation}\label{0rfyuiefioqwgrweugffgfgfg} u(x)=
2c \int_0^x \frac{d\eta}{1+ x^2}
= \frac{2}\pi\,\arctan x.
\end{equation}
This, \eqref{GAUSS:s:INAP} and~\eqref{0923883456LA} give that
$$ (-\Delta)^{1/2}u(x)=\frac{1}{\pi}\,\frac{2x}{1+x^2}=\frac{1}{\pi}\,\sin(2\arctan x)=
\frac{1}{\pi}\,\sin\big(\pi u( x)\big),
$$
that is~\eqref{78:87}, as desired.

\section{Proof of \eqref{LAP NONLOC}}\label{APP LAP NONLOC}
Due to translation invariance, we can reduce ourselves to proving~\eqref{LAP NONLOC}
at the origin.
We consider a measurable $u:\R^n\to\R$ such that
\begin{equation*}
\int_{\R^n}\frac{|u(y)|}{1+|y|^{n+2}}\ <\ +\infty.\end{equation*}
Assume first that
\begin{equation}\label{Lape}
{\mbox{$u(x)=0$ for any $x\in B_r$,}} 
\end{equation}
for some~$r>0$.
As a matter of fact, under these assumptions on $u$, the right-hand side of \eqref{LAP NONLOC}
vanishes at $0$ regardless the size of $r$.
Indeed,
\begin{align*}
&{\phantom{=}} \ \left.\int_{\R^n}\frac{u(x+2y)+u(x-2y)-4u(x+y)-4u(x-y)+6u(x)}{{|y|}^{n+2}}
\;dy\right|_{x=0} \\
& =\ \int_{\R^n}\frac{u(2y)+u(-2y)-4u(y)-4u(-y)}{{|y|}^{n+2}}\;dy\  
 =\ 2\int_{\R^n\setminus B_{r/2}}\frac{u(2y)}{{|y|}^{n+2}}\;dy
 -8\int_{\R^n\setminus B_r}\frac{u(y)}{{|y|}^{n+2}}\;dy\\ &
 =\ 2\int_{\R^n\setminus B_{r}}\frac{2^{n+2}\, u(Y)}{2^n\,{|Y|}^{n+2}}\;dY
 -8\int_{\R^n\setminus B_r}\frac{u(y)}{{|y|}^{n+2}}\;dy
 \ =\ 0.
\end{align*}
This proves~\eqref{LAP NONLOC}
under the additional assumption in~\eqref{Lape}, that we are now going to remove.
To this end, for $r\in(0,1)$, denote by $\chi_r$ the characteristic function of $B_r$, i.e.
$\chi_r(x)=1$ if $x\in B_r$ and $\chi_r(x)=0$ otherwise.
Consider now $u\in C^{2,\alpha}(B_r)$, for some $\alpha\in(0,1)$, with
\begin{equation}\label{9woeddifi123456fi}u(0)=|\nabla u(0)|=0\end{equation}
for simplicity (note that one can always modify $u$ by considering
$\tilde u(x)=u(x)-u(0)-\nabla u(0)\cdot x$ and without affecting 
the operators in \eqref{LAP NONLOC}). 
Then, the right hand side of~\eqref{LAP NONLOC} becomes in this case
\begin{align*}
& \phantom{=}\ \int_{\R^n}\frac{u(2y)+u(-2y)-4u(y)-4u(-y)}{{|y|}^{n+2}}\;dy\ =\
2\int_{\R^n}\frac{u(2y)-4u(y)}{{|y|}^{n+2}}\;dy\ =\\ 
& =\ 2\int_{\R^n}\frac{u(2y)\chi_r(2y)-4u(y)\chi_r(y)}{{|y|}^{n+2}}\;dy\ +\ 
2\int_{\R^n}\frac{u(2y)(1-\chi_r(2y))-4u(y)(1-\chi_r(y))}{{|y|}^{n+2}}\;dy.
\end{align*}
The second addend is trivial for any $r\in(0,1)$, in view of the above remark,
since $u(1-\chi_r)$ is constantly equal to $0$ in $B_r$.
For the first one, we have
\begin{align}\label{2930ofoetu304868P0}
& \int_{\R^n}\frac{u(2y)\chi_r(2y)-4u(y)\chi_r(y)}{{|y|}^{n+2}}\;dy\ =\
\int_{B_{r/2}}\frac{u(2y)-4u(y)}{{|y|}^{n+2}}\;dy
-4\int_{B_r\setminus B_{r/2}}\frac{u(y)}{{|y|}^{n+2}}\;dy.
\end{align}
Now, we recall~\eqref{9woeddifi123456fi} and we notice that
\begin{align*}
|u(2y)-4u(y)|\leq\|u\|_{C^{2,\alpha}(B)}|y|^{2+\alpha},
\end{align*}
which in turn
implies that
\begin{align}\label{2930ofoetu304868P}
\left|\int_{B_{r/2}}\frac{u(2y)-4u(y)}{{|y|}^{n+2}}\;dy\right|\ \leq\ \const\|u\|_{C^{2,\alpha}(B)}r^\alpha.
\end{align}
On the other hand, a Taylor expansion and~\eqref{9woeddifi123456fi} yield
\begin{equation}\label{0234oottoo394588686}
\begin{split}
& {\phantom{=}} \ \int_{B_r\setminus B_{r/2}}\frac{u(y)}{{|y|}^{n+2}}\;dy\ =\ 
\int_{r/2}^r\frac{1}{\rho^{n+2}}\int_{\partial B_\rho}u(y)\;dy\;d\rho\  \\
& =\ \int_{r/2}^r\frac{1}{\rho^{3}}\int_{\partial B_1}u(\rho\theta)\;d\theta\;d\rho\ =\ 
\int_{r/2}^r\frac{1}{2\rho}\int_{\partial B_1}\left( D^2u(0)\theta\cdot\theta+\eta(\rho\theta)\right)\;d\theta\;d\rho \\
& =\ \const\Delta u(0)+\int_{r/2}^r\frac{1}{2\rho}\int_{\partial B_1}\eta(\rho\theta)\;d\theta\;d\rho
\end{split}
\end{equation}
in view of \eqref{DEL:CLA}, for some $\eta:B_r\to\R$ such that $|\eta(x)|\leq c|x|^\alpha$.
From this, \eqref{2930ofoetu304868P0} and~\eqref{2930ofoetu304868P} we deduce that
\begin{align*}
& \int_{\R^n}\frac{u(2y)+u(-2y)-4u(y)-4u(-y)}{{|y|}^{n+2}}\;dy\ =\
-\const\lim_{r\searrow 0}\int_{B_r\setminus B_{r/2}}\frac{u(y)}{{|y|}^{n+2}}\;dy 
\ =\ -\const\Delta u(0)
\end{align*}
which finally justifies \eqref{LAP NONLOC}.\medskip

It is interesting to remark that the main contribution to prove~\eqref{LAP NONLOC}
comes in this case
from the ``intermediate ring'' in~\eqref{0234oottoo394588686}.

\section{Proof of~\eqref{NO:REGO}}\label{NO:REGOA}

Take for instance $\Omega$ to be the unit ball and~$\bar u=1-|x|^2$.
Suppose that~$\|\bar u-v_\e\|_{C^2(\Omega)}\leq\e$.
Then, for small~$\e$, if~$x\in\R^n\setminus B_{1/2}$ it holds that
$$ v_\e(x)\le \bar u(x)+\e= 1-|x|^2+\e\le \frac34+\e\le\frac45,$$
while
$$ v_\e(0)\ge\bar u(0)-\e=1-\e\ge\frac56.$$
This implies that there exists~$x_\e\in \overline{B_{1/2}}$ such that
$$ v_\e(x_\e)=\sup_{B_1} v_\e \ge\frac56>\frac45\ge\sup_{B_1\setminus B_{1/2}}v_\e.$$
As a result,
\begin{eqnarray*}
\PV \int_{\Omega} \frac{v_\e(x_\e)-v_\e(y)}{|x_\e-y|^{n+2s}}\,dy\ge
\int_{B_1\setminus B_{3/4}} \frac{v_\e(x_\e)-v_\e(y)}{|x_\e-y|^{n+2s}}\,dy
\ge
\int_{B_1\setminus B_{3/4}} 
\left(\frac56-\frac45\right)\,dy
\ge\const.
\end{eqnarray*}
This says that~$(-\Delta)^s v_\e$ cannot vanish at~$x_\e$ and so~\eqref{NO:REGO} is proved.

\section{Proof of~\eqref{XAB}}\label{09wdy8gfuiigsr3489trfhhdghhghghg}

Let us first notice that the identity 
\begin{equation}\label{1}
\lambda^s\ =\ \frac{s}{\Gamma(1-s)}
\int_0^\infty\frac{1-e^{-t\lambda}}{t^{1+s}}\;dt
\end{equation}
holds for any $\lambda>0$ and~$s\in(0,1)$, because
\begin{align*}
\int_0^\infty\frac{1-e^{-t}}{t^{1+s}}\;dt
=\left.\frac{1-e^{-t}}{-s\,t^s}\right|_0^\infty+\frac1s\int_0^\infty\frac{e^{-t}}{t^s}\;dt=\frac{\Gamma(1-s)}s.
\end{align*}
We also observe that
when $u\in C^\infty_0(\Omega)$, the coefficients $\hat u_j$ decay fast as $j\to\infty$: indeed
\[
\hat u_j=-\frac1{\mu_j}\int_\Omega u\,\Delta\psi_j=-\frac1{\mu_j}\int_\Omega \psi_j\,\Delta u
=\ldots=(-1)^k\frac1{\mu_j^k}\int_\Omega \psi_j\,\Delta^k u.
\]
Therefore, applying equality \eqref{1} to the $\mu_j$'s in \eqref{spec-neu}
we obtain\footnote{The representation in~\eqref{boch-neu}
makes sense for a larger class of functions with respect 
to \eqref{spec-neu}, so in a sense \eqref{boch-neu} can 
be interpreted as an extension of definition \eqref{spec-neu}.}
\begin{equation}\label{boch-neu}
{(-\Delta)^s_{N,\Omega}} u\ =\ 
\frac{s}{\Gamma(1-s)}\sum_{j=0}^{+\infty}\int_0^\infty\frac{\hat u_j\psi_j-e^{-t\mu_j}\hat u_j\psi_j}{t^{1+s}}\;dt
=\frac{s}{\Gamma(1-s)}\int_0^\infty
\frac{u-e^{t{\Delta_{N,\Omega}}u}}{t^{1+s}}\;dt,\quad u\in C^\infty_0(\Omega)
\end{equation}
where $\{ e^{t \Delta_{N,\Omega} }\}_{t>0}$ stands for the 
heat semigroup associated to ${\Delta_{N,\Omega}}$. i.e. $e^{t{\Delta_{N,\Omega}}}u$ solves
$$ \left\{
\begin{matrix}
\partial_t v(x,t)=\Delta v(x,t) {\mbox{ in }}\Omega\times(0,\infty) \\
\partial_\nu v(x,t)=0 \text{ on }\partial\Omega\times(0,\infty) \\
v(x,0)=u(x)\text{ on }\Omega\times\{0\}.
\end{matrix}
\right.$$
To check \eqref{boch-neu}, it is sufficient to notice that
\[
\partial_t\left(\sum_{j=0}^{+\infty}e^{-t\mu_j}\hat u_j\psi_j\right)=
-\sum_{j=0}^{+\infty}\mu_je^{-t\mu_j}\hat u_j\psi_j=
\sum_{j=0}^{+\infty}e^{-t\mu_j}\hat u_j\Delta\psi_j=
\Delta\left(\sum_{j=0}^{+\infty}e^{-t\mu_j}\hat u_j\psi_j\right)
\]
and that
\[
\left.\left(\sum_{j=0}^{+\infty}e^{-t\mu_j}\hat u_j\psi_j\right)\right|_{t=0}
=\sum_{j=0}^{+\infty}\hat u_j\psi_j=u.
\]
Under suitable regularity assumptions on $\Omega$,
write now the heat semigroup in terms of the heat kernel~$p_{N}^\Omega$ as
\begin{equation}
e^{t{\Delta_{N,\Omega}}}u(x)\ =\ \int_\Omega p_{{N}}^\Omega(t,x,y)\,u(y)\;dy,
\qquad x\in\Omega,\ t>0
\end{equation}
where the following two-sided estimate on $p_{{N}}^\Omega$ holds
(see, for example, \cite[Theorem 3.1]{MR2648271})
\begin{equation}\label{heatk-est}
\frac{c_1\:e^{-c_2|x-y|^2/t}}{t^{n/2}}
\ \leq\ p_{{N}}^\Omega(t,x,y)\ \leq\ 
\frac{c_3\:e^{-c_4|x-y|^2/t}}{t^{n/2}},
\qquad x,y\in\Omega,\ t,c_1,c_2,c_3,c_4>0.
\end{equation}
Recall also that $p_{N}^\Omega(t,x,y)=p_{N}^\Omega(t,y,x)$ 
for any $t>0$ and $x,y\in\Omega$, and that
\begin{equation}\label{heatk-int}
\int_\Omega p_{{N}}^\Omega(t,x,y)\;dy\ =\ 1,
\qquad x\in\Omega,\ t>0,
\end{equation}
which follows from noticing that, for any $u\in C^\infty_0(\Omega)$,
\[
\partial_t\int_\Omega e^{t
{\Delta_{N,\Omega}}
}u=\int_\Omega\partial_t e^{t
{\Delta_{N,\Omega}}
}u=
\int_\Omega\Delta e^{t
{\Delta_{N,\Omega}}
}u=-\int_{\partial\Omega}\partial_\nu e^{t
{\Delta_{N,\Omega}}
}u=0
\]
and therefore
\[
\int_\Omega u(x)\;dx=\int_\Omega e^{t
{\Delta_{N,\Omega}}
}u(x)\;dx
=\int_\Omega\int_\Omega p_{{N}}^\Omega(t,x,y)\,u(y)\;dy\;dx
=\int_\Omega u(y) \int_\Omega p_{N}^\Omega(t,x,y)\;dx\;dy.
\]
By \eqref{heatk-int}, for any $x\in\Omega$ and $t>0$,
\begin{align*}
& u(x)-e^{t
{\Delta_{N,\Omega}}
}u(x)=\int_\Omega p_{{N}}^\Omega(t,x,y)
\left(u(x)-u(y)\right)dy
\end{align*}
and, exchanging the order of integration 
in~\eqref{boch-neu} (see below for a justification of this passage), one gets 
\begin{align*}
& {(-\Delta)^s_{N,\Omega}} u(x) =
\frac{s}{\Gamma(1-s)}\int_0^\infty\frac{u(x)-e^{t
{\Delta_{N,\Omega}}
}
u(x)}{t^{1+s}}\;dy
= \frac{s}{\Gamma(1-s)}\int_0^\infty\frac{\int_\Omega p_{{N}}^\Omega(t,x,y)
\left(u(x)-u(y)\right)dy}{t^{1+s}}\;dt \\
&= \frac{s}{\Gamma(1-s)}\;\PV\int_\Omega\left(u(x)-u(y)\right)
\int_0^\infty\frac{p_{{N}}^\Omega(t,x,y)}{t^{1+s}}\;dt\;dy=\PV\int_\Omega\frac{\left(u(x)-u(y)\right)k(x,y)}{{|x-y|}^{n+2s}}\;dy
,\end{align*}
where, in view of~\eqref{heatk-est}, we have
\begin{align*}
k(x,y) :=\frac{s}{\Gamma(1-s)}\;|x-y|^{n+2s}
\int_0^\infty\frac{p_{{N}}^\Omega(t,x,y)}{t^{1+s}}\;dt \simeq |x-y|^{n+2s}
\int_0^\infty\frac{e^{-|x-y|^2/t}}{t^{n/2+1+s}}\,dt
\simeq\int_0^\infty\frac{e^{-1/t}}{t^{n/2+1+s}}\,dt
\simeq 1.
\end{align*}
These considerations establish~\eqref{XAB}.
Note however
that in the above computations there is a limit exiting the integral in the $t$ variable,
namely:
\begin{equation}\label{0wodfkfk3455fkf}
\int_0^\infty\frac{\int_\Omega p_{{N}}^\Omega(t,x,y)
\left(u(x)-u(y)\right)dy}{t^{1+s}}\;dt=
\lim_{\eps\searrow 0}\int_0^\infty\frac{\int_{\Omega\setminus B_\eps(x)} p_{{N}}^\Omega(t,x,y)
\left(u(x)-u(y)\right)dy}{t^{1+s}}\;dt.
\end{equation}
To properly justify this 
we are going to build an integrable majorant in $t$ and independent of $\eps$ of
the integrand
\begin{equation}\label{INTE:G}
\frac{\int_{\Omega\setminus B_\eps(x)} p_{{N}}^\Omega(t,x,y)
\left(u(x)-u(y)\right)dy}{t^{1+s}}.
\end{equation}
To this end,
first of all we observe that,
by the boundedness of $u$ and \eqref{heatk-int},
\[
\left|\frac{\int_{\Omega\setminus B_\eps(x)} p_{{N}}^\Omega(t,x,y)
\left(u(x)-u(y)\right)dy}{t^{1+s}}\right|\leq 
\frac{2\|u\|_{L^\infty(\Omega)}}{t^{1+s}}
\int_{\Omega\setminus B_\eps(x)} p_{{N}}^\Omega(t,x,y)\;dy\leq 
\frac{2\|u\|_{L^\infty(\Omega)}}{t^{1+s}}
\]
which is integrable at infinity.
So, to obtain an integrable bound for~\eqref{INTE:G},
we can now focus on small values of $t$, say $t\in(0,1)$.
For this, we denote by~$p$ the heat kernel in~$\R^N$ and we write
\begin{multline*}
\int_{\Omega\setminus B_\eps(x)} p_{{N}}^\Omega(t,x,y)
\left(u(x)-u(y)\right)dy\ =\\
=\ \int_{\Omega\setminus B_\eps(x)} p(t,x,y)
\left(u(x)-u(y)\right)dy-
\int_{\Omega\setminus B_\eps(x)}(p_{{N}}^\Omega(t,x,y)-p(t,x,y))(u(x)-u(y))\;dy=:A+B.
\end{multline*}
We first manipulate $A$.
We reformulate $u$ as
\[
u(y)=u(x)+\nabla u(x)\cdot(y-x)+\eta(y)|x-y|^2,
\qquad y\in \R^n,\ \|\eta\|_{L^\infty(\R^n)}\leq \|u\|_{C^2(\Omega)},
\]
so that
\begin{align}\label{ALI}
\begin{split}
& \int_{\Omega\setminus B_\eps(x)} p(t,x,y)\left(u(x)-u(y)\right)dy
=\int_{\R^n\setminus B_\eps(x)} p(t,x,y)\left(u(x)-u(y)\right)dy
-u(x)\int_{\R^n\setminus \Omega} p(t,x,y)\;dy \\
& =\int_{\R^n\setminus B_\eps(x)} p(t,x,y)\nabla u(x)\cdot(x-y)\;dy
-\int_{\R^n\setminus B_\eps(x)} p(t,x,y)\,\eta(y)|x-y|^2\;dy-u(x)\int_{\R^n\setminus \Omega} p(t,x,y)\;dy.
\end{split}\end{align}
In the last expression, the first integral on the right-hand side is $0$
by odd symmetry, while for the second one
\begin{align}\label{902309}
\begin{split}
& \left|\int_{\R^n\setminus B_\eps(x)} p(t,x,y)
\,\eta(y)|x-y|^2\;dy\right|
\leq \|u\|_{C^2(\Omega)}\int_{\R^n\setminus B_\eps(x)} p(t,x,y)|x-y|^2\;dy \\
& \leq \const\|u\|_{C^2(\Omega)}t^{-/2}\int_{\R^n\setminus B_\eps(x)} e^{-|x-y|^2/(4t)}|x-y|^2\;dy
\leq \const\|u\|_{C^2(\Omega)}t\int_{\R^n\setminus B_{\eps/\sqrt{4t}}} e^{-|z|^2}|z|^2\;dz \\
& \leq \const\|u\|_{C^2(\Omega)}t.
\end{split}
\end{align}
As for
the last integral in~\eqref{ALI}, we have that
\begin{align}\label{6427890}
\begin{split}
& |u(x)|\int_{\R^n\setminus\Omega}p(t,x,y)\;dy\leq\const|u(x)|t^{-n/2}
\int_{\R^n\setminus\Omega} e^{-|x-y|^2/(4t)}\;dy\leq \\
& \leq\const|u(x)|t^{-n/2}
\int_{\R^n\setminus B_{\text{dist}(x,\partial\Omega)}} e^{-|y|^2/(4t)}\;dy
\leq\const|u(x)|
\int_{\R^n\setminus B_{\text{dist}(x,\partial\Omega)/\sqrt{4t}}} e^{-|z|^2}\;dz \\
& \leq \const|u(x)|e^{-\text{dist}(x,\partial\Omega)/\sqrt{4t}}.
\end{split}
\end{align}
Equations \eqref{902309} and \eqref{6427890} imply that 
\[
\frac{|A|}{t^{1+s}}\leq\const t^{-s},\qquad t\in(0,1),
\]
which is integrable for~$t\in(0,1)$.

We turn now to the estimation of $B$ which we rewrite as
\[
B=\int_{\Omega}\big(p_{{N}}^\Omega(t,x,y)-p(t,x,y)\big)\big(u(x)-u(y)\big)\chi_{\Omega\setminus B_\eps(x)}(y)\;dy
\]
where $\chi_U$ stands for the characteristic function of a set $U\subset\R^n$.
By definition, $B$ solves the heat equation in $\Omega$ with zero initial condition. Moreover,
since $u$ is supported in a compact subset $K$ of $\Omega$, $B$ is satisfying the 
(lateral) boundary condition
\begin{align*}
& \Big|B\big|_{\partial B}\Big|\leq \int_{\Omega}\big|p_{{N}}^\Omega(t,x,y)-p(t,x,y)\big||u(y)|\chi_{\Omega\setminus B_\eps(x)}(y)\;dy\leq \const t^{-n/2}\int_K e^{-c_1|x-y|^2/t}|u(y)|\;dy \\
& \leq \const\|u\|_{L^1(\Omega)}t^{-n/2} e^{-c_2/t}
\end{align*}
for some $c_1,c_2>0$, in view of \eqref{heatk-est} and that, for $x\in\partial \Omega$ and $y\in K$,
$|x-y|\geq\text{dist}(K,\partial\Omega)>0$. 
Then, by the parabolic maximum principle (see, for example, Section 7.1.4 in~\cite{MR1625845}),
\[
\frac{|B|}{t^{1+s}}\leq \const t^{-n/2-1-s} e^{-c_2/t},
\]
which again is integrable for~$t\in(0,1)$.
These observations
provide an integrable bound
for the integrand in~\eqref{ALI},
thus completing the justification of the claim in~\eqref{0wodfkfk3455fkf}, as desired.

\section{Proof of \eqref{TORUS}}\label{TORUS:A}

If~$u$ is periodic, we can write it in Fourier series as
$$ u(x)=\sum_{k\in\Z^n} u_k\,e^{2\pi ik\cdot x},$$
and the Fourier basis is also a basis of eigenfunctions. We have that
\begin{eqnarray*}
&&\int_{\R^n} \frac{u(x+y)+u(x-y)-2u(x)}{|y|^{n+2s}}\,dy=
\sum_{k\in\Z^n}
\int_{\R^n} \frac{u_k e^{ 2\pi ik\cdot(x+y)}+
u_ke^{2\pi ik\cdot(x-y)}-2u_ke^{2\pi ik\cdot x}}{|y|^{n+2s}}\,dy
\\&&\qquad
=\sum_{k\in\Z^n} u_k e^{ 2\pi ik\cdot x}
\int_{\R^n} \frac{e^{ 2\pi ik\cdot y}+
e^{-2\pi ik\cdot y}-2}{|y|^{n+2s}}\,dy
=\sum_{k\in\Z^n} u_k e^{ 2\pi ik\cdot x}\,|k|^{2s}
\int_{\R^n} \frac{e^{ 2\pi i\frac{k}{|k|}\cdot Y}+
e^{-2\pi i\frac{k}{|k|}\cdot Y}-2}{|Y|^{n+2s}}\,dY\\
\\&&\qquad
=\sum_{k\in\Z^n} u_k e^{ 2\pi ik\cdot x}\,|k|^{2s}
\int_{\R^n} \frac{e^{ 2\pi iY_1}+
e^{-2\pi iY_1}-2}{|Y|^{n+2s}}\,dY
=\const\sum_{k\in\Z^n} u_k e^{ 2\pi ik\cdot x}\,|k|^{2s}
\end{eqnarray*}
and this shows~\eqref{TORUS}.

\section{Proof of~\eqref{BARU:NO}}\label{BARU:NOA}

We fix~$\bar k\in\N$.
We consider the~$\bar k$th eigenvalue~$\lambda_{\bar k}>0$
and the corresponding normalized eigenfunction~$\phi_{\bar k}=:\bar u$.
We argue by contradiction and suppose that for any~$\e>0$ we can find~$v_\e$ such that~$
\|\bar u-v_\e\|_{C^2(B_1)}\le\e$, with~$(-\Delta)^s_{D,\Omega} v_\e=0$ in~$B_1$.

Using the notation in~\eqref{EXPA}, we have that~$\bar u_k=\delta_{k\bar k}$ and therefore
\begin{equation}\label{CO:ba}
\begin{split}
& \left\| (-\Delta)^s_{D,\Omega} \bar u-(-\Delta)^s_{D,\Omega} v_\e\right\|^2_{L^2(\Omega)}=
\left\| (-\Delta)^s_{D,\Omega} \bar u\right\|^2_{L^2(\Omega)}
=\left\| (-\Delta)^s_{D,\Omega}\phi_{\bar k}
\right\|^2_{L^2(\Omega)}\\
&\qquad=\left\| \lambda_{\bar k}^s\,\phi_{\bar k}
\right\|^2_{L^2(\Omega)}=\lambda_{\bar k}^{2s}.
\end{split}\end{equation}
Furthermore
\begin{eqnarray*}
&& \left\| (-\Delta)^s_{D,\Omega} \bar u-(-\Delta)^s_{D,\Omega} v_\e\right\|^2_{L^2(\Omega)}
=\left\| (-\Delta)^s_{D,\Omega} (\bar u- v_\e)\right\|^2_{L^2(\Omega)}
=\left\| 
\sum_{k=0}^{+\infty}\lambda_k^s (\bar u-v_\e)_k\, \phi_k
\right\|^2_{L^2(\Omega)}\\
&&\qquad\qquad= 
\sum_{k=0}^{+\infty}\lambda_k^{2s} (\bar u-v_\e)_k^2\le
\const\sum_{k=0}^{+\infty}\lambda_k^{2} (\bar u-v_\e)_k^2
=\const\left\| \Delta(\bar u- v_\e)\right\|^2_{L^2(\Omega)}\\
&&\qquad\qquad\le\const \| \bar u- v_\e\|^2_{C^2(\Omega)}\le\const\e.
\end{eqnarray*}
Comparing this with~\eqref{CO:ba}, we obtain that~$\lambda_{\bar k}^{2s}\le\const\e$,
which is a contradiction for small~$\e$. Hence, the proof of~\eqref{BARU:NO} is complete.

\section{Proof of~\eqref{LAP:L}}\label{LAP:LA}

Let
$$ v(t):=\int_0^t \frac{\dot u (\tau) }{(t-\tau)^{s}} \, d\tau.$$
Then, by~\eqref{L7849576AP:L}
and writing~$\vartheta:=\omega\,({t-\tau})$, we see that
\begin{eqnarray*}
{\mathcal{L}} v(\omega)&=&
\int_{0}^{+\infty } \left[
\int_0^t \frac{\dot u (\tau) }{(t-\tau)^{s}} \, d\tau\right]
e^{-\omega t}\,dt
=\int_0^{+\infty} \left[
\int_\tau^{+\infty}
\frac{\dot u (\tau) e^{-\omega t}}{(t-\tau)^{s}}
\,dt\right]\,d\tau\\
&=&\omega^{s-1} \int_0^{+\infty} \left[
\int_0^{+\infty}
\frac{\dot u (\tau) e^{-\omega \tau}\,e^{-\vartheta}}{\vartheta^{s}}
\,d\vartheta\right]\,d\tau
=\Gamma(1-s)\,
\omega^{s-1} \int_0^{+\infty} 
\dot u (\tau) e^{-\omega \tau}\,d\tau\\
&=& \Gamma(1-s)\,
\omega^{s-1} \int_0^{+\infty} \left(\frac{d}{d\tau}\big(
u (\tau) e^{-\omega \tau}\big)
+\omega u (\tau) e^{-\omega \tau}\right)
\,d\tau\\
&=& \Gamma(1-s)\,
\omega^{s-1} 
\left( -u (0)+ \omega
\int_0^{+\infty} 
u (\tau) e^{-\omega \tau}
\,d\tau \right)
\\ &=& \Gamma(1-s)\,
\omega^{s-1}
\left( -u (0)+ \omega{\mathcal{L}}u(\omega)\right),
\end{eqnarray*}
where~$\Gamma$ denotes here the Euler's Gamma Function.
This and~\eqref{0we0wy83fguguqgferighierhg0011} 
give~\eqref{LAP:L}, up to neglecting normalizing constants,
as desired.

It is also worth pointing out that, as~$s\nearrow1$,
formula~\eqref{LAP:L} recovers the classical derivative, since, by~\eqref{L7849576AP:L},
\begin{eqnarray*}
{\mathcal{L}}\dot u(\omega)&=&
\int_{0}^{+\infty } \dot u(t)e^{-\omega t}\,dt\\
&=&   
\int_{0}^{+\infty } \left(\frac{d}{dt} \big( u(t)e^{-\omega t}\big)
+\omega u(t)e^{-\omega t}
\right)\,dt\\
&=& -u(0)+
\omega\int_{0}^{+\infty } 
u(t)e^{-\omega t}\,dt\\
&=& -u(0)+\omega\,{\mathcal{L}}u(\omega),
\end{eqnarray*}
which is~\eqref{LAP:L} when~$s=1$.

\section{Proof of~\eqref{LAP:L2}}\label{LAP:LA2}

First, we compute the Laplace Transform 
of the constant function. Namely,
by~\eqref{L7849576AP:L}, for any~$b\in\R$,
\begin{equation}\label{09wieufrye8fgiro3tyghifdhg0101}
{\mathcal{L}} b(\omega)=
b\,\int_{0}^{+\infty } e^{-\omega t}\,dt=\frac{b}{\omega}.
\end{equation}
We also set
$$ \Psi(t):=
\int_0^t \frac{f(\tau)}{(t-\tau)^{1-s}}\,d\tau$$
and we use~\eqref{L7849576AP:L}
and the substitution~$\vartheta:=\omega\,({t-\tau})$
to calculate that
\begin{eqnarray*}
{\mathcal{L}} \Psi(\omega)&=&
\int_{0}^{+\infty } \left[
\int_0^t \frac{f(\tau)}{(t-\tau)^{1-s}}\,d\tau
\right]\,e^{-\omega t}\,dt\\
&=&
\int_{0}^{+\infty } \left[
\int_\tau^{+\infty} \frac{f(\tau)
\,e^{-\omega t}
}{(t-\tau)^{1-s}}\,dt
\right]\,d\tau\\
&=& 
\omega^{-s}\int_{0}^{+\infty } \left[
\int_0^{+\infty} \frac{f(\tau)
\,e^{-\omega \tau}\,e^{-\vartheta}
}{\vartheta^{1-s}}\,d\vartheta
\right]\,d\tau\\
&=& \Gamma(s)\,
\omega^{-s}\,\int_{0}^{+\infty } f(\tau)
\,e^{-\omega \tau}\,d\tau =
\Gamma(s)\,
\omega^{-s}\,{\mathcal{L}} f(\omega),
\end{eqnarray*}
where~$\Gamma$ denotes here the Euler's Gamma Function.

Exploiting this and~\eqref{09wieufrye8fgiro3tyghifdhg0101},
and making use also of~\eqref{LAP:L}, we can write the expression~$\partial^s_{C,t} u=f$
in terms of the Laplace Transform as
\begin{eqnarray*}
&& \omega^s \Big( {\mathcal{L}} u(\omega)-
{\mathcal{L}} b(\omega)\Big)=
\omega^s
{\mathcal{L}} u(\omega)-\omega^{s-1} u(0)=
{\mathcal{L}} (\partial^s_{C,t} u)(\omega)\\
&&\qquad\qquad=
{\mathcal{L}} f(\omega) = \frac{\omega^s}{\Gamma(s)} {\mathcal{L}} \Psi(\omega),
\end{eqnarray*}
with~$b:=u(0)$. Hence, dividing by~$\omega^s$ and inverting the Laplace Transform,
we obtain that
$$ u(t)-b= \frac{1}{\Gamma(s)} \Psi(t),
$$
which is~\eqref{LAP:L2}.

\section{Proof of~\eqref{MSDES}}\label{MSDES:A}

We take~$G$ to be the fundamental solution of the operator ``identity minus Laplacian'',
namely
\begin{equation}\label{9203dioeu45849022} 
G-\Delta G =\delta_0 \quad{\mbox{ in }}\R^n,\end{equation}
being~$\delta_0$ the Dirac's Delta.
The study of this fundamental solution can be done
by Fourier Transform in the sense of distributions,
and this leads to an explicit
representation
in dimension~$1$ 
recalling~\eqref{209rewifhweot894tfhkfheo485749}; we give here a general argument,
valid in any dimension, based on the heat kernel
$$ g(x,t):= {\frac {1}{(4\pi t)^{n/2}}}e^{-\frac{|x|^{2}}{4t}}.$$
Notice that~$\partial_t g=\Delta g$ and~$g(\cdot,0)=\delta_0(\cdot)$.
Let also
\begin{equation}\label{9203dioeu4584902} G(x):=\int_0^{+\infty} e^{-t} g(x,t)\,dt.\end{equation}
Notice that
\begin{eqnarray*}&& \Delta G(x)=
\int_0^{+\infty} e^{-t} \Delta g(x,t)\,dt
=\int_0^{+\infty} e^{-t} \partial_t g(x,t)\,dt
=\int_0^{+\infty} \Big(\partial_t(e^{-t} g(x,t))+e^{-t} g(x,t)\Big)\,dt
\\&&\qquad= -\delta_0(x)+
\int_0^{+\infty} e^{-t} g(x,t)\,dt
= -\delta_0(x)+G(x),\end{eqnarray*}
hence $G$, as defined in~\eqref{9203dioeu4584902}
solves~\eqref{9203dioeu45849022}.

Notice also that~$G$ is positive and bounded, due to~\eqref{9203dioeu4584902}.
We also claim that
\begin{equation}\label{FUBOU0100}
{\mbox{for any~$x\in\R^n\setminus B_1$, it holds that $G(x)\le C e^{-c|x|}$,}}
\end{equation}
for some~$c$, $C>0$. To this end, let us fix~$x\in\R^n\setminus B_1$
and distinguish two regimes. If~$t\in [0,\,|x|]$, we have that~$\frac{|x|^{2}}{t}\ge
|x|$ and thus
$$ g(x,t)\le
{\frac {1}{(4\pi t)^{n/2}}}e^{-\frac{|x|^{2}}{8t}}e^{-\frac{|x|}8}.$$
Consequently, using the substitution~$\rho:=\frac{|x|^{2}}{8t}$,
\begin{equation}\label{09847697065djr95607-302}
\int_0^{|x|} e^{-t} g(x,t)\,dt\le
\int_0^{|x|} {\frac {1}{(4\pi t)^{n/2}}}e^{-\frac{|x|^{2}}{8t}}e^{-\frac{|x|}8}\,dt=
\int_{|x|/8}^{+\infty}\frac{C\rho^{n/2}}{|x|^n}\,e^{-\rho}
\,e^{-\frac{|x|}8}\,\frac{|x|^2\,d\rho}{\rho^2}\le C|x|\,e^{-\frac{|x|}8},
\end{equation}
for some~$C>0$ possibly varying from line to line.
Furthermore
$$ \int_{|x|}^{+\infty} e^{-t} g(x,t)\,dt\le
\int_{|x|}^{+\infty} e^{-\frac{|x|}2}\,e^{-\frac{t}2} g(x,t)\,dt\le
e^{-\frac{|x|}2}\int_1^{+\infty} {\frac {e^{-\frac{t}2}}{(4\pi t)^{n/2}}}\,dt\le
C\,e^{-\frac{|x|}2},$$
for some~$C>0$. This and~\eqref{09847697065djr95607-302} give that
$$ \int_0^{+\infty} e^{-t} g(x,t)\,dt\le C|x|\,e^{-\frac{|x|}8},$$
up to renaming~$C$, which implies~\eqref{FUBOU0100}
in view of~\eqref{9203dioeu4584902}.

Now we compute the Laplace Transform of~$t^s$: namely, by~\eqref{L7849576AP:L},
\begin{equation}\label{78:-09er8}
{\mathcal{L}}(t^s)(\omega)=
\int_0^{+\infty} t^s e^{-\omega t}\,dt=\omega^{-1-s}
\int_0^{+\infty} \tau^s e^{-\tau}\,d\tau=C\omega^{-1-s}.
\end{equation}
We compare this result with the Laplace Transform of 
the mean squared displacement
related to the diffusion operator in~\eqref{MSDES}.
For this, we take~$u$ to be as in~\eqref{MSDES}
and, in the light of~\eqref{MSD:MSD:1},
we consider the function
\begin{equation}\label{9092rfj289f7289e1301-0} v(\omega):={\mathcal{L}}\left(
\int_{\R^n} |x|^2\, u(x,t)\,dx
\right)(\omega)=\int_{\R^n} |x|^2\, {\mathcal{L}} u(x,\omega)\,dx
.\end{equation}
In addition, by taking the Laplace Transform (in the variable~$t$, for a fixed~$x\in\R^n$) of
the equation in~\eqref{MSDES},
making use of~\eqref{LAP:L} we find that
\begin{equation} \label{OMElap}
\omega^s
{\mathcal{L}} u(x,\omega)-\omega^{s-1} \delta_0(x)=\Delta {\mathcal{L}} u(x,\omega).\end{equation}
Now, we let
\begin{equation} \label{4646}
W(x,\omega):= \omega^{1-\frac{sn}2} \,{\mathcal{L}} u(\omega^{-s/2}x,\omega).\end{equation}
{F}rom~\eqref{OMElap}, we have that
\begin{eqnarray*} \Delta W(x,\omega)&=&
\omega^{1-\frac{sn}2} \,\omega^{-s}\,\Delta{\mathcal{L}} u(\omega^{-s/2}x,\omega)
\\&=&\omega^{1-\frac{sn}2} \,\omega^{-s}\,\Big(
\omega^s
{\mathcal{L}} u(\omega^{-s/2} x,\omega)-\omega^{s-1} \delta_0(\omega^{-s/2} x)
\Big)\\
&=& W(x,\omega)-\omega^{-\frac{sn}2}\delta_0(\omega^{-s/2} x)\\
&=& W(x,\omega)-\delta_0(x),
\end{eqnarray*}
and so, comparing with~\eqref{9203dioeu45849022}, we have that~$W(x,\omega)=G(x)$.

Accordingly, by~\eqref{4646},
$$ 
{\mathcal{L}} u(x,\omega)=
\omega^{\frac{sn}2-1} \,W(\omega^{s/2}x,\,\omega)= \omega^{\frac{sn}2-1} \,G(\omega^{s/2}x).
$$
We insert this information into~\eqref{9092rfj289f7289e1301-0}
and we conclude that
$$ v(\omega)=\omega^{\frac{sn}2-1} \,
\int_{\R^n} |x|^2\, G(\omega^{s/2}x)\,dx=
\omega^{-1-s} \,
\int_{\R^n} |y|^2\, G(y)\,dy.
$$
We remark that the latter integral is finite, thanks to~\eqref{FUBOU0100},
hence we can write that
$$ v(\omega)=C\omega^{-1-s},$$
for some~$C>0$.

Therefore, we can compare this result with~\eqref{78:-09er8}
and use the inverse Laplace Transform to obtain that
the mean squared displacement in this case is proportional to~$t^s$, as desired.

\section{Memory effects of Caputo type}\label{MEMORY}

It is interesting to observe that the Caputo derivative
models a simple memory effect that the classical derivative cannot comprise.
For instance, integrating a classical derivative of a function~$u$
with~$u(0)=0$, one obtains
the original function ``independently on the past'', namely if we set
\begin{equation}\label{TIME1}
M_u(t):=\int_0^t \dot u(\vartheta)\,d\vartheta ,
\end{equation}
we just obtain in this case that~$M_u(t)=u(t)-u(0)=u(t)$.
On the other hand, an expression as in~\eqref{TIME1}
which takes into account the Caputo derivative does ``remember the past''
and takes into account the preceding events in such a way that
recent events ``weight'' more than far away ones.
To see this phenomenon, we can modify~\eqref{TIME1}
by defining, for every~$s\in(0,1)$,
\begin{equation}\label{TIME2} M_u^s(t):=
\int_0^t \partial^s_{C,t} u(\vartheta)\,d\vartheta .
\end{equation}
To detect the memory effect, for the sake of concreteness,
we take a large time~$t:=N\in\N$
and we suppose that the function~$u$ is constant on unit intervals,
that is~$u=u_k$ in~$[k-1,k)$, for each~$k\in\{1,\dots,N\}$,
with~$u_k\in\R$, and~$u(0)=u_1=0$. We see that~$M_u^s$ in this case does
not produce just the final outcome~$u_N$, but a weighted
average of the form
\begin{equation}\label{TIME3}
M_u^s(N)=\sum_{k=0}^{N-1} c_k\, u_{N-k},
\qquad{\mbox{with $c_j>0$ decreasing and }}c_j\simeq\frac1{j^s}
{\mbox{ for large $j$.}}
\end{equation}
To check this, we notice that, for all~$\tau\in(0,N)$,
$$\dot u(\tau)=\sum_{k=2}^N (u_{k}-u_{k-1})\delta_{k-1}(\tau),$$
and hence we exploit~\eqref{0we0wy83fguguqgferighierhg0011} and~\eqref{TIME2}
to see that
\begin{eqnarray*}
M_u^s(N)&=& \int_0^N \left[
\int_0^\vartheta \frac{\dot u (\tau) }{(\vartheta-\tau)^{s}} \, d\tau\right]
\,d\vartheta 
\\ &=&
\sum_{k=2}^N 
\int_0^N\left[
\int_0^\vartheta 
(u_{k}-u_{k-1})\delta_{k-1}(\tau)
\frac{d\tau}{(\vartheta-\tau)^{s}}\right]
\,d\vartheta\\
&=&
\sum_{k=2}^N
\int_{k-1}^N
\frac{(u_{k}-u_{k-1})}{(\vartheta-k+1)^{s}}
\,d\vartheta
\\
&=&
\sum_{k=2}^N u_{k}
\int_{k-1}^N
\frac{d\vartheta}{(\vartheta-k+1)^{s}}
-
\sum_{k=2}^N u_{k-1}
\int_{k-1}^N
\frac{d\vartheta}{(\vartheta-k+1)^{s}}
\\&=&
\sum_{k=2}^N u_{k}
\int_{k-1}^N
\frac{d\vartheta}{(\vartheta-k+1)^{s}}
-
\sum_{k=1}^{N-1} u_{k}
\int_k^N
\frac{d\vartheta}{(\vartheta-k)^{s}}\\
&=&
\sum_{k=1}^N u_{k} 
\left[
\int_{k-1}^N
\frac{d\vartheta}{(\vartheta-k+1)^{s}}
-
\int_k^N
\frac{d\vartheta}{(\vartheta-k)^{s}}
\right]\\
&=&
\sum_{k=1}^N u_{k} \,\frac{(N-k+1)^{1-s}-(N-k)^{1-s}}{1-s}\\
&=& \sum_{k=2}^N c_{N-k} \,u_{k}\\
&=&
\sum_{k=0}^{N-2} c_{k} \,u_{N-k} 
,\end{eqnarray*}
with
\begin{eqnarray*} c_{j}& :=&\frac{(j+1)^{1-s}-j^{1-s}}{1-s}
\,.\end{eqnarray*}
This completes the proof of the memory effect claimed in~\eqref{TIME3}.

\section{Proof of \eqref{DIV:EQ}}\label{DIV:FORM}

Since~$M$ is bounded and positive and $u$ is bounded, it holds that
\begin{equation}\label{FUORI}
\int_{\R^n\setminus B_1} \frac{| u(x)-u(x-y)| }{|M(x-y,y)\, y|^{n+2s}}\,dy\le\const
\int_{\R^n\setminus B_1} \frac{dy }{|y|^{n+2s}}\,dy\le \frac{\const}{s}.
\end{equation}
Moreover, for~$y\in B_1$,
\begin{equation}\label{T10} u(x-y)=u(x)-\nabla u(x)\cdot y +\frac12 D^2 u(x)\,y\cdot y+O(|y|^3).\end{equation}
To simplify the notation,
we now fix~$x\in\R^n$ and we define~${\mathcal{M}}(y):= M(x-y,y)$.
Then, for $y\in B_1$, we have that
\[
M(x-y,y)\, y = {\mathcal{M}}(y)\,y=
{\mathcal{M}}(0)\,y + \sum_{i=1}^n \partial_i {\mathcal{M}}(0) \, y\,y_i+O(|y|^3)\]
and so
\begin{eqnarray*} |M(x-y,y)\, y|^2 &=& 
|{\mathcal{M}}(0)\,y|^2 + 
2 \sum_{i=1}^n ({\mathcal{M}}(0)\,y)\cdot(\partial_i {\mathcal{M}}(0) \, y)\,y_i
+O(|y|^4).\end{eqnarray*}
Consequently, since~${\mathcal{M}}(0)=M(x,0)$ is non-degenerate,
we can write
$$ {\mathcal{E}}(y):=
2 \sum_{i=1}^n ({\mathcal{M}}(0)\,y)\cdot(\partial_i {\mathcal{M}}(0) \, y)\,y_i
=
O(|y|^3)$$
and
\begin{equation}\label{T9}
\begin{split}
& |M(x-y,y)\, y|^{-n-2s} \\
=\;&
\left(
|{\mathcal{M}}(0)\,y|^2 + 
{\mathcal{E}}(y)
+O(|y|^4) \right)^{-\frac{n+2s}2}
\\ =\;&
|{\mathcal{M}}(0)\,y|^{-n-2s} \left(
1 + 
|{\mathcal{M}}(0)\,y|^{-2}\,{\mathcal{E}}(y)
+O(|y|^2) \right)^{-\frac{n+2s}2}
\\ =\;&
|{\mathcal{M}}(0)\,y|^{-n-2s} \left(
1 -\frac{n+2s}2\,
|{\mathcal{M}}(0)\,y|^{-2}\,{\mathcal{E}}(y)
+O(|y|^2) \right)
\\ =\;&
|{\mathcal{M}}(0)\,y|^{-n-2s} -\frac{n+2s}2\,
|{\mathcal{M}}(0)\,y|^{-n-2s-2}\,{\mathcal{E}}(y)
+O(|y|^{2-n-2s}) .
\end{split}\end{equation}
Hence (for smooth and bounded functions~$u$, and~$y\in B_1$) we obtain that
\begin{eqnarray*}
&& \frac{
u(x)-u(x-y)
}{|M(x-y,y)\,y|^{n+2s}}\\
&=& 
\frac{
u(x)-u(x-y)
}{|{\mathcal{M}}(0)\,y|^{n+2s} }
-\frac{n+2s}2\,
\frac{
\big( u(x)-u(x-y)\big)\,{\mathcal{E}}(y)
}{ |{\mathcal{M}}(0)\,y|^{n+2s+2}}
+O(|y|^{3-n-2s}).
\end{eqnarray*}
Thus, since the map~$y\mapsto\frac{\nabla u(x)\cdot y}{|{\mathcal{M}}(0)\,y|^{n+2s} }$
is odd, recalling~\eqref{T10} we conclude that
\begin{equation}\label{thVER:01}
\begin{split}
&\int_{B_1} \frac{u(x)-u(x-y)}{|M(x-y,y)\,y|^{n+2s}}\,dy\\
=\;& \int_{B_1}\left( \frac{
u(x)-u(x-y)
}{|{\mathcal{M}}(0)\,y|^{n+2s} }
-\frac{n+2s}2\,
\frac{
\big( u(x)-u(x-y)\big)\,{\mathcal{E}}(y)
}{ |{\mathcal{M}}(0)\,y|^{n+2s+2}}
+O(|y|^{3-n-2s})\right)\,dy\\
=\;&
\int_{B_1}\left( \frac{
u(x)-u(x-y)-\nabla u(x)\cdot y
}{|{\mathcal{M}}(0)\,y|^{n+2s} }
-\frac{n+2s}2\,
\frac{
\big( u(x)-u(x-y)\big)\,{\mathcal{E}}(y)
}{ |{\mathcal{M}}(0)\,y|^{n+2s+2}}
+O(|y|^{3-n-2s})\right)\,dy
\\
=\;&
\int_{B_1}\left( -\frac{
D^2u(x) \,y\cdot y}{2\,|{\mathcal{M}}(0)\,y|^{n+2s} }
-(n+2s)\,\sum_{i=1}^n
\frac{ (\nabla u(x)\cdot y)
\,\big( ({\mathcal{M}}(0)\,y)\cdot(\partial_i {\mathcal{M}}(0) \, y)\big)\,y_i
}{ |{\mathcal{M}}(0)\,y|^{n+2s+2}}
+O(|y|^{3-n-2s})\right)\,dy
\end{split}\end{equation}
Now we observe that, for any~$\alpha\ge0$,
\begin{equation}\label{GEN:PHI}
\begin{split}
&{\mbox{if $\varphi$ is positively
homogeneous of degree $2+\alpha$ and~$T\in{\rm Mat}(n\times n)$, then}}
\\ & (1-s)\,\int_{B_1} \frac{\varphi(y)}{|Ty|^{n+2s+\alpha}}\,dy=
\frac12\,\int_{S^{n-1}} \frac{\varphi(\omega)}{|T\omega|^{n+2s+\alpha}}\,d{\mathcal{H}}^{n-1}_{\omega}.
\end{split}\end{equation}
Indeed, using polar coordinates and the fact that~$\varphi(\rho\omega)=\rho^{2+\alpha}
\varphi(\omega)$, for any~$\rho\ge0$ and~$\omega\in S^{n-1}$, thanks to the homogeneity,
we see that
\begin{eqnarray*}
&& \int_{B_1} \frac{\varphi(y)}{|Ty|^{n+2s+\alpha}}\,dy=
\iint_{(0,1)\times S^{n-1}} \frac{\rho^{n-1}\,\varphi(\rho\omega)}{\rho^{n+2s+\alpha}\,
|T\omega|^{n+2s+\alpha}}\,d\rho\,d{\mathcal{H}}^{n-1}_{\omega}\\
&&\qquad=
\iint_{(0,1)\times S^{n-1}} \frac{\rho^{1-2s}\varphi(\omega)}{
|T\omega|^{n+2s+\alpha}}\,d\rho\,d{\mathcal{H}}^{n-1}_{\omega}=
\frac{1}{2(1-s)}
\int_{S^{n-1}} \frac{\varphi(\omega)}{
|T\omega|^{n+2s+\alpha}}\,d{\mathcal{H}}^{n-1}_{\omega},
\end{eqnarray*}
which implies~\eqref{GEN:PHI}.

Using~\eqref{GEN:PHI} (with~$\alpha:=0$ and~$\alpha:=2$), we obtain that
\begin{eqnarray*}
&& \lim_{s\nearrow1} \,(1-s)\,
\int_{B_1} \frac{D^2u(x) \,y\cdot y}{|{\mathcal{M}}(0)\,y|^{n+2s} }\,dy=\frac12\,
\int_{S^{n-1}} \frac{D^2u(x) \,\omega\cdot \omega}{
|{\mathcal{M}}(0)\,\omega|^{n+2} }\,d{\mathcal{H}}^{n-1}_{\omega}\end{eqnarray*}
and
\begin{eqnarray*}
&&\lim_{s\nearrow1} \,(1-s)\,
\int_{B_1}
\frac{ (\nabla u(x)\cdot y)
\,\big( ({\mathcal{M}}(0)\,y)\cdot(\partial_i {\mathcal{M}}(0) \, y)\big)\,y_i
}{ |{\mathcal{M}}(0)\,y|^{n+2s+2}}\,dy\\ &&\qquad\qquad=\frac12\,
\int_{S^{n-1}}
\frac{ (\nabla u(x)\cdot \omega)
\,\big( ({\mathcal{M}}(0)\,\omega)\cdot(\partial_i {\mathcal{M}}(0) \, \omega)\big)\,\omega_i
}{ |{\mathcal{M}}(0)\,\omega|^{n+4}}
\,d{\mathcal{H}}^{n-1}_{\omega}
\end{eqnarray*}
Thanks to this, \eqref{FUORI} and~\eqref{thVER:01}, we find that
\begin{equation}\label{thVER:02}
\begin{split}
&\lim_{s\nearrow1} \,\int_{\R^n} \frac{u(x)-u(x-y)}{|M(x-y,y)\,y|^{n+2s}}\,dy\\
=\;&\lim_{s\nearrow1} \,\int_{B_1} \frac{u(x)-u(x-y)}{|M(x-y,y)\,y|^{n+2s}}\,dy\\
=\;&-\frac14\,
\int_{S^{n-1}} \frac{D^2u(x) \,\omega\cdot \omega}{
|{\mathcal{M}}(0)\,\omega|^{n+2} }\,d{\mathcal{H}}^{n-1}_{\omega}
-\frac{n+2}{2}\,\sum_{i=1}^n
\int_{S^{n-1}}
\frac{ (\nabla u(x)\cdot \omega)
\,\big( ({\mathcal{M}}(0)\,\omega)\cdot(\partial_i {\mathcal{M}}(0) \, \omega)\big)\,\omega_i
}{ |{\mathcal{M}}(0)\,\omega|^{n+4}}
\,d{\mathcal{H}}^{n-1}_{\omega}\\
=\;& -\sum_{i,j=1}^n a_{ij}(x) \partial^2_{ij} u(x)-\sum_{j=1}^n b_{j} \partial_j u(x),
\end{split}\end{equation}
with
\begin{eqnarray*}
&& a_{ij}(x):= \frac14\,\int_{S^{n-1}} \frac{\omega_i\, \omega_j}{
|{\mathcal{M}}(0)\,\omega|^{n+2} }\,d{\mathcal{H}}^{n-1}_{\omega}
=\frac14\,\int_{S^{n-1}} \frac{\omega_i\, \omega_j}{
|M(x,0)\,\omega|^{n+2} }\,d{\mathcal{H}}^{n-1}_{\omega}
\\{\mbox{and }}&&
b_j(x):=\frac{n+2}{2}\,\sum_{i=1}^n
\int_{S^{n-1}}
\frac{ \omega_i\,\omega_j
\,\big( ({\mathcal{M}}(0)\,\omega)\cdot(\partial_i {\mathcal{M}}(0) \, \omega)\big)
}{ |{\mathcal{M}}(0)\,\omega|^{n+4}}
\,d{\mathcal{H}}^{n-1}_{\omega}
.\end{eqnarray*}
We observe that
\begin{equation}\label{SP:001}
b_j=\sum_{i=1}^n \partial_i a_{ij}(x).
\end{equation}
To check this, we first compute that
\begin{equation}\label{7o9e32tu6o51ueh}
\begin{split}
\sum_{i=1}^n \partial_i a_{ij}(x) \;&=\frac14\,\sum_{i=1}^n \partial_{x_i}
\left( \int_{S^{n-1}} \frac{\omega_i\, \omega_j}{
|M(x,0)\,\omega|^{n+2} }\,d{\mathcal{H}}^{n-1}_{\omega}
\right)\\
&=-\frac{n+2}{4}\,\sum_{i=1}^n
\int_{S^{n-1}} \frac{\omega_i\, \omega_j\,\big( (M(x,0)\,\omega)\cdot(\partial_{x_i}
M(x,0)\,\omega)\big)}{
|M(x,0)\,\omega|^{n+4} }\,d{\mathcal{H}}^{n-1}_{\omega}
.\end{split}
\end{equation}
Now, we write a Taylor expansion of~$M(x,y)$ in the variable~$y$ of the form
$$ M_{\ell m}(x,y)=A_{\ell m}(x) + B_{\ell m}(x)\cdot y+O(y^2),$$
for some~$A_{\ell m}:\R^n\to\R$ and~$B_{\ell m}:\R^n\to\R^n$. We notice that
\begin{equation}\label{98ygikkdhweuidfqwuwuigd63etryeu}
\partial_{x_i} M_{\ell m}(x,0)=\partial_{x_i} A_{\ell m}(x).
\end{equation}
Also,
\begin{equation}\label{9weuf8247279fhidfh}
\begin{split}
\partial_{i} {\mathcal{M}}_{\ell m}(0)\,&=\lim_{y\to0} \partial_{y_i}\big( M_{\ell m}(x-y,y) \big)
\\ &=\lim_{y\to0} \partial_{y_i}\big( A_{\ell m}(x-y) + B_{\ell m}(x-y)\cdot y+O(y^2) \big)\\&=
-\partial_{x_i}A_{\ell m}(x) + B_{\ell m}(x)\cdot e_i.\end{split}
\end{equation}
Furthermore, we
use the structural assumption~\eqref{DIU8234j:IPOT}, and we see that
\begin{eqnarray*}
&& 
A_{\ell m}(x) - B_{\ell m}(x)\cdot y+O(y^2)
=M(x,-y)\\
&& \qquad=M(x-y,y)=
A_{\ell m}(x-y) + B_{\ell m}(x-y)\cdot y+O(y^2)
\\ &&\qquad=A_{\ell m}(x) -\nabla A_{\ell m}(x)\cdot y
+ B_{\ell m}(x)\cdot y+O(y^2).
\end{eqnarray*}
Comparing the linear terms, this gives that
$$ 2B_{\ell m}(x)=\nabla A_{\ell m}(x).$$
This and~\eqref{9weuf8247279fhidfh} imply that
\[ \partial_{i} {\mathcal{M}}_{\ell m}(0)=
-\partial_{x_i}A_{\ell m}(x) + \frac12\nabla A_{\ell m}(x)\cdot e_i
=-\frac12 \partial_{x_i}A_{\ell m}(x) 
.\]
Comparing this with~\eqref{98ygikkdhweuidfqwuwuigd63etryeu},
we see that
\[ \partial_{x_i} M_{\ell m}(x,0)=-2\partial_{i} {\mathcal{M}}_{\ell m}(0).
\]
So, we insert this information into~\eqref{7o9e32tu6o51ueh}
and we conclude that
\[
\sum_{i=1}^n \partial_i a_{ij}(x)
=\frac{n+2}{2}\,\sum_{i=1}^n
\int_{S^{n-1}} \frac{\omega_i\, \omega_j\,\big( (M(x,0)\,\omega)\cdot(\partial_{i}
{\mathcal{M}}(0)\,\omega)\big)}{
|M(x,0)\,\omega|^{n+4} }\,d{\mathcal{H}}^{n-1}_{\omega}.\]
This establishes~\eqref{SP:001}, as desired.

\begin{figure}
    \centering
    \includegraphics[width=13cm]{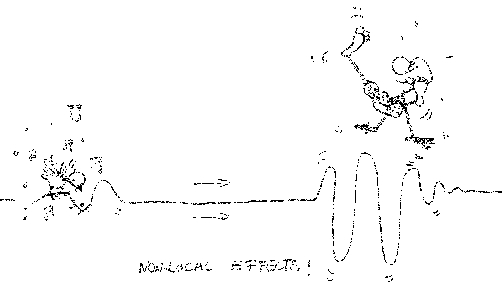}
    \caption{\it {{A nice representation of nonlocal effects.}}}
    \label{HIRSCH}
\end{figure}

Then, plugging~\eqref{SP:001} into~\eqref{thVER:02}, we obtain the equation in
divergence form\footnote{A slightly different approach as that in~\eqref{DIV:EQ}
is to consider the energy functional in~\eqref{E:2 en}
and prove, e.g. by Taylor expansion, that it converges to
the energy functional
$$ \const \int_{\R^n}a_{ij}(x) \,\partial_i u(x)\,\partial_j u(x)\,dx.$$
On the other hand, a different proof of~\eqref{DIV:EQ},
that was nicely pointed out to us by Jonas Hirsch (who has also acted as a skilled
cartoonist for Figure~\ref{HIRSCH})
after a lecture,
can be performed by taking into account the weak form of the operator in~\eqref{DIU8234j:3},
i.e. integrating such expression against a test function~$\varphi\in C^\infty_0(\R^n)$,
thus finding
\begin{eqnarray*}
&& (1-s)\,\iint_{\R^n\times\R^n} \frac{ \big(u(x)-u(x-y)\big)\,\varphi(x)}{
|M(x-y,y)\, y|^{n+2s}}\,dx\,dy
\\ &=&
(1-s)\,\iint_{\R^n\times\R^n} \frac{ \big(u(x)-u(z)\big)\,\varphi(x)}{
|M(z,x-z)\,(x-z)|^{n+2s}}\,dx\,dz
\\ &=&
(1-s)\,\iint_{\R^n\times\R^n} \frac{ \big(u(z)-u(x)\big)\,\varphi(z)}{
|M(x,z-x)\,(x-z)|^{n+2s}}\,dx\,dz
\\ &=&
-(1-s)\,\iint_{\R^n\times\R^n} \frac{ \big(u(x)-u(z)\big)\,\varphi(z)}{
|M(z,x-z)\,(x-z)|^{n+2s}}\,dx\,dz
,\end{eqnarray*}
where the structural condition~\eqref{DIU8234j:IPOT} has been used in the last line.
This means that the weak formulation of
the operator in~\eqref{DIU8234j:3} can be written as
$$ \frac{1-s}{2}\,\iint_{\R^n\times\R^n} \frac{ \big(u(x)-u(z)\big)\,\big(\varphi(x)-\varphi(z)\big)}{
|M(z,x-z)\,(x-z)|^{n+2s}}\,dx\,dz.$$
So one can expand this expression and take the limit as~$s\nearrow1$, to obtain
$$ \const \int_{\R^n}a_{ij}(x) \,\partial_i u(x)\,\partial_j\varphi(x)\,dx,$$
which is indeed the weak formulation of the classical divergence form operator.}
which was claimed in~\eqref{DIV:EQ}.

\section{Proof of \eqref{NONDIV:EQ}}\label{NONDIV:FORM}

First we observe that
\begin{equation}\label{PRI}
\int_{\R^n\setminus B_1} \frac{
|u(x)-u(x-y)|
}{|M(x,y)\,y|^{n+2s}}\,dy\le \const
\int_{\R^n\setminus B_1} \frac{
dy
}{|y|^{n+2s}}\le \frac{\const}{s}.
\end{equation}
Furthermore, for~$y\in B_1$,
\[ M(x,y)\,y = M(x,0)\,y+O(|y|^2).\]
Consequently,
\[ |M(x,y)\,y|^2=|M(x,0)\,y|^2+O(|y|^3)\]
and so, from the non-degeneracy of~$M(\cdot,\cdot)$,
\begin{eqnarray*}&& |M(x,y)\,y|^{-n-2s}
=\big( |M(x,0)\,y|^2+O(|y|^3)\big)^{-\frac{n+2s}2}
\\&&\qquad=|M(x,0)\,y|^{-n-2s}\big( 1+O(|y|)\big)^{-\frac{n+2s}2}=
|M(x,0)\,y|^{-n-2s}\big( 1-O(|y|)\big).
\end{eqnarray*}
Using this and the expansion in~\eqref{T10}, we see that,
for~$y\in B_1$,
\begin{eqnarray*}
&& \frac{
u(x)-u(x-y)-\nabla u(x)\cdot y
}{|M(x,y)\,y|^{n+2s}}\\
&=&|M(x,0)\,y|^{-n-2s}\big( 1-O(|y|)\big)\left( -\frac12 D^2 u(x)\,y\cdot y+O(|y|^3)\right)\\
&=& |M(x,0)\,y|^{-n-2s}\left( -\frac12 D^2 u(x)\,y\cdot y+O(|y|^3)\right).
\end{eqnarray*}
Thus,
since, in the light of~\eqref{DIU8234j:IPOT:2},
we know that the map~$y\mapsto\frac{
\nabla u(x)\cdot y
}{|M(x,y)\,y|^{n+2s}}$ is odd,
we can write that
\begin{eqnarray*}
\int_{B_1} \frac{
u(x)-u(x-y)
}{|M(x,y)\,y|^{n+2s}}\,dy&=&
\int_{B_1} \frac{
u(x)-u(x-y)-\nabla u(x)\cdot y
}{|M(x,y)\,y|^{n+2s}}\,dy\\
&=&-\frac12
\int_{B_1}\frac{D^2 u(x)\,y\cdot y}{ |M(x,0)\,y|^{n+2s} }\,dy
+\frac{O(1)}{3-2s}\\
&=&-\frac{\const}{1-s}
\int_{ S^{n-1}}\frac{
D^2 u(x)\,\omega\cdot\omega}{ |M(x,0)\,\omega|^{n+2s} }\,d{\mathcal{H}}^{n-1}_{\omega}
+\frac{O(1)}{3-2s}
,\end{eqnarray*}
where the last identity follows
by using~\eqref{GEN:PHI} (with~$\alpha:=0$).
{F}rom this and~\eqref{PRI} we obtain that
\begin{eqnarray*}
\lim_{s\nearrow1} (1-s)\,\int_{\R^n} \frac{
u(x)-u(x-y)
}{|M(x,y)\,y|^{n+2s}}\,dy &=&
\lim_{s\nearrow1} (1-s)\,\int_{B_1} \frac{
u(x)-u(x-y)
}{|M(x,y)\,y|^{n+2s}}\,dy \\ &=&
-\const
\int_{ S^{n-1}}\frac{
D^2 u(x)\,\omega\cdot\omega}{ |M(x,0)\,\omega|^{n+2} }\,d{\mathcal{H}}^{n-1}_{\omega}
\\&=&-\const\sum_{i,j=1}^n
\int_{ S^{n-1}}\frac{
\omega_i \omega_j}{ |M(x,0)\,\omega|^{n+2} }\,d{\mathcal{H}}^{n-1}_{\omega}\,\partial^2_{ij} u(x),
\end{eqnarray*}
which gives~\eqref{NONDIV:EQ}.

\end{appendix}

\vfill

\end{document}